\documentclass[10pt]{article}
\usepackage{mathrsfs}
\usepackage{amsmath,amssymb}
\usepackage{amsthm}
\usepackage[colorlinks,linkcolor=red,anchorcolor=blue,citecolor=green]{hyperref}
\usepackage{colonequals}
\usepackage{graphics}
\usepackage{graphicx}
\usepackage{color}
\usepackage{enumerate}
\usepackage[titletoc]{appendix}
\allowdisplaybreaks

\usepackage{tikz}
\usepackage{indentfirst}

\usepackage[margin=1in]{geometry}
\usepackage[T1]{fontenc}
\IfFileExists{lmodern.sty}{\usepackage{lmodern}}{}
\IfFileExists{microtype.sty}{
     }{}
\usepackage{amsmath,amssymb,amsthm,mathtools}
\usepackage{enumitem}
\usepackage{xcolor}
\usetikzlibrary{arrows.meta,calc,positioning}
\usepackage{float}
\usepackage{hyperref}

\newtheorem{theorem}{Theorem}[section]
\newtheorem{proposition}[theorem]{Proposition}
\newtheorem{lemma}[theorem]{Lemma}
\newtheorem{corollary}[theorem]{Corollary}
\newtheorem{conjecture}[theorem]{Conjecture}
\theoremstyle{definition}
\newtheorem{definition}[theorem]{Definition}
\newtheorem{example}[theorem]{Example}
\theoremstyle{remark}
\newtheorem{remark}[theorem]{Remark}

\renewcommand{\theequation}{\arabic{section}.\arabic{equation}}
\newcommand{\G}{\mathcal G}
\newcommand{\Hh}{\mathcal H}
\newcommand{\X}{\mathcal X}
\newcommand{\Y}{\mathcal Y}
\newcommand{\Pp}{\mathbb P}
\newcommand{\E}{\mathbb E}
\newcommand{\Ber}{\operatorname{Ber}}

\newcommand{\one}{\mathbf 1}
\newcommand{\cC}{\mathcal C}
\newcommand{\cE}{\mathcal E}

\newcommand{\Aut}{\operatorname{Aut}}

\newcommand{\Ninf}{N_{\infty}}

\begin{document}
\begin{center}
{{\large \bf Strict monotonicity of weighted inhomogeneous\\
percolation thresholds under covering maps }\footnote{The project is supported partially by CNNSF (No.\,12671180) and by
Hu Xiang Gao Ceng Ci Ren Cai Ju Jiao Gong Cheng-Chuang Xin Ren Cai (No.\,2019RS1057).}}
\end{center}
\begin{center}
Jing Liu$^a$, Kainan Xiang$^a$, Lang Zou$^b$\\
\footnotesize{$^a$Hunan Key Laboratory for Computation and Simulation in Science and Engineering \&\\
   Key Laboratory of Intelligent Computing and Information Processing of Ministry of Education \&\\
   School of Mathematics and Computational Science, Xiangtan University\\
   Xiangtan City 411105, Hunan Province, P. R. China\\
   $^b$School of Artificial Intelligence and Statistics, Xiangtan University\\
   Xiangtan City 411105, Hunan Province, P. R. China}\\
\footnotesize{Emails:  jingliu0813@126.com (Liu), \ \ kainan.xiang@xtu.edu.cn (Xiang)}\\
\footnotesize{\hskip 0mm langzou@xtu.edu.cn} (Zou)
\end{center}

\begin{abstract}
We study Bernoulli bond percolation under weak covering maps for a class of
weighted inhomogeneous product measures such that each edge $e$ carries a positive
weight $w(e)$ and is open at intensity $\lambda\geq 0$ with probability
$1-\exp(-\lambda w(e))$. We prove a strict monotonicity of the weighted critical intensities $\lambda_c$ in the covering compatibility setting under bounded degree, uniform ellipticity, tame fibre and disjoint tree-lifting.
While for the weighted uniform connection thresholds $\lambda_{\mathrm{conn}}$, we establish their strict monotonicity for non-injective weak covering maps under bounded degree, covering compatibility, uniform ellipticity, bounded-fibre
and disjoint tree-lifting. Moreover, we verify a general-profile theorem on strict monotonicities of the two thresholds $\lambda_c$ and $\lambda_{\mathrm{conn}}$ respectively under explicit compact-interval
finite-energy, derivative bounds and aforementioned geometric conditions.\\


\end{abstract}

\section{Introduction}
\setcounter{equation}{0}
\noindent
The strict monotonicity of thresholds (critical parameters) for independent Bernoulli percolations is of significant interests (\cite{Be-Sc1996, Gr1999, Ha-Jo2006, Ly-Pe2016}). This paper devotes to study the strict threshold monotonicity of weighted inhomogeneous Bernoulli bond percolations under covering maps.

Given any locally finite connected infinite graph $\X=(V(\X),E(\X))$ with vertex set $V(\X)$ and edge set $E(\X).$ The homogeneous Bernoulli bond (resp.\,\,site)
percolation $\X[p]$ (introduced by Broadbent and Hammersley (1957) \cite{Br-Ha1957}) declares every edge (resp.\,\,vertex) open independently with a common probability
$p\in [0,1]$. Denote by $p_c(\X):=p_c^{\rm bond}(\X)$ (resp.\,\,$p_c^{\rm site}(\X)$) the critical probability for existence of infinite open-edge (resp.\,\,open-vertex) clusters in the percolation.
Write $p_u(\X)$ for the critical probability of uniqueness of infinite open cluster in the bond percolation on $\X.$
\begin{enumerate}
\item[{\bf (i)}] By Hammersley (1961) \cite{Ham1961}, $p_c(\X)\leq p_c^{\rm site}(\X).$ Clearly $p_c(\X)\leq p_u(\X).$ Fundamental problems arise naturally: When does $p_c(\X)< p_c^{\rm site}(\X)?$  When does $p_c(\X)< p_u(\X)?$
\item[{\bf (ii)}] For any connected infinite subgraph $\Y$ of $\X,$ $p_c(\X)\leq p_c(\Y)$. When does $p_c(\X)< p_c(\Y)$?
\end{enumerate}
Note that $p_c(\X)=p_c^{\rm site}(\X)$ for trees $\X;$ and $p_c(\X)=p_u(\X)$ for all amenable quasi-transitive graphs $\X$ (\cite{Ga-Ke1992, Bu-Ke1989}), and
$\X[p]$ almost surely has a unique infinite open cluster for any quasi-transitive graph $\X$ and any $p>p_u(\X)$ (\cite{Hag-Pe1999, Ha-Pe1999, Sc1999, Ly-Sc1999}).

For the strict monotonicity $p_c(\X)< p_c^{\rm site}(\X),$ work in the 1980s in two dimensions was obtained in \cite{Hi1982, Ke1982, To1985}. Then Grimmett and Stacey (1998) \cite{Gr-St1998} proved it holds for a general class of graphs including $\mathbb{Z}^d$ with $d\geq 2$ by the enhancement technique of Aizenman and Grimmett (1991) \cite{Ai-Gr1991}. This technique was extended to abstract monotonic graph perturbations by van den Berg (1993) \cite{van1993} and has since become a powerful tool for establishing strict threshold monotonicity. Whereas for the strict monotonicity $p_c(\X)< p_u(\X),$ Benjamini and Schramm (1996) \cite{Be-Sc1996} proposed the following famous conjecture:

\begin{conjecture}[\mbox{\cite[Conjecture 6]{Be-Sc1996}}]\label{conj:pu}
For nonamenable quasi-transitive graph $\X$, $p_c(\X)<p_u(\X).$
\end{conjecture}
\noindent Hutchcroft (2020) \cite{Hu2020} verified a remarkable result that $p_c(\X)<p_u(\X)$ holds for any $\X$ with its automorphism group ${\rm Aut}(\X)$ having a quasi-transitive nonunimodular subgroup. For historical progresses on
Conjecture \ref{conj:pu}, see \cite{Hu2020} and references therein.

Let $p_{\rm ac}(\X)$ be the threshold of existence of absolutely continuous spectrums for infinite open clusters in percolations $\X[p]$. Clearly $p_c(\X)\leq p_{\rm ac}(\X).$
Note the following new strict monotonicity on $p_c(\X)$ and $p_{\rm ac}(\X)$ was proved recently by Becker and Oltman (2026) \cite{Be-Ol2026}: For any $d$-regular tree $\X$ with $d\geq 3,$ $p_c(\X)<p_{\rm ac}(\X).$ It is interesting to extend this strict monotonicity
to quasi-transitive graphs.
\smallskip

Now we focus on describing the strict monotonicity {\bf (ii)} and its extended version. For progresses of this strict monotonicity made in the 1980s and 1990s, see \cite[p.\,30, pp.\,57-71]{Gr1999}. Particularly,
 $\mathbb{Z}^d$ can be viewed naturally as a subgraph of $\mathbb{Z}^{d+1}$ for any $d\in\mathbb{N}:=\{1,2,3,\cdots\},$ and
$$p_c\left(\mathbb{Z}^{d+1}\right)<p_c\left(\mathbb{Z}^d\right).$$
Benjamini and Schramm (1996) \cite{Be-Sc1996} proved that $p_c(\X)\leq p_c(\Y)$ when $\Y$ is a quotient graph of $\X,$ and asked when the corresponding strict inequality holds. They believed that $p_c(\X)< p_c(\Y)$ holds for quasi-transitive graphs $\X$ and $\Y$ such that $\X$ covers but is not isomorphic to $\Y$, and $p_c(\X) < 1;$ which was confirmed by Martineau and Severo (2019) \cite[Theorem 0.1]{Ma-Se2019} (Indeed, they proved stronger conclusions).

 Martineau and Severo \cite{Ma-Se2019} proved that, under suitable geometric hypotheses, a
nontrivial covering operation strictly increases the critical probabilities $p_c$ and $p_u$:
\begin{enumerate}
\item[{\bf (a)}] Given any two locally finite connected infinite graphs $\X=(V(\X),E(\X))$ and  $\Y=(V(\Y),E(\Y))$. If $\pi:\ V(\X)\to V(\Y)$ is a weak covering map between bounded-degree graphs, the
fibres of $\pi$ are tame, and $\pi$ has the disjoint tree-lifting property, then
\[
  p_c(\X)<1 \ \Longrightarrow\ p_c(\X)<p_c(\Y).
\]
\item[{\bf (b)}] Assume $\X=(V(\X),E(\X))$ and $\Y=(V(\Y),E(\Y))$ are quasi-transitive and weak covering map $\pi:\ V(\X)\rightarrow V(\Y)$ has bounded fibres. Then
$p_u(\X)\leq p_u(\Y)\ \mbox{if}\ p_u(\X)<1,$ and
$$p_u(\X)< p_u(\Y)\ \mbox{if the map}\ \pi\ \mbox{has the disjoint tree-lifting property and}\ p_u(\X)<1.$$

\end{enumerate}
Their proof combines a lifted exploration of the cluster with an exploratory
version of the Aizenman-Grimmett strategy \cite{Ai-Gr1991}.

After the work of Martineau and Severo \cite{Ma-Se2019}, Martineau, Poudevigne-Auberon and Rax (2026) \cite{Ma-Po2025} developed
a general stochastic-domination theory for lifts of random variables and
applied it to percolation under graph fibrations. Their framework allows substantially more general lifted random variables than
identically distributed Bernoulli coordinates and revisits the non-strict monotonicity of $p_c$. The $p_c$-strict part of \cite{Ma-Po2025} is independent of their
abstract domination theorem: it gives a coupling proof of the homogeneous bond $p_c$-strict-monotonicity result of \cite{Ma-Se2019} without using essential
enhancements or differential inequalities.
\smallskip

The present paper concerns inhomogeneous one-parameter analogues of both {\bf (a)} and {\bf (b)}. We need that a candidate vector of edge probabilities has a distinguished scalar critical
parameter. So we specify the exponential family
\[
  p_e^{\X}(\lambda)=1-\exp\{-\lambda w_{\X}(e)\},\ e\in E(\X),\ \lambda\in [0,\infty),
\]
where $w_{\X}(e)>0$ is an edge weight. This parametrisation is particularly
convenient for the lifted enhancement argument: replacing an edge of weight
$w_{\X}$ by $N$ independent auxiliary copies of weight $w_{\X}/N$ leaves its physical
opening probability unchanged. For the analogue of {\bf (a)}, we consider the threshold $\lambda_c$ for existence of infinite open-edge clusters in the inhomogeneous percolation. See Theorem~\ref{thm:main}.
To extend {\bf (b)}, we define the threshold $\lambda_{\mathrm{conn}}$, motivated by
Schonmann's homogeneous large-ball characterization \cite{Sc1999}, to characterize the uniform large-ball connections in the inhomogeneous setting. Refer to Theorems~\ref{thm:conn-monotone} and \ref{thm:conn-strict}.
Note for homogeneous percolations, $\lambda_{\mathrm{conn}}$ is
just the uniqueness threshold $\lambda_u.$ We believe this is still true for uniformly elliptic weighted quasi-transitive exponential models, see Conjecture~\ref{conj:stability}.
Moreover, we also prove a general-profile theorem, Theorem \ref{thm:general-profile}, under explicit compact-interval
finite-energy and derivative bounds.

Two assumptions play a central role in our strictness arguments, Theorems~\ref{thm:main} and \ref{thm:conn-strict}. First, the
weights are uniformly elliptical (i.e.\,\,uniformly bounded above and away from zero).  This supplies the
position-independent finite-energy and local-success constants needed by the
enhancement construction.
Second, a lifted edge in the covering graph is at
least as likely to open as the base edge that it lifts. This covering compatibility permits
a coordinatewise monotone coupling.
The latter type of
comparison cannot be omitted for arbitrary deterministic weight profiles; a
simple counterexample is given in Remark~\ref{rem:compatibility}~(i). Uniform
ellipticity is used as a convenient sufficient hypothesis
for strictness rather than an optimal one (see Proposition~\ref{prop:ellipticity} and Example~\ref{ex:vanishing-weight}).

To prove the strictness arguments,
we keep the exploratory-enhancement mechanism (introduced by Aizenman and Grimmett \cite{Ai-Gr1991})
in order to compare an edge-dependent
one-parameter family under the compatibility condition and to obtain a uniform differential comparison that is later reusable for the two thresholds.
On one hand, the strict weighted arguments require $\Hh$-auxiliary coordinates, the corresponding $\G$-auxiliary coordinates, new bookkeeping for edge-dependent
auxiliary coordinates, and the corresponding enhancement marks to construct a uniform thinning coupling. On the other hand, the strict weighted arguments still require
a deterministic pivotal conversion and a uniform pivotal comparison on weighted
Margulis-Russo derivatives to make the thresholds strictly smaller. Our proofs employ the exploratory-enhancement framework of \cite{Ma-Se2019}. Contrast to the homogeneous setting \cite{Ma-Se2019},
where genuine innovation of our proofs lies is in the weighted lifted coupling (in Section~\ref{sec:weighted-lifted-coupling}) for proving strictness arguments in Theorems~\ref{thm:main}, \ref{thm:conn-strict} and \ref{thm:general-profile}.

\smallskip
Note that percolations studied by all above mentioned papers including the present one are nearest-neighbour percolations; and short-range percolation can be viewed as a nearest-neighbour one by changing suitably the underlying graph structure.
Refer to Subsection~\ref{subsec:short-range} for conjectures on the short-range inhomogeneous percolations.

For independent long-range percolation models on locally finite transitive graphs, Bethuelsen and M\"{o}nch (2025) \cite{Bet-Mo2025} proved a nice
strict monotonicity of the critical points with respect to local perturbations in weight functions when the weight functions are shift-invariant and summable. They did not studied the covering type strict-monotonicity, and did not use
the Aizenman-Grimmett approach due to its shortcomings in the long-range percolations. Instead, they developed a new coupling argument, which relies on stochastic domination techniques rather than differential inequalities,
to obtain strict inequality of critical values. For a discussion on the strict monotonicity of the thresholds and conjectures in the long-range setting, see Subsection~\ref{subsec:long-range}.
\smallskip

\section{Definitions and main results}
\setcounter{equation}{0}
\noindent
All graphs are assumed nonempty, connected, locally finite, simple and
unoriented.  For a graph $\X$, write $V(\X)$, $E(\X)$ and $d_{\X}$ respectively for its vertex and edge
sets and graph distance.
Let
\begin{eqnarray}\label{eq:ball}
  && B_r^{\X}(x):=\{y\in \X:\, d_{\X}(x,y)\le r\},\  S_r^{\X}(x):=\{y\in \X:\, d_{\X}(x,y)=r\},\\
\label{eq:e-ball}  && B_r^{\X}(e):=\{y\in \X:\, d_{\X}(y,e)\leq r\},\ x\in V(\X),\ e\in E(\X),\ r\in\mathbb{Z}_{+}=\{0,1,2,\cdots\}.
\end{eqnarray}

\begin{definition}[Weak covering]\label{def:weak-covering}
A map $\pi:V(\G)\to V(\Hh)$ is a \emph{weak covering map} if it is
$1$-Lipschitz and has the weak lifting property: for every $x\in V(\G)$ and every
neighbour $u$ of $\pi(x)$, there exists a neighbour $y$ of $x$ such that
$\pi(y)=u$.
\end{definition}

Under our assumptions that $\G$ is nonempty and $\Hh$ is connected,
every weak covering map is automatically surjective. {\it Indeed},
$\pi(V(\G))$ is nonempty, and the weak lifting property implies that whenever a
vertex belongs to $\pi(V(\G))$, all of its neighbours do as well. By connectivity
of $\Hh$, $\pi(V(\G))=V(\Hh).$ In particular, every vertex and every finite vertex set of $\Hh$ have a lift in
$\G$.

A lift of a subtree $T\subseteq\Hh$ is a subtree $T'\subseteq\G$ on which
$\pi$ induces a graph isomorphism from $T'$ to $T$.

\begin{definition}[Disjoint tree lifting]

The map $\pi$ has the \emph{disjoint tree-lifting property} if, whenever
$x\ne y$, $\pi(x)=\pi(y)\in V(T)$ and $T$ is a subtree of $\Hh$, there exist
vertex-disjoint lifts $T_x,T_y$ of $T$ containing $x$ and $y$ respectively.
\end{definition}

\begin{definition}[Tame fibres]
The fibres of $\pi$ are \emph{tame} if there exists $R<\infty$ such that for
every $x\in V(\G),$ there is a $y\in V(\G)$ satisfying
\[
  \pi(y)=\pi(x),\  0<d_{\G}(x,y)\le R.
\]
\end{definition}

\begin{definition}[Strong covering and bounded fibres]
A weak covering map is a \emph{strong covering map} if the lifted neighbour is
unique: for every $x\in V(\G)$ and every neighbour $u$ of $\pi(x)$, there is a
unique neighbour $y$ of $x$ with $\pi(y)=u$.  The fibres of $\pi$ are
\emph{bounded} if there is a $K<\infty$ such that
\[
  \pi(x)=\pi(y)\ \Longrightarrow\ d_{\G}(x,y)\le K.
\]
Call $K$ a fibre-diameter bound. Note bounded fibre diameter automatically implies
\[
  \pi^{-1}(\pi(x))\subseteq B_K^{\G}(x)\ \mbox{for}\ x\in V(\G).
\]
Consequently $\pi^{-1}(A)$ is finite for any finite $A\subseteq V(\Hh).$
\end{definition}

Let $w_{\X}:E(\X)\to(0,\infty)$ be a weight function.  At intensity
$\lambda\ge0$, weighted Bernoulli percolation on $(\X,w_{\X})$ is the product
measure
\begin{equation}\label{eq:weighted-probability}
\Pp_{\lambda}^{\X}:=\Pp_{\lambda,w_\X}^{\X}\ \mbox{under which}\ \Pp_{\lambda}^{\X}(e\text{ is open})
  =p_e^{\X}(\lambda)=1-\exp\{-\lambda w_{\X}(e)\}.
\end{equation}
The corresponding open cluster of $A\subseteq V(\X)$ is denoted by
$\cC_{\X}^{\lambda}(A)$.

\begin{definition}[Weighted critical intensity]\label{def:lambda-c}
For a root $o\in V(\X)$, define
\[
  \lambda_c(\X,w_{\X})
 :=\inf\left\{\lambda\ge 0:\
       \Pp_{\lambda}^{\X}\left(\left|\cC_{\X}^{\lambda}(o)\right|=\infty\right)>0\right\}\ \left(\inf\varnothing=\infty\right).
\]
\end{definition}

This definition is independent of the root. {\it To see this}, fix $\lambda>0$
and vertices $x,y\in V(\X)$, and let $\gamma$ be a finite path joining them. Write  $A_z:=\left\{\left|\cC_{\X}^{\lambda}(z)\right|=\infty\right\},\ z\in \X.$ Then by Harris-FKG inequality,
\[
 \Pp_{\lambda}^{\X}(A_y)\geq \Pp_{\lambda}^{\X}
 \bigl(A_x\cap\{\gamma\text{ is open}\}\bigr)
 \ge
 \Pp_{\lambda}^{\X}(A_x)
 \prod_{e\in E(\gamma)}p_e^{\X}(\lambda)>0\ \mbox{if}\ \Pp_{\lambda}^{\X}(A_x)>0.
\]
Similarly, $\Pp_{\lambda}^{\X}(A_x)>0$ if $\Pp_{\lambda}^{\X}(A_y)>0.$ Note $\Pp_{0}^{\X}(A_x)=\Pp_{0}^{\X}(A_y)=0.$ The claim holds.

\begin{definition}[Uniform ellipticity and covering compatibility]
Weights $w_{\G}$ and $w_{\Hh}$ are \emph{uniformly elliptic} if there exist
$0<a\le b<\infty$ such that
\begin{equation}\label{eq:ellipticity}
  a\le w_{\G}(f)\le b \ (f\in E(\G)),   \
  a\le w_{\Hh}(e)\le b \ (e\in E(\Hh)).
\end{equation}
They are \emph{compatible with $\pi$} if every edge $f=\{x,y\}$ with $\pi(x)\not=\pi(y)$ satisfies
\begin{equation}\label{eq:compatibility}
  w_{\G}(f)\ge w_{\Hh}(\{\pi(x),\pi(y)\}).
\end{equation}
\end{definition}

\begin{definition}[Weighted quasi-transitivity]\label{def:weighted-quasi-transitive}
The weighted graph $(\X,w_{\X})$ is \emph{weighted quasi-transitive} if
\[
  \Aut(\X,w_{\X})
  :=\{\gamma\in\Aut(\X):w_{\X}(\gamma e)=w_{\X}(e)
       \text{ for every }e\in E(\X)\}
\]
has finitely many orbits on $V(\X)$.
\end{definition}

\paragraph{Terminology convention.}
Throughout the paper, \emph{quasi-transitive graph} refers only to the
underlying unweighted graph, i.e.\,\,to the action of $\Aut(\X)$.  The stronger
phrase \emph{weighted quasi-transitive} always has the meaning of Definition \ref{def:weighted-quasi-transitive}.  When a particular subgroup is needed, we say explicitly that
$\Gamma\le \Aut(\X,w_{\X})$ is a \emph{quasi-transitive weight-preserving
subgroup}. Weighted quasi-transitivity implies quasi-transitivity of the
underlying graph, but the converse may not be true.
\smallskip

When $(\X,w_{\X})$ is infinite and weighted quasi-transitive, every orbit of
$\Aut(\X,w_{\X})$ on $V(\X)$ is infinite, and  $\Pp_\lambda^{\X}$ is ergodic for any $\lambda\geq 0;$
and $\Ninf^{\X}(0)=0,$ and for any $\lambda>0,$
\begin{eqnarray}\label{eq:cluster number}
  \Ninf^{\X}(\lambda)\ \mbox{is a constant in}\ \{0,1,\infty\}\ \text{almost surely}.
\end{eqnarray}
Here $\Ninf^{\X}(\lambda)$ is the number of infinite open clusters in the percolation $\Pp_\lambda^{\X}$ ($\lambda\ge 0$).

{\it Indeed}, {\bf (i)} if a vertex orbit $O$ is finite, then $d_{\X}(\,\cdot\,,O)$ is invariant under
$\Aut(\X,w_{\X})$ and takes arbitrarily large integer values, producing infinitely many vertex orbits. A contradiction! {\bf (ii)} Clearly, $\Pp_0^{\X}$ is ergodic and  $\Ninf^{\X}(0)=0.$ {\bf (iii)}
For $\lambda>0$, $\Pp_\lambda^{\X}$ is invariant under $\Aut(\X,w_{\X})$. The proof of the
ergodicity criterion in \cite[Proposition~7.3]{Ly-Pe2016} uses only
independence, invariance, and that any vertex orbit is infinite, so
it applies to $\Pp_\lambda^{\X}.$ Hence $\Pp_\lambda^{\X}$ is ergodic. Clearly  $\Pp_\lambda^{\X}$ is insertion tolerant.
Therefore, \cite[Theorem~7.5]{Ly-Pe2016} (which uses ergodicity and insertion
tolerance) implies \eqref{eq:cluster number}.
\smallskip

By \eqref{eq:cluster number}, without presupposing the uniqueness phase is an interval, we make the following definition:

\begin{definition}[Weighted uniqueness intensity]\label{def:uniqueness threshold}
For an infinite weighted quasi-transitive graph, define
\[
  \lambda_u(\X,w_{\X})
  =\inf\left\{\lambda\ge0:\Ninf^{\X}(\lambda)=1
       \text{ almost surely}\right\}  \ \left(\inf\varnothing=\infty\right).
\]
\end{definition}

\begin{definition}[Uniform connection intensity]\label{def:lambda-conn}
For any infinite weighted graph $\X$ and $\lambda\ge 0$, set
\[
 a_\ell^{\X}(\lambda)
 :=\inf_{x,y\in V(\X)}
   \Pp_\lambda^{\X}\left[B_\ell^{\X}(x)\leftrightarrow B_\ell^{\X}(y)\right].
\]
Since $B_\ell^{\X}(x)\subseteq B_{\ell+1}^{\X}(x)$ for any $x$, $a_\ell^{\X}(\lambda)\le a_{\ell+1}^{\X}(\lambda),$
$\lim\limits_{\ell\to\infty}a_\ell^{\X}(\lambda)$ exists in $[0,1]$. Define
\begin{equation}\label{eq:def-lambda-conn}
 \lambda_{\mathrm{conn}}(\X,w_{\X})
 :=\inf\left\{\lambda\ge0:
 \lim_{\ell\to\infty}a_\ell^{\X}(\lambda)=1\right\}\ \left(\inf\varnothing=\infty\right).
\end{equation}
This is a genuine threshold because $\lim\limits_{\ell\to\infty}a_\ell^{\X}(\lambda)$ is increasing in $\lambda.$
\end{definition}

\begin{definition}[Stability property]\label{def:stability}
A weighted quasi-transitive exponential model has property $\mathsf{Stab}$ if:
\begin{enumerate}[label=\textnormal{(S\arabic*)},leftmargin=2.8em]
 \item uniqueness at intensity $\lambda$ implies uniqueness at any $\lambda'>\lambda$;
 \item if the limit in \eqref{eq:def-lambda-conn} equals one at $\lambda$, then
       every $\lambda'>\lambda$ has a unique infinite cluster.
\end{enumerate}
\end{definition}

In homogeneous Bernoulli percolation,
(S1) follows from the cluster-containment and
merging results in the canonical coupling (Schonmann~\cite{Sc1999}, H\"{a}ggstr\"{o}m-Peres-Schonmann
\cite{Ha-Pe1999}); and (S2) is a consequence of Schonmann's homogeneous stability theorem (\cite{Sc1999}).

For the exponential weighted model, the natural clock coupling does provide
uniform sprinkling whenever $w(e)\ge a$: conditionally on an edge being closed
at $\lambda$, its probability of being added by $\lambda'>\lambda$ is at
least $1-e^{-(\lambda'-\lambda)a}$.  This lower bound does not by itself prove
either the global cluster-containment statement underlying (S1) or the
large-ball-to-uniqueness implication in (S2).

The weight-preserving assumption is substantive. Generally, $\Pp_{\lambda}^{\X}$ need not be invariant under any quasi-transitive automorphism group and
the number of infinite clusters need not be controlled by the standard ergodic
argument. Note we prove in Lemma \ref{lem:connection-criterion}(ii) that for infinite weighted quasi-transitive graph $\X$, when the model $\Pp_{\lambda}^{\X}$ ($\lambda\geq 0$) has the property $\mathsf{Stab},$
\begin{eqnarray}\label{eq:connection-uniqueness}
\lambda_{\mathrm{conn}}(\X,w)=\lambda_u(\X,w);
\end{eqnarray}
and think \eqref{eq:connection-uniqueness} is still true without $\mathsf{Stab}$ (see Conjecture \ref{conj:stability}).\\

Our main results are described as follows. The non-strict comparison, Theorem \ref{thm:pc-monotone}, on $\lambda_c$ holds under weak covering alone. While
Theorem \ref{thm:conn-monotone} on the non-strict monotonicity for $\lambda_{\mathrm{conn}}$ holds under weak covering with bounded fibres.
The first strict result, Theorem \ref{thm:main}, is a weighted inhomogeneous extension of \cite[Theorem~1.1]{Ma-Se2019}. The second strict result, Theorem \ref{thm:conn-strict}, can be viewed as a weighted inhomogeneous generalization of
\cite[Theorem~1.4]{Ma-Se2019} since $\lambda_{\mathrm{conn}}=\lambda_u$ holds for homogeneous Bernoulli bond percolation. We prove Theorems \ref{thm:pc-monotone} and \ref{thm:conn-monotone} in \ref{app-1}, Theorem
\ref{thm:main} in Sections \ref{sec:exploratory}-\ref{sec:strict}, and Theorem \ref{thm:conn-strict} in Section \ref{sec:uniqueness}. In Section \ref{sec:conjectures}, we give a general version of these theorems by Theorem \ref{thm:general-profile}.

\begin{theorem}[Weighted covering monotonicity]\label{thm:pc-monotone}
Let $\pi:V(\G)\to V(\Hh)$ be a weak covering map between connected locally
finite graphs, and suppose the weights are compatible with $\pi$.  Then
\[
  \lambda_c(\G,w_{\G})\le\lambda_c(\Hh,w_{\Hh}).
\]
\end{theorem}
\smallskip

\begin{theorem}[Weighted strict monotonicity]\label{thm:main}
Let $\G$ and $\Hh$ be graphs of bounded degree.  Suppose that
$\pi:V(\G)\to V(\Hh)$ is a weak covering map with tame fibres and the disjoint
tree-lifting property.  Let $w_{\G}$ and $w_{\Hh}$ be uniformly elliptic and
compatible with $\pi$.  If $\lambda_c(\G,w_{\G})<\infty,$
then
\[
  \lambda_c(\G,w_{\G})
  <\lambda_c(\Hh,w_{\Hh}).
\]
\end{theorem}
\smallskip

\begin{theorem}[Monotonicity of the uniform connection intensity]
\label{thm:conn-monotone}
Let $\G$ and $\Hh$ be infinite connected locally finite graphs.  Suppose that
$\pi:V(\G)\to V(\Hh)$ is a weak covering map with bounded fibres and that the
weights are compatible with $\pi$.  Then
\[
  \lambda_{\mathrm{conn}}(\G,w_{\G})
  \le \lambda_{\mathrm{conn}}(\Hh,w_{\Hh}).
\]
\end{theorem}
\smallskip

\begin{theorem}[Strict monotonicity of the uniform connection intensity]
\label{thm:conn-strict}
Let $\G$ and $\Hh$ be infinite quasi-transitive graphs of bounded degree.
Suppose that $\pi:V(\G)\to V(\Hh)$ is a non-injective weak covering map with
bounded fibres and the disjoint tree-lifting property.  Let the weights be
uniformly elliptic and compatible with $\pi$.  If
$\lambda_{\mathrm{conn}}(\G,w_{\G})<\infty,$
then
\[
  \lambda_{\mathrm{conn}}(\G,w_{\G})
  <\lambda_{\mathrm{conn}}(\Hh,w_{\Hh}).
\]
\end{theorem}
\smallskip

\begin{corollary}[Conditional uniqueness comparison]\label{cor:pu-conditional}
Assume in addition to Theorem~\ref{thm:conn-monotone} that both
$(\G,w_{\G})$ and $(\Hh,w_{\Hh})$ are weighted quasi-transitive and satisfy
$\mathsf{Stab}$.  Then
\[
 \lambda_u(\G,w_{\G})\le\lambda_u(\Hh,w_{\Hh}).
\]
In case of Theorem~\ref{thm:conn-strict}, if $(\G,w_{\G})$ and $(\Hh,w_{\Hh})$ are weighted quasi-transitive and satisfy
$\mathsf{Stab}$, then
\[
 \lambda_u(\G,w_{\G})<\lambda_u(\Hh,w_{\Hh}).
\]
\end{corollary}
\smallskip

\begin{remark}\label{rem:compatibility} {\bf (i)} Compatibility is necessary. For example, let $\Hh=\mathbb Z^2$ and $\G=\mathbb Z^2\times\{0,1\}$, and let
$\pi$ forget the $\{0,1\}$ coordinate. Then $\pi$ is a strong covering map; its fibres have diameter
one, and hence are tame. Note $p_c(\G)<p_c(\Hh).$ Define $
  \kappa(\X):=-\log(1-p_c(\X)),\ \X\in\{\G,\Hh\}.$
Put $w_{\Hh}\equiv1$ and $w_{\G}\equiv\varepsilon$. Then
$$\lambda_c(\G,w_{\G})=\frac{\kappa(\G)}{\varepsilon},\ \lambda_c(\Hh,w_{\Hh})=\kappa(\Hh).$$
Choosing $0<\varepsilon<\frac{\kappa(\G)}{\kappa(\Hh)}$ reverses the desired inequality. Thus a condition comparing the probabilities
of a base edge and its lifts is essential.

{\bf (ii)} Uniform ellipticity is a convenient sufficient hypothesis, but not the optimal one. For instance, as shown in Proposition~\ref{prop:ellipticity} and Example~\ref{ex:vanishing-weight},
when $w_{\G}$ is not of uniform ellipticity, the strict monotonicities can still hold. Uniform ellipticity supplies the uniform finite-energy and local-success constants used below; and without a uniform lower bound, the extra lifted route may
become arbitrarily unlikely at distant locations.

{\bf (iii)} Set $p=1-\exp\{-\lambda\}$ and $p_e(p)=1-(1-p)^{w(e)}.$
The map $\lambda\mapsto p=1-\exp\{-\lambda\}$ is strictly increasing from
$[0,\infty]$ onto $[0,1].$ So this is only a
reparametrisation of the same weighted family.  If
\[
  p_c^{w}(\X):=1-\exp\{-\lambda_c(\X,w_{\X})\},
\]
then every non-strict or strict comparison for $\lambda_c$ is equivalent to
the corresponding comparison for $p_c^{w}$.  In particular, by
Theorem~\ref{thm:main}, $p_c^{w}(\G)<p_c^{w}(\Hh).$
\end{remark}

\section{Exploratory enhancement}\label{sec:exploratory}
\setcounter{equation}{0}
\noindent
Fix $r\ge1$ and a root $o\in V(\Hh)$.  Let
\[
  \omega\in\{0,1\}^{E(\Hh)},\ \alpha\in\{0,1\}^{V(\Hh)}.
\]
Starting from $C_0=\{o\}$, define an increasing sequence of vertex sets.  Given
$C_{2n}$, let $C_{2n+1}$ be the union of the $\omega$-open clusters of all
vertices in $C_{2n}$.  Then let
\begin{equation}\label{eq:enhancement-step}
\begin{split}
  C_{2n+2}:=C_{2n+1}\cup\biggl\{v:\ &\text{there exists }u\in C_{2n+1}
  \text{ such that }d_{\Hh}(u,v)=r+1,\\
  &E\left(B_r^{\Hh}(u)\right)\subseteq\omega,\ \alpha_u=1\biggr\}.
\end{split}
\end{equation}
Set
\[
  C_o(\omega,\alpha):=\bigcup_{n\ge0}C_n.
\]

\begin{lemma}[Connectedness and sphere hitting]
\label{lem:enhanced-connected}
For every configuration $(\omega,\alpha)$, each explored set $C_n$ is connected
in the underlying graph $\Hh$. Consequently $C_o(\omega,\alpha)$ is connected,
and for every integer $L\ge1$,
\[
  C_o(\omega,\alpha)\cap S_L^{\Hh}(o)\ne\varnothing
  \ \Longleftrightarrow\  C_o(\omega,\alpha)\ \text{contains a vertex at distance at least}\ L\ \text{from}\ o.
\]
Particularly, the sphere-hitting events
\[
  \cE_L:=\left\{C_o(\omega,\alpha)\cap S_L^{\Hh}(o)\ne\varnothing\right\}
\]
are decreasing in $L$, and
\begin{equation}\label{eq:sphere-hitting-limit}
  \bigcap_{L\ge1}\cE_L
  =\{|C_o(\omega,\alpha)|=\infty\}.
\end{equation}
\end{lemma}

\begin{proof}
We argue inductively. Note $C_0=\{o\}$ is connected. If $C_{2n}$ is
connected, then every $\omega$-open cluster used to form $C_{2n+1}$ meets
$C_{2n}$; hence their union $C_{2n+1}$ is connected. Now suppose that
$v\in C_{2n+2}\setminus C_{2n+1}$ is inserted by an enhancement centred at
$u\in C_{2n+1}$. The condition $E(B_r^{\Hh}(u))\subseteq\omega$ implies that every
vertex of $B_r^{\Hh}(u)$ lies in the ordinary closure $C_{2n+1}$. Choose a geodesic
$u=x_0,x_1,\ldots,x_{r+1}=v$. Then $x_r\in B_r^{\Hh}(u)\subseteq C_{2n+1}$ and
$x_r$ is adjacent to $v$. Namely every newly inserted vertex is adjacent, in the
underlying graph, to $C_{2n+1}$, so $C_{2n+2}$ is connected as well. Hence, each $C_n$ is
connected, and so is the increasing union $C_o(\omega,\alpha)$.

If a connected set containing $o$ contains a vertex $y$ with $d_{\Hh}(o,y)\ge L$,
a path inside that set from $o$ to $y$ must contain a vertex at distance
exactly $L$ from $o$. This proves the equivalence and hence
$\cE_{L+1}\subseteq\cE_L$. Finally, local finiteness implies that every bounded
ball is finite. Therefore an infinite connected set containing $o$ is
unbounded and meets every sphere, while a set meeting every sphere is clearly
infinite. This proves \eqref{eq:sphere-hitting-limit}.
\end{proof}

Given any $(\lambda,s)\in [0,\infty)\times [0,1].$ Under the product measure $\mathbb{P}^{\Hh}_{\lambda,s}$ on
$$\{0,1\}^{E(\Hh)}\times\{0,1\}^{V(\Hh)}$$
for which the edges $e$ have probabilities
$p_e^{\Hh}(\lambda)=1-\exp\left\{-\lambda w_{\Hh}(e)\right\}$
and the variables $\alpha_x$ are independent
$\Ber(s)$ variables, write $\cC_{\Hh}^{\lambda,s}(o)$ for the distribution of $C_o(\omega,\alpha)$. Here
$$\Ber(s)=s\delta_1+(1-s)\delta_0.$$

The following two propositions imply Theorem~\ref{thm:main}.

\begin{proposition}[Weighted lifted coupling]\label{prop:coupling}
Assume that $\G$ and $\Hh$ have bounded degree, that
$\pi:V(\G)\to V(\Hh)$ is a weak covering map with tame fibres and the disjoint
tree-lifting property, and that $w_{\G}$ and $w_{\Hh}$ are uniformly elliptic
and compatible with $\pi$.  Then there exists $r\ge1$ such that for any
$\lambda_0>0$, there is an $s\in(0,1)$ for which
\[
  \cC_{\Hh}^{\lambda,s}(o)
  \preceq
  \pi\bigl(\cC_{\G}^{\lambda}(o')\bigr)\ \text{for any $\lambda\ge\lambda_0$ and any $o'\in\pi^{-1}(o)$.}
\]
Here $\preceq$ means that there is a coupling under which $C_o(\omega,\alpha)$ is contained in $\pi\bigl(\cC_{\G}^{\lambda}(o')\bigr)$ almost surely.
\end{proposition}

\begin{proposition}[Strictness of the enhancement]\label{prop:enhancement}
Assume that $\Hh$ has bounded degree and that there exist constants
$0<a\le b<\infty$ such that
\[
  a\le w_{\Hh}(e)\le b\ \,(e\in E(\Hh)).
\]
If $\lambda_c(\Hh,w_{\Hh})<\infty$, then for any $r\ge1$ and any $s\in(0,1],$ there exists
\[
  \lambda_s\in(0,\lambda_c(\Hh,w_{\Hh}))
\]
such that $C_o(\omega,\alpha)$ is infinite with positive probability under $\cC_{\Hh}^{\lambda,s}(o)$ for
all $\lambda\ge\lambda_s$.
\end{proposition}

\begin{proof}[Proof of Theorem~\ref{thm:main} by Propositions \ref{prop:coupling}-\ref{prop:enhancement}]
If $\lambda_c(\Hh,w_{\Hh})=\infty$, the theorem holds trivially. So assume
\[
  0<\lambda_c(\Hh,w_{\Hh})<\infty.
\]
Choose
\[
  0<\lambda_0<\lambda_c(\Hh,w_{\Hh}).
\]
Proposition~\ref{prop:coupling} supplies $r$ and $s\in(0,1)$ such that for any
$\lambda\ge\lambda_0$,
\[
  \cC_{\Hh}^{\lambda,s}(o)
  \preceq\pi\bigl(\cC_{\G}^{\lambda}(o')\bigr).
\]
By Proposition~\ref{prop:enhancement}, there is
\[
  \lambda_s<\lambda_c(\Hh,w_{\Hh})
\]
such that $C_o(\omega,\alpha)$ is infinite with positive probability under $\cC_{\Hh}^{\lambda,s}(o)$ for any
$\lambda\ge\lambda_s$. Set
\[
  \lambda=\max\{\lambda_0,\lambda_s\}
  <\lambda_c(\Hh,w_{\Hh}).
\]
Then  $C_o(\omega,\alpha)$ is infinite with positive probability under $\cC_{\Hh}^{\lambda,s}(o)$, and hence so is its dominating set
$\pi(\cC_{\G}^{\lambda}(o'))$. This implies $\cC_{\G}^{\lambda}(o')$ is infinite with positive probability. Therefore
\[
  \lambda_c(\G,w_{\G})
  \le\lambda
  <\lambda_c(\Hh,w_{\Hh}).
\]
\end{proof}

\section{The weighted lifted coupling}
\label{sec:weighted-lifted-coupling}
\setcounter{equation}{0}

\noindent
In this section we prove Proposition~\ref{prop:coupling}.  Throughout the
section, $\pi:\ V(\G)\to V(\Hh)$ satisfies the hypotheses of
Theorem~\ref{thm:main}: the graphs have bounded degree, the fibres of $\pi$
are tame, $\pi$ has the disjoint tree-lifting property, and the weights are
uniformly elliptic and compatible with $\pi$.

The graphs are connected and locally finite, their vertex and edge sets
are countable.  We therefore fix once and for all enumerations of
$V(\G),\, E(\G),\, V(\Hh)$ and $E(\Hh),$ and use them for every tie-breaking choice
below.  In particular, all spanning trees, lifted edges and duplicated exits
used in the construction are chosen deterministically from the information
already revealed at the time of the choice.

\subsection{The geometric input}
\noindent
We begin with a slightly strengthened form of the duplicated-exit lemma (see \cite[Lemma 2.3]{Ma-Se2019}). The
strengthening records that the two exits supplied for a boundary vertex are
lifts of the same base edge.

\begin{lemma}[Uniform duplicated exits]\label{lem:geometry}
There exists an integer $r\ge1$ such that the following holds.  For
$x\in V(\G)$, let
\[
 Z(x,r):=\text{the connected component of $x$ in }
 \pi^{-1}(B_r^{\Hh}(\pi(x)))\cap B_{3r}^{\G}(x).
\]
Then for any $v\in S_{r+1}^{\Hh}(\pi(x))$, there exist a vertex
$p_v\in S_r^{\Hh}(\pi(x))$ and two distinct edges
$f_v^{(1)},f_v^{(2)}\in E(\G)$ such that
\[
 \left\{\pi(g_i),\pi\left(g'_i\right)\right\}=\{p_v,v\} \ \text{for}\ f_v^{(i)}=\left\{g_i,g'_i\right\},\ i=1,2,
\]
and each $f_v^{(i)}$ has one endpoint in $Z(x,r)$ and the other endpoint in
$\pi^{-1}(v)$.
\end{lemma}

\begin{proof}
Let $R$ be the tameness constant and put $r=\lceil R/2\rceil$. Here for any $t\in\mathbb{R},$ $\lceil t\rceil$ is the smallest integer no less than $t.$ Fix
$x\in V(\G)$. By tameness, there exists $y\ne x$ such that
\[
 \pi(y)=\pi(x),\ d_{\G}(x,y)\le R\le2r.
\]
Let $T$ be a breadth-first spanning tree rooted at $\pi(x)$ of $B_{r+1}^{\Hh}(\pi(x))$. By the disjoint tree-lifting property, there are vertex-disjoint
lifts $T_x$ and $T_y$ of $T$ containing $x$ and $y$ respectively.

Let $\gamma$ be a geodesic from $x$ to $y$ of length $m\le 2r$. If $z$ is
the $k$th vertex of $\gamma$, then from the $1$-Lipschitz property and
$\pi(x)=\pi(y),$
\[
 d_{\Hh}(\pi(z),\pi(x))\le\min\{k,m-k\}\le r.
\]
Hence $\gamma\subseteq\pi^{-1}(B_r^{\Hh}(\pi(x)))$. The union of $\gamma$ with
the portions of $T_x\cup T_y$ lying above $B_r^{\Hh}(\pi(x))$ is connected,
contains $x$, and lies in $B_{3r}^{\G}(x)$. It is therefore contained in
$Z(x,r)$.

Now let $v\in S_{r+1}^{\Hh}(\pi(x))$. Since $T$ is a breadth-first tree, the
predecessor $p_v$ of $v$ in $T$ belongs to $S_r^{\Hh}(\pi(x))$.  The tree edge
$\{p_v,v\}$ has one lift in $T_x$ and one lift in $T_y$.  Their inner
endpoints lie in $Z(x,r)$, and their outer
ones belong to $\pi^{-1}(v)$.  Since $T_x$ and $T_y$ are vertex-disjoint,
these two lifted edges are distinct.
\end{proof}

\subsection{Auxiliary coordinates and deterministic reservations}
\noindent
Fix $r$ as in Lemma~\ref{lem:geometry}.  Put
\[
 M_r:=\sup_{u\in V(\Hh)}\left|B_r^{\Hh}(u)\right|<\infty,
\]
and choose an integer
\begin{equation}\label{eq:N-choice}
 N>2M_r.
\end{equation}
The finiteness of $M_r$ follows from the bounded-degree assumption.

For $e\in E(\Hh)$ and $f\in E(\G),$ define
\begin{equation}\label{eq:auxiliary-probabilities}
 q_e^{\Hh}(\lambda)
 :=1-\exp\!\left(-\frac{\lambda w_{\Hh}(e)}{N}\right),
 \  q_f^{\G}(\lambda)
 :=1-\exp\!\left(-\frac{\lambda w_{\G}(f)}{N}\right).
\end{equation}
Thus
\[
 1-\left(1-q_e^{\Hh}(\lambda)\right)^N=1-\exp\left\{-\lambda w_{\Hh}(e)\right\},\ \text{and similarly for $\G$}.
\]

For an edge $f=\{x,y\}\in E(\G),$ define the set of possible enhancement
centres that can use $f$ by
\[
 \mathcal U_r(f):=B_r^{\Hh}(\pi(x))\cup B_r^{\Hh}(\pi(y)).
\]
By \eqref{eq:N-choice},
\[
 |\mathcal U_r(f)|\le2M_r<N.
\]
Hence, for any $f\in E(\G)$, we may fix an injection
\begin{equation}\label{eq:reservation-map}
 \kappa_f:\ \mathcal U_r(f)\longrightarrow\{1,\dots,N\}.
\end{equation}
The coordinate $(f,\kappa_f(u))$ will be called the coordinate
\emph{reserved for centre $u$ on $f$}.

We work on a probability space carrying three mutually independent families
\[
 \{U_{e,k}:\ e\in E(\Hh),\ 1\le k\le N\},
 \
 \{V_{f,k}:\ f\in E(\G),\ 1\le k\le N\},
 \
 \{W_u:\ u\in V(\Hh)\},
\]
all consisting of independent $\operatorname{Unif}[0,1]$ random variables. Here $\operatorname{Unif}[0,1]$ is the uniform distribution on $[0,1].$
Their roles, including the information that may be used before a reveal,
are as follows.
\begin{center}
\small
\renewcommand{\arraystretch}{1.16}
\begin{tabular}{@{}p{0.13\textwidth}p{0.38\textwidth}p{0.42\textwidth}@{}}
\hline
Source & Use & Selection before inspection \\
\hline
$U_{e,k}$ & Base edge and its assigned ordinary lift share this uniform.
& The lifted edge is selected from the ordinary exploration history before
the entire $U_{e,\cdot}$-block is exposed. \\
$V_{f,k}$ & Reserved enhancement coordinate, or final value of an
unassigned lifted coordinate.
& The requested label is selected from the history; an unused filler is
never inspected to decide whether it will be used. \\
$W_u$ & Thinning to a mark of common parameter $s$, or a mark at an
ineligible or never processed centre.
& Eligibility and the trial success probability are determined without
inspecting $W_u$. \\
\hline
\end{tabular}
\end{center}

The $\Hh$-auxiliary coordinates are
\[
 \xi^{\Hh}(e,k):=\one_{\left\{U_{e,k}\le q_e^{\Hh}(\lambda)\right\}},
\]
and the associated physical base configuration is
\[
 \omega(e):=\max_{1\le k\le N}\xi^{\Hh}(e,k).
\]
Consequently $\omega$ has exactly the law $\Pp_\lambda^{\Hh}$.

The $\G$-auxiliary coordinates $\xi^{\G}(f,k)$ will be revealed adaptively (see \eqref{eq:G-coordinate-1}, \eqref{eq:G-coordinate-2} and \eqref{eq:G-coordinate-3} below).
A physical edge $f\in E(\G)$ will
be called \emph{certified open} as soon as one of its already revealed
$\Hh$-auxiliary coordinates $\xi^{\Hh}(\pi(f),k)$ has value one.  During the ordinary exploration, every
tested base edge $e$ is assigned exactly one lifted edge, denoted by
$\varphi(e)\in E(\G)$.  The map $\varphi$ is defined only on base edges that
have actually been ordinarily tested.

\subsection{Ordinary exploration}
\noindent
We construct simultaneously a set $C\subseteq V(\Hh)$ and a set
$C^\uparrow\subseteq V(\G)$.  Initially,
\[
 C=\{o\},\ C^\uparrow=\{o'\},\ o'\in\pi^{-1}(o).
\]
We maintain
\begin{equation}\label{eq:projection-invariant}
 \pi(C^\uparrow)=C,
\end{equation}
and every vertex of $C^\uparrow$ is connected to $o'$ by physical edges of
$\G$ already certified to be open.

An ordinary phase explores every not-yet-tested base edge incident to the
current set $C$.  We use a fixed fair exploration schedule: every edge which
becomes incident to $C$ is eventually tested.  Suppose that a previously
untested edge $e=\{u,v\}\in E(\Hh)$ is tested from an endpoint $u\in C$.
Choose, using the fixed tie-breaking rule, a vertex
\[
 x\in C^\uparrow\cap\pi^{-1}(u).
\]
By the weak lifting property, there exists an edge $f=\{x,y\}\in E(\G)$ such
that $\pi(y)=v$.  We shall prove below that, at this moment, no $\G$-auxiliary
coordinate of any lift of $e$ has previously been used.  We therefore set
$\varphi(e)=f$ and reveal the whole block
\begin{eqnarray}\label{eq:G-coordinate-1}
 \xi^{\G}(f,k):=\one_{\left\{U_{e,k}\le q_f^{\G}(\lambda)\right\}},
 \ 1\le k\le N\ \left(\mbox{the $\G$-auxiliary coordinates of $f=\varphi(e)$}\right).
\end{eqnarray}
By compatibility,
\[
 q_f^{\G}(\lambda)\ge q_e^{\Hh}(\lambda),
\]
and therefore
\begin{equation}\label{eq:coordinatewise-domination}
 \xi^{\Hh}(e,k)\le\xi^{\G}(f,k), \ 1\le k\le N.
\end{equation}
If $\omega(e)=1$, then some $\xi^{\Hh}(e,k)$ equals one, so
\eqref{eq:coordinatewise-domination} implies that $f$ is certified open.  We
then
$$\mbox{add $y$ to $C^\uparrow$ and $v$ to $C$}.$$

The ordinary phase continues until every base edge incident to the resulting
$C$ has been tested. If infinitely many tests are required, we take the
increasing union along the fixed fair enumeration.  The resulting $C$ is
precisely
$$\mbox{the union of the $\omega$-open clusters of the vertices present at
the beginning of the ordinary phase.}$$

\paragraph{Countable phases.}
An ordinary phase has at most countably many finite-block tests, and an
enhancement phase has at most countably many centre tests. A finite phase
can be padded by idle steps. Thus all reveals may be indexed by the fixed
countable well-order of pairs $(m,n)\in\mathbb N^2$, ordered first by
phase $m$ and then by within-phase step $n$.
At the end of a phase take increasing unions of the explored sets and the
sigma-field generated by the preceding histories. No new random variable is
inspected at such a limit. The same applies at the final union of all phases.
Lemma~\ref{lem:adaptive-deferred} below verifies the deferred-decisions
argument for precisely these countable limits.

\subsection{Freshness at an enhancement centre}
\noindent
The following lemma contains the bookkeeping needed in the enhancement phase.

\begin{lemma}[Fresh-coordinate bookkeeping]\label{lem:fresh-coordinate}
Suppose that an ordinary phase has just been completed and that $u\in C$
satisfies
\begin{equation}\label{eq:fully-open-ball}
 E\left(B_r^{\Hh}(u)\right)\subseteq\omega.
\end{equation}
Then the following statements hold.
\begin{enumerate}[label=\textnormal{(\roman*)},leftmargin=2.6em]
\item $B_r^{\Hh}(u)\subseteq C$, and every base edge having at least one endpoint in
$B_r^{\Hh}(u)$ has already been ordinarily tested.

\item Let $x\in C^\uparrow\cap\pi^{-1}(u)$ and $f$ be an edge of a spanning
tree of $Z(x,r)$. If $f=\varphi(e)$ for some base edge $e$, then $f$ is already
certified open.

\item For each $v\in S_{r+1}^{\Hh}(u)$, let $f_v^{(1)},f_v^{(2)}$ be the two exits
supplied by Lemma~\ref{lem:geometry}. Their common projected base edge has
already been ordinarily tested, and at least one of $f_v^{(1)},f_v^{(2)}$ is
not the assigned ordinary lift of that base edge.

\item If an enhancement centred at $u$ uses a lifted edge $f$ which is not
already certified open and is not an ordinary assigned lift, then the reserved
coordinate $(f,\kappa_f(u))$ has not previously been revealed.

\item No enhancement ever reveals an auxiliary coordinate of a lift of a base
edge before that base edge has been ordinarily tested.  Consequently, when a
previously untested base edge is first tested during an ordinary phase, all
$N$ auxiliary coordinates of every one of its lifts are still unrevealed.
\end{enumerate}
\end{lemma}

\begin{proof}
Since all edges of $B_r^{\Hh}(u)$ are open, every vertex of $B_r^{\Hh}(u)$ is connected to
$u$ by an open geodesic contained in $B_r^{\Hh}(u)$. The preceding ordinary phase
was run to exhaustion, so $B_r^{\Hh}(u)\subseteq C$. By definition of the ordinary
exploration, every edge incident to $C$, and hence every edge having an
endpoint in $B_r^{\Hh}(u)$, has already been tested. This proves (i).
\smallskip

For (ii), suppose that $f=\varphi(e)$ and that $f$ belongs to a spanning tree
of $Z(x,r)$. Since an ordinary assigned lift has endpoints with distinct
images, its projected base edge $e$ has both endpoints in $B_r^{\Hh}(u)$. Thus
$e\in E\left(B_r^{\Hh}(u)\right)$, and $\omega(e)=1$ by \eqref{eq:fully-open-ball}. At the
ordinary test of $e$, \eqref{eq:coordinatewise-domination} therefore certified
$\varphi(e)=f$ open.
\smallskip

For (iii), each exit edge has an inner endpoint whose image lies in $S_r^{\Hh}(u)$,
so its projected base edge is incident to $B_r^{\Hh}(u)$ and has already been tested
by (i).  The two exits are distinct lifts of the same base edge, whereas a
base edge receives only one assigned ordinary lift.  Hence at most one of the
two exits is the assigned lift.
\smallskip

We prove (iv) and (v) together.  If an enhancement centred at $z$ uses an edge
$g=\{g_1,g_2\}\in E(\G)$, then either $g$ lies in $Z(x_z,r)$ or it is one of the selected
exits, where $\pi(x_z)=z.$  In both cases
\[
 z\in B_r^{\Hh}(\pi(g_1))\cup B_r^{\Hh}(\pi(g_2))=\mathcal U_r(g).
\]
Hence the coordinate reserved for this use is $(g,\kappa_{g}(z))$. Different
enhancement centres use different reserved indices on the same physical edge,
because $\kappa_{g}$ is injective, and each centre is tested at most once. Thus
no previous enhancement centred at a different vertex can have revealed
$(f,\kappa_f(u))$. By the assumption in (iv), the edge $f$ is not an ordinary
assigned lift, so no ordinary exploration has revealed all of its coordinates.
This proves (iv).

Finally, every edge $f$ of $\G$ used by an enhancement centred at $z$ has a
projected base edge $e$ with an endpoint in $B_r^{\Hh}(z)$.  By (i), $e$ has already been ordinarily tested.
Hence an enhancement cannot touch a
lift of a still untested base edge. This proves (v).
\end{proof}

\subsection{The enhancement phase}
\noindent
Let an ordinary phase have been completed. Denote its terminal sets by
$C_{\mathrm{fr}}$ and $C^\uparrow_{\mathrm{fr}}$, and keep these two
snapshots fixed throughout the enhancement phase. In particular,
$C_{\mathrm{fr}}$ is closed under ordinary $\omega$-connectivity. The
list of centres for this phase consists exactly of the vertices of
$C_{\mathrm{fr}}$ whose enhancement variables have not yet been assigned.
Order this list by the fixed enumeration. Separately initialize cumulative
output sets by
\[
 C_{\mathrm{out}}:=C_{\mathrm{fr}},\
 C^\uparrow_{\mathrm{out}}:=C^\uparrow_{\mathrm{fr}}.
\]
Only these output sets grow during the phase. The revealed history and the
set of certified-open edges also grow as tests are performed; they are not
frozen. Consider a centre $u$ from the fixed list.

The eligibility test is measurable from the completed ordinary history.
If $B_r^{\Hh}(u)\subseteq C_{\mathrm{fr}}$, every edge of $B_r^{\Hh}(u)$ has already been tested.
Otherwise choose a geodesic inside $B_r^{\Hh}(u)$ from $u$ to a vertex outside
$C_{\mathrm{fr}}$. Its first edge leaving $C_{\mathrm{fr}}$ was tested and
is closed, since $C_{\mathrm{fr}}$ is ordinarily closed. That edge already certifies
$E\left(B_r^{\Hh}(u)\right)\not\subseteq\omega$. Thus deciding eligibility never inspects
an untested base uniform. This argument always uses the frozen snapshot,
even after other centres have enlarged $C_{\mathrm{out}}$.


If $E\left(B_r^{\Hh}(u)\right)\not\subseteq\omega$, set
\[
 \alpha_u:=\one_{\{W_u\le s\}}.
\]
This value can never influence the enhanced exploration at $u$, because the
condition of a fully open $r$-ball fails permanently.

Suppose now that $E\left(B_r^{\Hh}(u)\right)\subseteq\omega$.  Choose
\[
 x_u\in C^\uparrow_{\mathrm{fr}}\cap\pi^{-1}(u)
\]
according to the fixed tie-breaking rule, put $Z_u:=Z(x_u,r)$, and choose a
fixed spanning tree $\mathcal T_u$ of $Z_u$.

For every edge $f$ of $\mathcal T_u$ which is not already certified open,
Lemma~\ref{lem:fresh-coordinate}(ii) implies that $f$ is not an ordinary
assigned lift.  We therefore designate the fresh coordinate
$(f,\kappa_f(u))$.

Next, for every $v\in S_{r+1}^{\Hh}(u)$, choose the duplicated pair
$f_v^{(1)},f_v^{(2)}$ supplied by Lemma~\ref{lem:geometry}.  By
Lemma~\ref{lem:fresh-coordinate}(iii), at least one exit is not
the assigned ordinary lift of their common projected base edge.  Choose one
such edge and denote it by $g_{u,v}$.  If $g_{u,v}$ is already certified open,
no new coordinate is needed; otherwise designate
$(g_{u,v},\kappa_{g_{u,v}}(u))$.

Let $\mathcal D_u$ be the finite set of all designated coordinates.  By
Lemma~\ref{lem:fresh-coordinate}, all coordinates in $\mathcal D_u$ are fresh.
Moreover, if $u\ne u'$, then $\mathcal D_u\cap\mathcal D_{u'}=\varnothing$,
because different centres receive different reserved indices on every
physical edge.

Reveal, for $(f,k)\in\mathcal D_u$,
\begin{eqnarray}\label{eq:G-coordinate-2}
 \xi^{\G}(f,k):=\one_{\left\{V_{f,k}\le q_f^{\G}(\lambda)\right\}}\ \left(\mbox{the $\G$-auxiliary coordinate of designated $(f,k)\in\mathcal D_u$}\right),
\end{eqnarray}
and let $A_u$ be the event that all these newly revealed coordinates are open.
All physical edges involved in $\mathcal D_u$ lie in $B_{3r+1}^{\G}(x_u)$.  Hence,
with
\[
 K:=\sup_{x\in V(\G)}\left|E\left(B_{3r+1}^{\G}(x)\right)\right|<\infty,
\]
we have $|\mathcal D_u|\le K$.  Fix $\lambda_0>0$ and put
\[
 q_*:=1-\exp\!\left(-\frac{\lambda_0a}{N}\right)>0.
\]
For every $\lambda\ge\lambda_0$ and every fresh coordinate,
$q_f^{\G}(\lambda)\ge q_*$. Let
$\mathcal{U}=\sigma\left(U_{e,k}:\,e\in E(\Hh),\,1\le k\le N\right).$
Conditionally on the past construction and on $\mathcal{U},$
the coordinates in $\mathcal D_u$ are
independent, and therefore
\begin{equation}\label{eq:qu-lower-bound}
 q_u:=\Pp(A_u\mid\text{past},\mathcal{U})
 =\prod_{(f,k)\in\mathcal D_u}q_f^{\G}(\lambda)
 \ge q_*^K.
\end{equation}
Fix, once and for all,
\begin{equation}\label{eq:s-choice-coupling}
 s:=\frac12q_*^K\in(0,1).
\end{equation}
Define
\begin{equation}\label{eq:alpha-thinning}
 \alpha_u
 :=\one_{A_u}\,
   \one\!_{\left\{W_u\le\frac{s}{q_u}\right\}}.
\end{equation}
The quotient $s/q_u$ belongs to $(0,1/2]$ by
\eqref{eq:qu-lower-bound}.

If $\alpha_u=1$, then every edge of the spanning tree $\mathcal T_u$ is
certified open and every selected exit $g_{u,v}$ is certified open.
Consequently, every vertex of $Z_u$ is connected to $x_u$; and for each
$v\in S_{r+1}^{\Hh}(u),$ there is a certified-open path from $x_u$ to a vertex of
$\pi^{-1}(v)$. We therefore
$$\mbox{add all vertices of $Z_u$ and all selected outer
endpoints to $C^\uparrow_{\mathrm{out}}$, and add $S_{r+1}^{\Hh}(u)$ to $C_{\mathrm{out}}$.}$$
The already reached part above $B_r^{\Hh}(u)$ causes no
extra base addition $\left(\text{note $B_r^{\Hh}(u)\subseteq C_{\mathrm{fr}}$}\right).$

The sets $\mathcal D_u$ are pairwise disjoint. Processing an earlier
centre may certify additional edges and so reduce the requests at a later
centre; it does not change the centre list, its eligibility data or the
snapshot from which $x_u$ is chosen. Freshness continues to hold because
reserved indices for different centres are distinct. At the end of all
tests in this phase, set
\[
 C:=C_{\mathrm{out}},\ C^\uparrow:=C^\uparrow_{\mathrm{out}},
\]
and start the next ordinary phase. Newly added vertices become eligible
for the centre list only after that ordinary phase, exactly as in
Section~\ref{sec:exploratory}. For an infinite centre list, the output
sets at the end of the phase are the increasing unions of the finite-step
outputs, and the history is the sigma-field generated by the earlier
histories. No random coordinate is sampled at a limit stage.


\subsection{The marginal laws}
\noindent
We now verify explicitly the two probabilistic facts needed for the coupling:
the final $\G$-configuration is a product measure, and the enhancement marks
are i.i.d.\,\,and independent of the base configuration.

\begin{lemma}[Adaptive deferred decisions]\label{lem:adaptive-deferred}
Let $(Z_i)_{i\in\mathcal I}$ be independent random variables indexed by a countable
set $\mathcal I$, independent of an initial sigma-field $\mathcal F_0$. An
exploration chooses each finite block of previously uninspected labels
measurably from its revealed history, and then inspects that block. The
history contains only the initial information, the choices already made,
and functions of variables already inspected. Conditional on that history
and on the labels chosen, a fresh block has its prescribed product law.
The exploration may consist of countably many phases, each containing
countably many such steps; at a limit stage its history is the sigma-field
generated by the earlier histories.

More generally, suppose that output labels in a countable set $\mathcal J$
are each assigned at most once. Each finite set of output labels is chosen
before its values are sampled, and conditional on the preceding history,
the new values have the product of prescribed laws $(\nu_j)_{j\in\mathcal J}$
at those labels. Assume that no unused filling variable is inspected by the
exploration. After filling the unassigned labels independently with their
prescribed laws, the completed output array has law
$\bigotimes_{j\in\mathcal J}\nu_j$.
\end{lemma}

\begin{proof}
For a finite number of steps, condition successively on the previously
chosen labels and on the revealed values. On the event that the next block
is a specified finite set $D$, the decision to inspect $D$ has used no
variable with label in $D$. Independence of the source variables therefore
gives their original product conditional law on this event. Summing over
the countably many possible finite sets $D$ proves the assertion for a
random block. The conditional law may depend on $D$ when the coordinate
distributions are not identical.

The precise unused-coordinate identity, including phase limits, is stated as follows.
Write $\mu_i$ for the law of $Z_i$, let $\mathcal F_t$ denote the history
at a stage $t$ of the fixed countable well-order, and let $R_t$ be the set
of source labels already inspected by that stage. For a deterministic
finite $D\subseteq\mathcal I$, set
\[
 U_D(t):=\{D\cap R_t=\varnothing\}.
\]
For every $A\in\mathcal F_t$ and bounded measurable $g$ on the coordinates
in $D$, we have that
\begin{eqnarray}\label{eq:unused-coordinate}
 \E\!\left[\one_A\one_{U_D(t)}g(Z_D)\right]
 =\Pp\bigl(A\cap U_D(t)\bigr)
   \int g\,{\rm d}\!\left(\bigotimes_{i\in D}\mu_i\right).
\end{eqnarray}
Equivalently, conditional on $\mathcal F_t$, the original product law of
$Z_D:=(Z_i)_{i\in D}$ is retained on $U_D(t)$. The indicator is essential: \eqref{eq:unused-coordinate}
is about labels that remain uninspected, not labels whose values have
already entered the history.

{\it Indeed}, at a successor step the next block is selected from $\mathcal F_t$.
On the event that the block meets $D$, the new unused indicator vanishes.
On its complement, only coordinates outside $D$ are inspected; the fresh
block product law, applied jointly with $D$, preserves the identity \eqref{eq:unused-coordinate}. The base case follows from independence of the source family
and $\mathcal F_0$.

At a countable limit stage $t$, take an increasing cofinal sequence
$t_n<t$. Then $R_t=\bigcup_nR_{t_n}$ and
$U_D(t_n)\downarrow U_D(t)$. For any event $A$ from an earlier history,
the identity \eqref{eq:unused-coordinate} holds at all sufficiently large $t_n$. Bounded
convergence gives it at $t$ for such $A$. A monotone-class argument then
extends it to every $A\in\mathcal F_t$, since
$\mathcal F_t=\sigma(\bigcup_n\mathcal F_{t_n})$. This proves \eqref{eq:unused-coordinate} throughout all phases and at their final union. The same proof
also allows conditioning on the initial sigma-field throughout.

For the output assertion, enlarge the probability space by a family of
mutually independent filling variables, with law $\nu_j$ at label $j$,
independent of the entire source family and the initial sigma-field.
At any stage form a completed array by retaining
the assigned values and using these filling variables at all other labels.
At a successor step, condition on the revealed history. The chosen
unassigned labels currently contain independent filling variables with
their prescribed laws, independently of the history and of all other
unused fillers. Replacing this block by the newly sampled block with the
same conditional product law leaves the unconditional law of the
completed array unchanged. The history may acquire additional information
from this sampling, provided unused fillers are never inspected.

At a countable limit stage choose an increasing cofinal sequence of earlier
stages. For every fixed finite set of output labels, the completed arrays
along this sequence eventually agree with the limiting completed array:
each assigned coordinate was assigned at some earlier stage, and every
unassigned coordinate retains the same filling variable. Bounded
convergence preserves every finite-dimensional product identity. Induction
over the stages, including these limit steps, proves the product law of the
final completed array.

Finally, the original unused source variables can serve as the final fillers
whenever their conditional product law given the completed exploration is
as above. Conditioning on finite subsets and using the unused-variable
identity proves that this replacement has the same distribution. All
statements also hold conditionally on $\mathcal F_0$ when the prescribed
output laws do not depend on that sigma-field.
\end{proof}

%
\vskip 2mm

Recall \eqref{eq:G-coordinate-1} and \eqref{eq:G-coordinate-2}. For a coordinate $(f,k)$ with $f\in E(\G)$ and $1\leq k\leq N$ never revealed during the ordinary exploration
and the enhancement phases, set
\begin{eqnarray}\label{eq:G-coordinate-3}
 \xi^{\G}(f,k):=\one_{\left\{V_{f,k}\le q_f^{\G}(\lambda)\right\}}\ \left(\mbox{the $\G$-auxiliary coordinate for $(f,k)$ never revealed}\right).
\end{eqnarray}

\begin{lemma}[Product law of the lifted configuration]\label{lem:product-law}
After all unrevealed $\G$-auxiliary coordinates are filled using their
independent variables $V_{f,k}$, the family
\[
 \left\{\xi^{\G}(f,k):\ f\in E(\G),\ 1\le k\le N\right\}
\]
is independent and satisfies
\[
 \Pp\left(\xi^{\G}(f,k)=1\right)=q_f^{\G}(\lambda).
\]
Consequently the physical configuration
\[
 \eta(f):=\max_{1\le k\le N}\xi^{\G}(f,k),\ f\in E(\G)
\]
has law $\Pp_\lambda^{\G}$.
\end{lemma}

\begin{proof}
There are two kinds of revealed $\G$-auxiliary coordinates. See \eqref{eq:G-coordinate-1} and \eqref{eq:G-coordinate-2}.
First, when a base edge
$e$ is ordinarily tested, its assigned lift $f=\varphi(e)$ is chosen before
the variables $U_{e,1},\dots,U_{e,N}$ are exposed.  By
Lemma~\ref{lem:fresh-coordinate}(v), no coordinate of $f$ has previously been
revealed.  Conditional on the entire past before this test,
\[
\left(\xi^{\G}(f,1),\dots,\xi^{\G}(f,N)\right)
\]
therefore has the product law
\[
 \bigotimes_{k=1}^N\Ber\left(q_f^{\G}(\lambda)\right),
\]
at the selected block of labels. This conditional product identity
is used with the labels fixed by the preceding history.


Second, every coordinate revealed during an enhancement phase has the form
$(f,\kappa_f(u))$ for an edge $f$ which is not an ordinary assigned lift.
Lemma~\ref{lem:fresh-coordinate}(iv) guarantees that the corresponding
$V_{f,\kappa_f(u)}$ has never previously been exposed.  Hence, conditionally
on the past, this $\G$-auxiliary coordinate $\xi^{\G}(f,\kappa_f(u))$ is an independent $\Ber\left(q_f^{\G}(\lambda)\right)$
variable.  Distinct enhancement reveals use distinct $V_{f,k}$-variables.

In both cases the random choice of the next label or finite block is measurable
with respect to the already revealed history, while the underlying uniforms to
be used next are fresh.  Lemma~\ref{lem:adaptive-deferred} therefore applies:
every newly revealed finite block of $\G$-auxiliary coordinates, conditionally on the
preceding history, has precisely the correct product distribution, even though
the labels are chosen adaptively. The output-completion part of the lemma, applied to the fixed label
set $E(\G)\times\{1,\ldots,N\}$, proves the required finite-dimensional
product identities at successive steps and at countable phase limits.
For this marginal argument we do not condition on the entire family
$\mathcal{U}$: an ordinary step still requires its unused $U_{e,\cdot}$-block to be fresh.

Finally, for a coordinate never revealed during the construction, the $\G$-auxiliary is defined by \eqref{eq:G-coordinate-3}.
For every fixed finite collection of still unassigned labels, the
associated unused $V_{f,k}$-variables have their original conditional product law.
Whether those labels remain unused is decided without inspecting those
variables. The final completion clause of
Lemma~\ref{lem:adaptive-deferred} therefore gives the claimed law for the
entire auxiliary family.


Taking the maximum over the $N$ independent copies of a fixed edge $f$ gives
\[
 \Pp(\eta(f)=1)
 =1-\left(1-q_f^{\G}(\lambda)\right)^N
 =1-e^{-\lambda w_{\G}(f)}.
\]
Independence between physical edges follows from independence of their
auxiliary coordinate blocks.
\end{proof}

\begin{lemma}[Law of the enhancement marks]\label{lem:alpha-law}
The family $(\alpha_u)_{u\in V(\Hh)}$ can be completed so that it consists of
independent $\Ber(s)$ random variables and is independent of the base
configuration $\omega$.
\end{lemma}

\begin{proof}
Fix one enhancement phase.  Let $\mathcal H_0$ be the sigma-field generated by
the complete history before that phase together with $\mathcal{U}$, the entire base auxiliary family. Process the unassigned vertices of $C_{\mathrm{fr}}$ in
the deterministic order fixed at the beginning of Section~\ref{sec:weighted-lifted-coupling}, say
$u_1,u_2,\cdots$.  For $j\ge1$, let $\mathcal H_{j-1}$ be the sigma-field
obtained from $\mathcal H_0$ by adjoining all auxiliary variables and
enhancement marks revealed while processing $u_1,\ldots,u_{j-1}$ in the
current phase.  Thus every random choice made before processing $u_j$ is
$\mathcal H_{j-1}$-measurable.


If $E\left(B_r^{\Hh}(u_j)\right)\not\subseteq\omega$, then
$\alpha_{u_j}=\one_{\{W_{u_j}\le s\}}$ and hence
\[
 \Pp(\alpha_{u_j}=1\mid\mathcal H_{j-1})=s.
\]
Suppose instead that $E\left(B_r^{\Hh}(u_j)\right)\subseteq\omega$. At the moment $u_j$ is
processed, the designated set $\mathcal D_{u_j}$ is
$\mathcal H_{j-1}$-measurable.  By Lemma~\ref{lem:adaptive-deferred} and the
fresh-coordinate bookkeeping of Lemma~\ref{lem:fresh-coordinate}, the
$V_{f,k}$-variables indexed by $\mathcal D_{u_j}$ retain their original independent
laws conditionally on $\mathcal H_{j-1}$.  Hence
\[
 q_{u_j}:=\Pp(A_{u_j}\mid\mathcal H_{j-1})
 =\prod_{(f,k)\in\mathcal D_{u_j}}q_f^{\G}(\lambda),
\]
and $W_{u_j}$ is independent of $\mathcal H_{j-1}$ and of these fresh
coordinates.  Therefore
\[
 \Pp(\alpha_{u_j}=1\mid\mathcal H_{j-1})
 =q_{u_j}\frac{s}{q_{u_j}}=s.
\]
The corresponding identity holds when the value zero is prescribed. Within the frozen
phase the labels $u_1,u_2,\ldots$ are fixed by $\mathcal H_0$. Consequently, for
any finite initial segment $u_1,\ldots,u_m$ and any $\varepsilon_1,\ldots,\varepsilon_m\in\{0,1\}$, the conditional chain rule
gives
\[
 \Pp\!\left(
   \alpha_{u_1}=\varepsilon_1,\ldots,
   \alpha_{u_m}=\varepsilon_m
   \ \middle|\ \mathcal H_0
 \right)
 =\prod_{j=1}^m
   s^{\varepsilon_j}(1-s)^{1-\varepsilon_j}.
\]
Labels in later phases can depend on earlier marks. Therefore we use the
output-completion argument of Lemma~\ref{lem:adaptive-deferred}, conditional
on $\mathcal{U}$, to pass from these sequential identities to the full array with
fixed labels. It applies because each mark is assigned once, its value has
conditional law $\Ber(s)$, and the unused $V_{f,k}$- and $W_u$-variables remain
fresh even after a countably infinite phase.

For a vertex whose mark is never assigned during exploration, set
$\alpha_u:=\one_{\{W_u\le s\}}$. No such $W_u$ was previously inspected.
The completion clause yields, for every deterministic finite
$F\subseteq V(\Hh)$ and every $\varepsilon\in\{0,1\}^{F}$,
\[
 \Pp\!\left(\alpha_u=\varepsilon_u\text{ for all }u\in F
       \,\middle|\,\mathcal{U}\right)
 =\prod_{u\in F}s^{\varepsilon_u}(1-s)^{1-\varepsilon_u}
 \ \text{almost surely}.
\]
This proves that the completed marks are i.i.d.\,\,$\Ber(s)$ and independent
of $\mathcal{U}$, and hence independent of its measurable function $\omega$.
\end{proof}

\begin{proof}[Proof of Proposition~\ref{prop:coupling}]
Choose $r$ as in Lemma~\ref{lem:geometry} and fix $\lambda_0>0$.  Choose $N$
by \eqref{eq:N-choice} and $s$ by \eqref{eq:s-choice-coupling}.  These
constants depend only on $r$, the degree bounds, $a$ and $\lambda_0$, and not
on $\lambda\ge\lambda_0$ or on the choice of the root lift $o'$.

By Lemmas~\ref{lem:product-law} and~\ref{lem:alpha-law}, the configuration
$\omega$ has law $\Pp_\lambda^{\Hh}$, the physical configuration $\eta$ has
law $\Pp_\lambda^{\G}$, and the marks $\alpha_u$ are i.i.d.\,\,$\Ber(s)$ and
independent of $\omega$.

During an ordinary phase, whenever a base edge $e$ is open,
\eqref{eq:coordinatewise-domination} implies that its assigned lift
$\varphi(e)$ is open.  Hence every vertex added to $C$ during an ordinary phase
has a lift in $C^\uparrow$ connected to $o'$ by $\eta$-open edges.  During an
enhancement phase, if $\alpha_u=1$, then $A_u$ occurs.  All edges of
$\mathcal T_u$ and all selected exits are therefore certified open, so every
vertex added above $B_r^{\Hh}(u)$ or $S_{r+1}^{\Hh}(u)$ is again connected to $o'$. Thus
throughout the construction
\[
 \pi\left(C^\uparrow\right)=C,
 \
 C^\uparrow\subseteq\cC_{\G}^{\lambda}\left(o'\right).
\]
The ordinary and enhancement phases agree exactly with the alternating
definition of $C_o(\omega,\alpha)$ from Section~\ref{sec:exploratory}.  Taking increasing unions
over all phases therefore gives
\[
 C_o(\omega,\alpha)=\pi\left(C^\uparrow_\infty\right)
 \subseteq\pi\Bigl(\cC_{\G}^{\lambda}\left(o'\right)\Bigr)
 \ \text{almost surely}.
\]
Since the left-hand side has the law $\cC_{\Hh}^{\lambda,s}(o)$, this proves
\[
 \cC_{\Hh}^{\lambda,s}(o)
 \preceq\pi\Bigl(\cC_{\G}^{\lambda}\left(o'\right)\Bigr)
\]
for every $\lambda\ge\lambda_0$ and every $o'\in\pi^{-1}(o)$.
\end{proof}

\section{Strictness of the weighted enhancement}\label{sec:strict}
\setcounter{equation}{0}
\noindent
We now prove Proposition~\ref{prop:enhancement}.  Throughout this section,
let $\Delta<\infty$ be a degree bound for $\Hh$ and fix constants
$0<a\le b<\infty$ such that $a\le w_{\Hh}(e)\le b$ for every
$e\in E(\Hh)$, as in the proposition.  The deterministic surgery follows the
local strategy of the homogeneous argument, but we spell it out in
a form that records the first entrance, the attaching edge, and the resumption
of the exploration.  The analytic modification is the weighted differential
formula together with uniform finite-energy constants.

For simplicity, write $\Pp_{\lambda,s}:=\Pp_{\lambda,s}^{\Hh}.$ For $L\ge1$, define
\[
  \cE_L:=\left\{C_o(\omega,\alpha)\cap S_L^{\Hh}(o)\ne\varnothing\right\},
  \
  \theta_L(\lambda,s):=\Pp_{\lambda,s}(\cE_L).
\]
The event $\cE_L$ depends on finitely many coordinates. More explicitly, before the
first hit of $S_L^{\Hh}(o)$, Lemma~\ref{lem:enhanced-connected} shows that all reached
vertices, and hence all reached enhancement centres, lie in $B_{L-1}^{\Hh}(o)$.
Testing whether such a centre has a fully open $r$-ball uses only edge
coordinates in $E\left(B_{L+r}^{\Hh}(o)\right)$, while an ordinary step that first reaches
$S_L^{\Hh}(o)$ uses only edges in $E\left(B_L^{\Hh}(o)\right)$. Hence
\begin{equation}\label{eq:event-measurable}
\begin{split}
&\text{$\cE_L$ is measurable with respect to the edge coordinates in $E\left(B_{L+r}^{\Hh}(o)\right)$ together with}\\
&\text{enhancement
marks in $B_{L-1}^{\Hh}(o)$.}
\end{split}
\end{equation}
Note these sets are finite.

By Lemma~\ref{lem:enhanced-connected}, the events $\cE_L$ decrease with $L$ and
their intersection is precisely the event that the exploratory enhanced
cluster is infinite. This implies
\begin{equation}\label{eq:theta-infinite-cluster}
  \theta(\lambda,s):=\Pp_{\lambda,s}
  \bigl(|C_o(\omega,\alpha)|=\infty\bigr)
  =\lim_{L\to\infty}\theta_L(\lambda,s).
\end{equation}

An edge $e$ is \emph{edge-pivotal} for an increasing event if forcing it open
makes the event occur and forcing it closed makes the event fail.  A vertex $x$
is \emph{enhancement-pivotal} if the analogous statement holds when $\alpha_x$
is forced to one or zero.

\begin{lemma}[Inhomogeneous Margulis--Russo formula]
\label{lem:inhomogeneous-russo}
Let $\cE$ be an increasing event depending on finite sets of independent edge
coordinates and enhancement coordinates.  Suppose that edge $e$ is open with
probability $p_e(t)$ with $p_e$ being continuously differentiable, and that each
enhancement coordinate has parameter $s$.  Then
\begin{align}
 \frac{\partial}{\partial t}\Pp_{t,s}(\cE)
 &=\sum_e p'_e(t)\,
   \Pp_{t,s}(e\text{ is edge-pivotal for $\cE$}),
 \label{eq:general-russo-t}\\
 \frac{\partial}{\partial s}\Pp_{t,s}(\cE)
 &=\sum_x
   \Pp_{t,s}(x\text{ is enhancement-pivotal for $\cE$}).
 \label{eq:general-russo-s}
\end{align}
\end{lemma}

\begin{proof}
Condition on all coordinates except edge $e$.  The conditional probability
of $\cE$ is affine in $p_e$ and its derivative with respect to $p_e$ is the
indicator averaged over the other coordinates that $e$ is pivotal.  The chain
rule gives the contribution $p'_e(t)\,\Pp_{t,s}(e\text{ is pivotal for $\cE$})$. Summing over
the finitely many relevant edges proves \eqref{eq:general-russo-t}.  The same
conditioning argument for a Bernoulli enhancement coordinate gives
\eqref{eq:general-russo-s}.  This is the usual finite-dimensional proof of
Russo's formula \cite{Ru1981, Gr1999}.
\end{proof}

\begin{lemma}[Deterministic pivotal conversion]\label{lem:surgery}
Let $r\ge1$.  There are constants
\[
  R=3r+2,
  \
  L_0=2r+2
\]
such that the following holds.  If $L\ge L_0$ and $e$ is edge-pivotal for
$\cE_L$ in $(\omega,\alpha)$, then one can change $(\omega,\alpha)$ only
inside $B_R^{\Hh}(e)$ and obtain a configuration $(\omega',\alpha')$ for which some
$z\in B_R^{\Hh}(e)$ is enhancement-pivotal for $\cE_L$.  Moreover, the number of
edge and enhancement coordinates on which the surgery acts is bounded above by a
constant depending only on $r$ and the degree bound.
\end{lemma}

\begin{proof} This lemma is \cite[Lemma 2.4]{Ma-Se2019}. For completeness, a detailed proof is given in \ref{app-2} according to that of \cite[Lemma 2.4]{Ma-Se2019}.
\end{proof}

\begin{lemma}[Uniform pivotal comparison]\label{lem:pivotal-comparison}
Fix compact intervals
\[
  I=[\lambda_-,\lambda_+]\subset(0,\infty),
  \ J=[\delta,1-\delta]\subset(0,1).
\]
There exists $c_0>0$ such that for all $L\ge L_0$, all $e\in E(\Hh)$ and all
$(\lambda,s)\in I\times J$,
\begin{equation}\label{eq:local-pivotal-comparison}
  \sum_{z\in B_R^{\Hh}(e)}
  \Pp_{\lambda,s}(z\text{ is enhancement-pivotal for $\cE_L$})
  \ge c_0\,
  \Pp_{\lambda,s}(e\text{ is edge-pivotal for $\cE_L$}).
\end{equation}
\end{lemma}

\begin{proof}
By ellipticity, every edge probability on $I$ lies in the common interval
\[
  \left[1-e^{-\lambda_-a},\ 1-e^{-\lambda_+b}\right]
  \subset(0,1).
\]
The enhancement probabilities lie in $J$, and put
\[
 \rho:=\min\left\{1-e^{-\lambda_-a},e^{-\lambda_+b},\delta\right\}>0.
\]
Thus, conditionally on all other coordinates, either state of every relevant
edge or enhancement coordinate has probability at least $\rho$.

Fix $e$ and let $D_e$ be a deterministic local set of binary coordinates large
enough to contain every coordinate that can be changed by the surgery of
Lemma~\ref{lem:surgery}, including the possible normalization of the actual
state of $e$. Bounded degree and the fixed radius $R$ give
$|D_e|\le M_{\mathrm{loc}}$ uniformly in $e$ and $L$. Let $\mathcal P_e$ be
the event that $e$ is edge-pivotal for $\cE_L$. Using all the fixed tie-breaking rules from
Lemma~\ref{lem:surgery}, define a deterministic map
\[
  \Phi_e:\ \mathcal P_e\longrightarrow
  \bigcup_{z\in B_R^{\Hh}(e)}(\mathcal Q_z\times\{z\}),
\]
where $\mathcal Q_z$ is the event that $z$ is enhancement-pivotal for $\cE_L$: first
normalize $e$ to the open state, then perform the local surgery, and record
both the output configuration and the selected pivotal centre $z$. By Lemma~\ref{lem:surgery}, the output indeed belongs to $\mathcal Q_z$, and input and output agree
outside $D_e$.

Condition on all coordinates outside $D_e$. If $\xi$ is an input local pattern
and $\eta$ is its output local pattern, then the product structure and the
finite-energy bound imply
\begin{equation}\label{eq:local-pattern-ratio}
 \frac{\Pp_{\lambda,s}(\eta\mid D_e^c)}
      {\Pp_{\lambda,s}(\xi\mid D_e^c)}
 \ge \left(\frac{\rho}{1-\rho}\right)^{M_{\mathrm{loc}}}.
\end{equation}
For a fixed pair $(\eta,z),$ there are at most $2^{M_{\mathrm{loc}}}$ possible
input local patterns, because $D_e$ contains at most
$M_{\mathrm{loc}}$ binary coordinates. Therefore, still conditionally on the
outside configuration,
\[
 \Pp_{\lambda,s}\left(\mathcal P_e\mid D_e^c\right)
 \le
 2^{M_{\mathrm{loc}}}
 \left(\frac{1-\rho}{\rho}\right)^{M_{\mathrm{loc}}}
 \sum_{z\in B_R^{\Hh}(e)}
 \Pp_{\lambda,s}\left(\mathcal Q_z\mid D_e^c\right).
\]
Integrating over the outside coordinates yields
\eqref{eq:local-pivotal-comparison}, where
\[
 c_0=
 2^{-M_{\mathrm{loc}}}
 \left(\frac{1-\rho}{\rho}\right)^{-M_{\mathrm{loc}}}>0.
\]
This constant depends only on $a,\,b,\,I,\,J,\,R$ and the degree bound, and not on
$L$, $e$, $\lambda$ or $s$ in the indicated ranges.
\end{proof}

\begin{proof}[Proof of Proposition~\ref{prop:enhancement}]
By monotonicity in $s$, it suffices to treat $s\in(0,1)$.  For a finite increasing
event, Lemma~\ref{lem:inhomogeneous-russo}, applied with
$p_e(\lambda)=1-\exp\{-\lambda w_{\Hh}(e)\}\ (e\in E(\Hh))$, gives
\begin{equation}\label{eq:weighted-russo}
  \frac{\partial}{\partial\lambda}\theta_L(\lambda,s)
  =\sum_{e\in E(\Hh)}
    w_{\Hh}(e)\exp\{-\lambda w_{\Hh}(e)\}\,
    \Pp_{\lambda,s}(e\text{ is edge-pivotal for $\cE_L$}),
\end{equation}
and
\begin{equation}\label{eq:s-russo}
  \frac{\partial}{\partial s}\theta_L(\lambda,s)
  =\sum_{z\in V(\Hh)}
    \Pp_{\lambda,s}(z\text{ is enhancement-pivotal for $\cE_L$}).
\end{equation}

Each vertex belongs to $B_R^{\Hh}(e)$ for at most
\[
  C_R:=\sup_{x\in V(\Hh)}\left|E\left(B_{R+1}^{\Hh}(x)\right)\right|<\infty
\]
different edges $e$.  Summing \eqref{eq:local-pivotal-comparison} over $e$ and
using \eqref{eq:s-russo} yields
\begin{equation}\label{eq:pivital-comparison}
  \frac{\partial}{\partial s}\theta_L(\lambda,s)
  \ge\frac{c_0}{C_R}
      \sum_{e\in E(\Hh)}\Pp_{\lambda,s}(e\text{ is edge-pivotal for $\cE_L$}).
\end{equation}
Since $w_{\Hh}(e)\exp\{-\lambda w_{\Hh}(e)\}\le b$, by \eqref{eq:weighted-russo} and \eqref{eq:pivital-comparison},
\begin{equation}\label{eq:differential-inequality}
  \frac{\partial}{\partial s}\theta_L(\lambda,s)
  \ge c\frac{\partial}{\partial\lambda}\theta_L(\lambda,s),
  \  c:=\frac{c_0}{C_Rb}>0,
\end{equation}
on every compact rectangle $I\times J\subset(0,\infty)\times(0,1)$, for all
$L\ge L_0$.

Bounded degree and the upper weight bound imply
\begin{eqnarray}\label{eq:positivity}
\lambda_c:=\lambda_c(\Hh,w_{\Hh})>0.
\end{eqnarray}
{\it Indeed}, the number of self-avoiding paths of length $n$
from a fixed vertex is at most $\Delta(\Delta-1)^{n-1}$, while each such path is
open with probability at most $(\lambda b)^n$ (note $p_e(\lambda)\le\lambda b$). A union bound followed by
$n\to\infty$ excludes an infinite open path whenever $\lambda b(\Delta-1)<1$.

Let $\bar s=\min\{s,1/2\}$ and choose $0<\delta<\bar s$.  Apply
\eqref{eq:differential-inequality} on a compact rectangle around $\lambda_c$ and
with enhancement coordinate between $\delta$ and $\bar s$.  Choose $\varepsilon>0$
small enough that
\[
  s_1:=\delta+\frac{2\varepsilon}{c}\le\bar s
\]
and the line segment $\left(\lambda(t),\,s(t)\right)_{0\le t\le 1}$ joining
\[
  (\lambda_c+\varepsilon,\delta)\
  \text{to}\
  (\lambda_c-\varepsilon,s_1)
\]
lies inside the rectangle. Notice that
$$\lambda'(t)=-c\,s'(t),\ s'(t)=\frac{2\varepsilon}{c}>0,\ t\in [0,1].$$
By the chain rule and \eqref{eq:differential-inequality},
\[
  \frac{{\rm d}}{{\rm d}t}\theta_L(\lambda(t),s(t))\ge0.
\]
Passing to $L\to\infty$ and using monotonicity in $s$ gives
\begin{align*}
  \theta(\lambda_c-\varepsilon,s)
  \ge\theta(\lambda_c-\varepsilon,s_1)
  \ge\theta(\lambda_c+\varepsilon,\delta)
  \ge\theta(\lambda_c+\varepsilon,0)>0.
\end{align*}
Thus $\lambda_s:=\lambda_c-\varepsilon<\lambda_c$ is supercritical for the
enhanced process.  Monotonicity in $\lambda$ completes the proof.
\end{proof}

\section{The uniform connection intensity}\label{sec:uniqueness}
\setcounter{equation}{0}
\noindent
This section proves Theorem~\ref{thm:conn-strict} (and Corollary \ref{cor:pu-conditional}). For vertex sets $A,B\subseteq V(\G)\ (\mbox{resp.}\ V(\Hh))$, we write
$A\leftrightarrow B$ when an open path meets both sets. And for a finite nonempty $A\subseteq V(\Hh)$, define $C_A(\omega,\alpha)$ by the
same alternating exploration as in Section~\ref{sec:exploratory}, but start from $C_0=A$. Write $C_A(\omega,\alpha):=C_A^{\lambda,s}(\omega,\alpha)$ for this random set under $\Pp_{\lambda,s}=\Pp_{\lambda,s}^{\Hh}$, and
\[
 A\xRightarrow{\lambda,s}B
 \ \text{for the event}\  C_A(\omega,\alpha)\cap B\ne\varnothing.
\]


\begin{lemma}[\text{\cite[Lemma~4.2]{Ma-Se2019}}]
\label{lem:automatic-tameness}
Let $\G$ and $\Hh$ be quasi-transitive graphs.  If a non-injective weak
covering map $\pi:\ V(\G)\to V(\Hh)$ has the disjoint tree-lifting property, then
its fibres are tame.
\end{lemma}

\begin{proposition}[Finite-start lifted enhancement]\label{prop:finite-start}
Under the assumptions of Theorem~\ref{thm:conn-strict}, there is an $r\ge1$ such
that for any $\lambda_0>0,$ there is an $s\in(0,1)$ with
the following property.  For any $\lambda\ge\lambda_0$ and any finite
nonempty $A_0\subseteq V(\G)$,
\[
 C_{\pi(A_0)}^{\lambda,s}\preceq
 \pi\bigl(\cC_{\G}^{\lambda}(A_0)\bigr).
\]
Consequently, for all finite $A,B\subseteq V(\Hh)$,
\begin{equation}\label{eq:enhanced-connection-lift}
 \Pp_\lambda^{\G}\left[\pi^{-1}(A)\leftrightarrow\pi^{-1}(B)\right]
 \ge \Pp_{\lambda,s}^{\Hh}[A\xRightarrow{}B].
\end{equation}
\end{proposition}

\begin{proof}
Lemma~\ref{lem:automatic-tameness} supplies tame fibres.  Run the construction
of Proposition~\ref{prop:coupling} with initial lifted set $A_0$ and initial
base set $\pi(A_0)$ in place of $(\{o'\},\{o\})$.  The exploration invariants
are interpreted componentwise: every vertex of the lifted explored set is
connected by certified-open edges to at least one vertex of $A_0$, and the
projection of the lifted explored set equals the base explored set.

No coordinate bookkeeping changes are required.  The reservation maps
$\kappa_f$ of \eqref{eq:reservation-map} are fixed independently of the
starting set; Lemma~\ref{lem:fresh-coordinate} applies verbatim because its
proof uses only exhaustion of the preceding ordinary phase and the local
geometry around an eligible enhancement centre.  Likewise, the proof of
Lemma~\ref{lem:product-law} is unchanged: every ordinary block is selected
before its shared uniforms are exposed, and every enhancement reveal uses a
fresh reserved $V_{f,k}$-coordinate.  Finally, the sequential thinning argument of
Lemma~\ref{lem:alpha-law} is unchanged and gives the same i.i.d.\,\,$\Ber(s)$
enhancement field, independent of the base configuration.  In particular,
the constants $r,\,N,\,K,\,q_*$ and $s$ depend only on the geometric and ellipticity
parameters appearing in Proposition~\ref{prop:coupling}, not on
$|A_0|$.

The deterministic inclusion produced by the construction is therefore
\[
 C_{\pi(A_0)}^{\lambda,s}
 \subseteq
 \pi\bigl(\cC_{\G}^{\lambda}(A_0)\bigr)
 \ \text{under the coupling}.
\]
Bounded fibres and local finiteness make $\pi^{-1}(A)$ finite when $A$ is
finite.  Taking $A_0=\pi^{-1}(A)$ gives
\eqref{eq:enhanced-connection-lift}.
\end{proof}

Recall from \cite[p.\,4129]{Ma-Se2019}, a set $B\subseteq V(\Hh)$ is \emph{$r$-nice} if
\begin{enumerate}
\item[] it is finite and nonempty, and for any $x\notin B$, there is a vertex $z\in V(\Hh)$ such that
$$x\in B_r^{\Hh}(z),\ B_r^{\Hh}(z)\cap B=\varnothing.$$
\end{enumerate}

For $D\subseteq V(\Hh)$ and $x\in V(\Hh)$, note $d_{\Hh}(x,D)=\inf\limits_{v\in D}d_{\Hh}(x,v).$
For an edge $e=\{x,y\}\in E(\Hh),$ recall from \eqref{eq:ball}-\eqref{eq:e-ball},  $B_R^{\Hh}(e)=B_R^{\Hh}(x)\cup B_R^{\Hh}(y);$
and set
\[
 d_{\Hh}(e,D):=\min\left\{d_{\Hh}(x,D),\,d_{\Hh}(y,D)\right\}.
\]
For a vertex set $D\subseteq V(\Hh)$, let
\[
 \partial_E D:=\left\{\{x,y\}\in E(\Hh):\ x\in D,\ y\notin D\right\}.
\]

Fix $o\in V(\Hh)$.  For $L\ge1$ and a configuration
$(\omega,\alpha)$, let $\left(\omega^{[L]},\alpha^{[L]}\right)$ be obtained by forcing
closed every edge not entirely contained in $B_L^{\Hh}(o)$ and forcing to zero every
enhancement coordinate outside $B_L^{\Hh}(o)$.  For finite nonempty
$A,B\subseteq B_L^{\Hh}(o),$ define
\begin{equation}\label{eq:finite-volume-event}
 \cE_L^{A,B}
 :=\left\{C_A\left(\omega^{[L]},\alpha^{[L]}\right)\cap B\ne\varnothing\right\}.
\end{equation}
This event depends on finitely many independent coordinates.

\begin{lemma}[Finite-set deterministic pivotal conversion]
\label{lem:finite-set-pivotal-conversion}
Fix $r\ge1$ and assume that $\Hh$ has maximum degree $\Delta<\infty$. There are constants
\[
 R=3r+2,
 \
 M=M(r,\Delta)<\infty
\]
with the following property.  Let $A,B\subseteq V(\Hh)$ be finite and
nonempty, suppose that $B$ is $r$-nice and $d_{\Hh}(A,B)>3r$, and choose $L$ so
large that
\begin{equation}\label{eq:finite-volume-separation}
 A\cup B\subseteq B_L^{\Hh}(o),
 \
 d_{\Hh}\left(A\cup B,S_L^{\Hh}(o)\right)>3r.
\end{equation}
If an edge $e$ is edge-pivotal for $\cE_L^{A,B}$ in a configuration
$(\omega,\alpha)$, then by modifying only coordinates in $B_R^{\Hh}(e)$ one obtains
a configuration $(\omega',\alpha')$ and a vertex $z\in B_R^{\Hh}(e)$ such that $z$
is enhancement-pivotal for $\cE_L^{A,B}$.  The total number of edge and vertex
coordinates in the region on which the surgery may act is at most $M$, and the
construction may be chosen deterministically.
\end{lemma}

\begin{proof} This lemma is \cite[Lemma 3.3]{Ma-Se2019}. For completeness, a detailed proof is given in \ref{app-3} according to that of \cite[Lemma 3.3]{Ma-Se2019}.
\end{proof}

\begin{proposition}[Enhancement beats an intensity increment]
\label{prop:finite-set-enhancement}
Assume that $\Hh$ has maximum degree $\Delta<\infty$ and that
$0<a\le w_{\Hh}(e)\le b<\infty$ for every $e\in E(\Hh)$.  Fix
$r\ge1$, $s\in(0,1]$, and a compact interval
$I=[\lambda_-,\lambda_+]\subset(0,\infty)$.  There exists
\[
 \delta=\delta(r,s,I,a,b,\Delta)>0
\]
such that whenever $[\lambda,\lambda+\delta]\subseteq I$ and $A,B$ are finite
nonempty sets with $B$ $r$-nice and $d_{\Hh}(A,B)>3r$,
\begin{equation}\label{eq:enhancement-beats-increment}
 \Pp_{\lambda,s}^{\Hh}[A\xRightarrow{}B]
 \ge \Pp_{\lambda+\delta}^{\Hh}[A\leftrightarrow B].
\end{equation}
The constant $\delta$ is independent of $A$ and $B$.
\end{proposition}

\begin{proof}
Put $\bar s:=\min\{s,1/2\}$ and
\[
 J:=[\bar s/2,\bar s]\subset(0,1).
\]
Since the enhanced connection probability is increasing in the enhancement
parameter, it suffices to prove the comparison with $\bar s$ in place of $s$.

Choose $L$ sufficiently large that \eqref{eq:finite-volume-separation} holds.  Write
\[
 \theta_L^{A,B}(\mu,t):=\Pp_{\mu,t}^{\Hh}\left(\cE_L^{A,B}\right).
\]
We first derive a pivotal comparison uniform in $A,B,L$.

For $(\mu,t)\in I\times J$, uniform ellipticity gives
\[
 1-\exp\{-\lambda_-a\}\le p_e^{\Hh}(\mu)
 \le 1-\exp\{-\lambda_+b\}.
\]
Thus every edge coordinate, in either state, has probability at least
\[
 \rho_E:=\min\left\{1-\exp\{-\lambda_-a\},\exp\{-\lambda_+b\}\right\}>0,
\]
and every enhancement coordinate, in either state, has probability at least
\[
 \rho_V:=\min\{\bar s/2,1-\bar s\}>0.
\]
Put $\rho:=\min\{\rho_E,\rho_V\}>0$.

By Lemma~\ref{lem:finite-set-pivotal-conversion}, including its normalization
of the actual state of $e$, there is a deterministic map from every
edge-pivotal input to a pair consisting of a locally modified output and a
vertex $z\in B_R^{\Hh}(e)$ which is enhancement-pivotal in that output. Let $D_e$
be a deterministic coordinate set containing the entire possible surgery
region. We may take $|D_e|\le M$ uniformly in $A,B,L$ and $e$. Conditional on
all coordinates outside $D_e$, the input and output are two binary patterns on
$D_e$, and the finite-energy bound gives
\[
 \frac{\Pp_{\mu,t}^{\Hh}\left(\text{output pattern}\mid D_e^c\right)}
      {\Pp_{\mu,t}^{\Hh}\left(\text{input pattern}\mid D_e^c\right)}
 \ge\left(\frac{\rho}{1-\rho}\right)^M.
\]
For a fixed output--centre pair, there are at most $2^M$ input patterns. Thus,
after summing over outputs and over $z\in B_R^{\Hh}(e)$ and then integrating over
the outside coordinates, one obtains a constant
\[
 c_0=c_0(r,\Delta,I,J,a,b)>0
\]
$\left(\text{for example $c_0=2^{-M}((1-\rho)/\rho)^{-M}$}\right)$ such that for every relevant
edge $e$ and every $(\mu,t)\in I\times J$,
\begin{equation}\label{eq:finite-set-pivotal-comparison}
 \sum_{z\in B_R^{\Hh}(e)}
 \Pp_{\mu,t}^{\Hh}
 \Bigl(z\text{ is enhancement-pivotal for }\cE_L^{A,B}\Bigr)
 \ge
 c_0\,\Pp_{\mu,t}^{\Hh}
 \Bigl(e\text{ is edge-pivotal for }\cE_L^{A,B}\Bigr).
\end{equation}
The constant $c_0$ is independent of $A,\, B,\, L,\, e,\, \mu$ and $t$ in the indicated
ranges.

Let
\[
 C_R:=\sup_{z\in V(\Hh)}
 \Bigl|\Bigl\{e\in E(\Hh):\ z\in B_R^{\Hh}(e)\Bigr\}\Bigr|<\infty.
\]
Summing \eqref{eq:finite-set-pivotal-comparison} over $e$ gives
\begin{equation}\label{eq:sum-finite-set-pivotal}
 \sum_z\Pp_{\mu,t}^{\Hh}\left(z\text{ is enhancement-pivotal for $\cE_L^{A,B}$}\right)
 \ge
 \frac{c_0}{C_R}
 \sum_e\Pp_{\mu,t}^{\Hh}\left(e\text{ is edge-pivotal for $\cE_L^{A,B}$}\right).
\end{equation}
Since $\cE_L^{A,B}$ is finite-dimensional, Lemma~\ref{lem:inhomogeneous-russo}
gives
\begin{align}
 \frac{\partial}{\partial t}\theta_L^{A,B}(\mu,t)
 &=\sum_z\Pp_{\mu,t}^{\Hh}\left(z\text{ is enhancement-pivotal for $\cE_L^{A,B}$}\right),
 \label{eq:finite-set-russo-s}\\
 \frac{\partial}{\partial\mu}\theta_L^{A,B}(\mu,t)
 &=\sum_e w_{\Hh}(e)\exp\{-\mu w_{\Hh}(e)\} \,
   \Pp_{\mu,t}^{\Hh}\left(e\text{ is edge-pivotal for $\cE_L^{A,B}$}\right).
 \label{eq:finite-set-russo-lambda}
\end{align}
Because $w_{\Hh}(e)\exp\{-\mu w_{\Hh}(e)\}\le b$, we obtain
\begin{equation}\label{eq:finite-set-differential}
 \frac{\partial}{\partial t}\theta_L^{A,B}(\mu,t)
 \ge c\frac{\partial}{\partial\mu}\theta_L^{A,B}(\mu,t),
 \  c=\frac{c_0}{C_Rb}>0,
\end{equation}
uniformly in $A,\, B,\, L$ and $(\mu,t)\in I\times J$.

Choose
\begin{equation}\label{eq:finite-set-delta}
 0<\delta\le
 \min\left\{\frac{c\bar s}{2},\frac{\lambda_+-\lambda_-}{2}\right\}.
\end{equation}
Let $[\lambda,\lambda+\delta]\subseteq I$; and for $t\in[\bar s/2,\bar s]$, define
\[
 \mu(t):=\lambda+\delta
 -\frac{2\delta}{\bar s}\left(t-\frac{\bar s}{2}\right).
\]
Then $\mu(\bar s/2)=\lambda+\delta$, $\mu(\bar s)=\lambda$, and
$\mu'(t)=-2\delta/\bar s\ge-c$.  By
\eqref{eq:finite-set-differential},
\[
 \frac{{\rm d}}{{\rm d}t}\theta_L^{A,B}(\mu(t),t)
 =\partial_t\theta_L^{A,B}
  +\mu'(t)\partial_\mu\theta_L^{A,B}\ge0.
\]
Consequently
\begin{equation}\label{eq:integrated-finite-set}
 \theta_L^{A,B}(\lambda,\bar s)
 \ge\theta_L^{A,B}(\lambda+\delta,\bar s/2)
 \ge\theta_L^{A,B}(\lambda+\delta,0).
\end{equation}

As $L\to\infty$, the events $\cE_L^{A,B}$ increase to the full enhanced
connection event $\{A\xRightarrow{}B\}$. {\it Indeed}, if the enhanced exploration
reaches the finite set $B$, then the first reached vertex of $B$ has a finite
witness consisting of finitely many ordinary open paths, enhancement centres,
and finite $r$-balls, and this witness is contained in $B_L^{\Hh}(o)$ for all
sufficiently large $L$.  Monotone convergence in
\eqref{eq:integrated-finite-set} therefore yields
\[
 \Pp_{\lambda,\bar s}^{\Hh}[A\xRightarrow{}B]
 \ge\Pp_{\lambda+\delta,0}^{\Hh}[A\xRightarrow{}B]
 =\Pp_{\lambda+\delta}^{\Hh}[A\leftrightarrow B].
\]
Finally $\bar s\le s$ and monotonicity in the enhancement parameter imply \eqref{eq:enhancement-beats-increment}.
\end{proof}

\begin{proof}[Proof of Theorem~\ref{thm:conn-strict}]
When $\lambda_{\mathrm{conn}}(\Hh,w_{\Hh})=\infty$, the theorem is trivially true. Now set
$$\lambda_{\mathrm{conn}}
:=\lambda_{\mathrm{conn}}(\Hh,w_{\Hh})<\infty.$$
We first note that any parameter satisfying the uniform connection criterion
is at least critical for percolation.  Let $\mu$ belong to the defining set in
Definition~\ref{def:lambda-conn}.  Fix $x$ and $\ell$, and choose $y_n$ with
$d_{\Hh}(x,y_n)\to\infty$.  The event
$B_\ell^{\Hh}(x)\leftrightarrow B_\ell^{\Hh}(y_n)$ implies that the open cluster of
$B_\ell^{\Hh}(x)$ has radius at least $d_{\Hh}(x,y_n)-2\ell$. Hence
\[
 \Pp_\mu^{\Hh}\left[B_\ell^{\Hh}(x)\leftrightarrow\infty\right]
 \ge \limsup_{n\to\infty}
 \Pp_\mu^{\Hh}\left[B_\ell^{\Hh}(x)\leftrightarrow B_\ell^{\Hh}(y_n)\right].
\]
Since the uniform connection lower bound tends to $1$ as $\ell\to\infty$,
some finite ball has positive probability of connection to infinity.  At least
one vertex in that ball then percolates with positive probability.  By the
root-independence established after Definition~\ref{def:lambda-c},
$\mu\ge\lambda_c(\Hh,w_{\Hh})$. By taking the infimum,
\[
 \lambda_{\mathrm{conn}}(\Hh,w_{\Hh})
 \ge\lambda_c(\Hh,w_{\Hh})>0\ (\mbox{noting}\ \eqref{eq:positivity}).
\]

Choose $0<\lambda_0<\lambda_{\mathrm{conn}}$.
Proposition~\ref{prop:finite-start} supplies $r$ and $s\in(0,1)$ for all
intensities at least $\lambda_0$.  Apply
Proposition~\ref{prop:finite-set-enhancement} on a compact interval $I$
containing $\lambda_{\mathrm{conn}}$ in its interior and obtain $\delta>0$.
Choose $\lambda<\lambda_{\mathrm{conn}}$ sufficiently close to
$\lambda_{\mathrm{conn}}$ that
\[
 \lambda_0<\lambda<\lambda_{\mathrm{conn}}<q:=\lambda+\delta,
 \ [\lambda,q]\subseteq I.
\]
Then for every finite $A,B\subseteq V(\Hh)$ with $B$ $r$-nice and
$d_{\Hh}(A,B)>3r$,
\begin{equation}\label{eq:strict-connection-chain}
 \Pp_\lambda^{\G}\left[\pi^{-1}(A)\leftrightarrow\pi^{-1}(B)\right]
 \ge \Pp_{\lambda,s}^{\Hh}[A\xRightarrow{}B]
 \ge \Pp_q^{\Hh}[A\leftrightarrow B].
\end{equation}

Let $K$ be a fibre-diameter bound.  Fix $\ell\ge1$ and $x,y\in V(\G)$ with
\[
 d_{\G}(x,y)>L(\ell):=2\ell+4r+K.
\]
Put $u=\pi(x)$, $v=\pi(y)$, $A=B_\ell^{\Hh}(u)$, and
\begin{equation}\label{eq:nice-annulus-set}
 B:=V(\Hh)\setminus
 \bigcup_{z:\,d_{\Hh}(z,v)>\ell+r}B_r^{\Hh}(z).
\end{equation}
Clearly $B$ is $r$-nice. Moreover,
\[
 B_\ell^{\Hh}(v)\subseteq B\subseteq B_{\ell+r}^{\Hh}(v),\ \mbox{which implies}\ B\not=\varnothing\ \mbox{and}\ \vert B\vert<\infty.
\]
{\it Indeed}, for the first inclusion, if $x'\in B_\ell^{\Hh}(v)$ and $x'\in B_r^{\Hh}(z)$, then
$d_{\Hh}(z,v)\le r+\ell$, and thus $x'\in B$. For the second, if $d_{\Hh}(x',v)>\ell+r$, take $z=x'$, then $x'\in B_r^{\Hh}(z)$ belongs to
the removed union, so $x'\notin B$.

By Lemma~\ref{lem:metric-comparison},
$d_{\Hh}(u,v)>2\ell+4r$, so $d_{\Hh}(A,B)>3r$; and
\[
 \pi^{-1}(A)\subseteq B_{\ell+K}^{\G}(x)\subseteq B_{L(\ell)}^{\G}(x),
 \
 \pi^{-1}(B)\subseteq B_{\ell+r+K}^{\G}(y)\subseteq B_{L(\ell)}^{\G}(y).
\]
By \eqref{eq:strict-connection-chain} and the inclusion
$B_\ell^{\Hh}(v)\subseteq B$,
\begin{equation}\label{eq:strict-ball-comparison}
 \Pp_\lambda^{\G}\left[B_{L(\ell)}^{\G}(x)\leftrightarrow B_{L(\ell)}^{\G}(y)\right]
 \ge
 \Pp_q^{\Hh}\left[B_\ell^{\Hh}(\pi(x))\leftrightarrow B_\ell^{\Hh}(\pi(y))\right].
\end{equation}
Note when $d_{\G}\left(x',y'\right)\le L(\ell)$, $B_{L(\ell)}^{\G}\left(x'\right)\cap B_{L(\ell)}^{\G}\left(y'\right)\not=\varnothing$ and $\Pp_\lambda^{\G}\left[B_{L(\ell)}^{\G}\left(x'\right)\leftrightarrow B_{L(\ell)}^{\G}\left(y'\right)\right]=1.$
By taking the infimum over $x,y$ in
\eqref{eq:strict-ball-comparison},
\[
 a_{L(\ell)}^{\G}(\lambda)
 \ge a_\ell^{\Hh}(q).
\]
Since $q>\lambda_{\mathrm{conn}}(\Hh,w_{\Hh})$, the right-hand side tends to one.
Hence $a_{L(\ell)}^{\G}(\lambda)\to1$ along the diverging sequence
$L(\ell)\to\infty$. By that the function $R\mapsto a_R^{\G}(\lambda)$ is
increasing and bounded above by $1$,
\[
 \lim_{R\to\infty}a_R^{\G}(\lambda)=1.
\]
Thus
\[
 \lambda_{\mathrm{conn}}(\G,w_{\G})
 \le\lambda<\lambda_{\mathrm{conn}}=\lambda_{\mathrm{conn}}(\Hh,w_{\Hh}).
\]
\end{proof}

\begin{proof}[Proof of Corollary~\ref{cor:pu-conditional}]
By Lemma~\ref{lem:connection-criterion}(ii), property $\mathsf{Stab}$ identifies
$\lambda_u$ with $\lambda_{\mathrm{conn}}$ on each weighted graph. By Theorems~\ref{thm:conn-monotone} and
\ref{thm:conn-strict}, the corollary holds.
\end{proof}

\section{Strong coverings, quotients and site models}
\label{sec:corollaries}
\setcounter{equation}{0}
\noindent
We use the following explicit quotient convention. If group $\Gamma$ acts on graph $\G=(V(\G),E(\G))$ by graph
automorphisms, then $\G/\Gamma$ has vertex set the $\Gamma$-orbits, and two
\emph{distinct} orbits are adjacent when some representatives are adjacent in
$\G$. Edges whose endpoints lie in the same orbit produce no loop, and multiple
edges between the same two orbits are identified. Thus $\G/\Gamma$ is again a graph.
Clearly  $$\mbox{the quotient map}\ \pi:\ x\in\G \rightarrow \Gamma x\in \G/\Gamma\ \mbox{is a weak covering map}.$$

\begin{lemma}[\text{\cite[Lemma 4.3]{Ma-Se2019}}]
\label{lem:strong-disjoint}
Every strong covering map has the disjoint tree-lifting property.
\end{lemma}

\begin{lemma}[\text{\cite[Lemma 4.1]{Ma-Se2019}}]
\label{lem:free-action}
Assume group $\Gamma$ acts freely on $V(\G)$ by graph automorphisms. Let
$\Hh=\G/\Gamma$ and $\pi(x)=\Gamma x,\ x\in V(\G).$ Then $\pi$ has the disjoint tree-lifting property.
\end{lemma}

\begin{corollary}[Weighted covering corollaries]\label{cor:covering}
Assume that the weight pairs $(w_{\G},w_{\Hh})$ are uniformly elliptic and
compatible with the covering map $\pi .$
\begin{enumerate}[label=\textnormal{(\alph*)}]
  \item If $\G$ and $\Hh$ have bounded degree, and
  $\pi:\ V(\G)\to V(\Hh)$ is a strong covering map with tame fibres and
  $\lambda_c(\G,w_{\G})<\infty$, then
  \begin{equation}\label{eq-strict-cover-1}
    \lambda_c(\G,w_{\G})<\lambda_c(\Hh,w_{\Hh}).
  \end{equation}
  \item Let a nontrivial group $\Gamma$ act freely on $V(\G)$ by graph
  automorphisms.  If $\G$ and $\Hh=\G/\Gamma$ are quasi-transitive and
  $\lambda_c(\G,w_{\G})<\infty$, then \eqref{eq-strict-cover-1} is true.
  \item If $\G$ and $\Hh$ are quasi-transitive graphs, $\pi$ is a non-injective
  strong covering map with bounded fibres, and
  $\lambda_{\mathrm{conn}}(\G,w_{\G})<\infty$, then
  \begin{equation}\label{eq-strict-cover-2}
    \lambda_{\mathrm{conn}}(\G,w_{\G})
    <\lambda_{\mathrm{conn}}(\Hh,w_{\Hh}).
  \end{equation}
  \item Let $\Gamma$ act freely as in (b), if the quotient fibres have
  uniformly bounded diameter and $\lambda_{\mathrm{conn}}(\G,w_{\G})<\infty,$ then \eqref{eq-strict-cover-2} holds.
  \item If in (c) or (d), $(\G,w_{\G})$ and $(\Hh,w_{\Hh})$ are
  weighted quasi-transitive and have $\mathsf{Stab}$, $\lambda_u(\G,w_{\G})<\infty,$ then
  $$\lambda_u(\G,w_{\G})<\lambda_u(\Hh,w_{\Hh}).$$
\end{enumerate}
\end{corollary}

\begin{proof}
Part (a) follows from Lemma~\ref{lem:strong-disjoint} and
Theorem~\ref{thm:main}. Part (b) is a consequence of Lemmas~\ref{lem:free-action} and~\ref{lem:automatic-tameness}, and Theorem~\ref{thm:main}.  Part (c) follows from
Lemma~\ref{lem:strong-disjoint} and Theorem~\ref{thm:conn-strict}. Part (d) is implied by Lemma~\ref{lem:free-action} and
Theorem~\ref{thm:conn-strict}. Part (e) follows from
Corollary~\ref{cor:pu-conditional}.
\end{proof}

\begin{remark} A finite quotient group $\Gamma$ gives a uniform bound on the \emph{cardinality} of the
fibres of the original quotient map. This is not the same as a uniform bound
on their diameters in $\G$.  In the homogeneous argument of
\cite{Ma-Se2019}, a compactness construction replaces the original map
by another weak covering whose fibres have uniformly bounded diameter. For nonconstant weights, however, the replacement map need not send a chosen
covering edge $f$ to the same base edge as the original quotient map.  Thus
\[
  w_{\G}(f)\ge w_{\Hh}(\pi(f))
\]
may be lost precisely when one passes from the original finite-cardinality
quotient map to the compactness-produced bounded-diameter map.  Consequently,
``$\Gamma$ finite'' alone is not a safe substitute for bounded fibre diameter
in the weighted connection or uniqueness corollaries. It is sufficient
if compatibility is separately verified for the replacement map, or under the
stronger map-independent domination
$$\inf_{f\in E(\G)} w_{\G}(f)\ge\sup_{e\in E(\Hh)} w_{\Hh}(e).$$
\end{remark}

\begin{remark}[Weighted site-percolation adaptation]
\label{rem:site} Recall the homogeneous site adaptation of the strict argument is discussed
in \cite{Ba-Bo2014, Ma-Se2019}. The site counterparts of Theorems~\ref{thm:main} and~\ref{thm:conn-strict} should be true.

To begin, assign vertex weights $v_{\X}(x)\in[a,b]$ and declare $x$ open with probability $1-e^{-\lambda v_{\X}(x)}.$
The natural compatibility condition is
\[
  v_{\G}(x)\ge v_{\Hh}(\pi(x)),\ x\in V(\G).
\]
At the level of the non-strict critical comparison, the bond exploration has a
direct site analogue: shared uniforms couple a base vertex with a selected lift,
and compatibility makes every successful base discovery a successful lifted
discovery.

For the weighted strict construction, one splits a physical vertex of
weight $v$ into $N$ independent auxiliary Bernoulli coordinates, each with
opening probability $1-e^{-\lambda v/N},$ so that
\[
  1-\Bigl(e^{-\lambda v/N}\Bigr)^N=1-e^{-\lambda v}.
\]
Uniform ellipticity would again give a common positive lower bound for every
fresh auxiliary coordinate; and for a finite increasing site event $\cE,$ the
Margulis--Russo derivative takes the form
\[
 \frac{\partial}{\partial\lambda}\Pp_\lambda^{\Hh}(\cE)
 =\sum_x v(x)e^{-\lambda v(x)}\,
   \Pp_\lambda^{\Hh}(x\text{ is site-pivotal for $\cE$}).
\]
A full site version of the strict $\lambda_c$ and
$\lambda_{\mathrm{conn}}$ comparisons will nevertheless require site-specific
analogues of the fresh-coordinate lemma, the finite-start lifted coupling, and
the deterministic pivotal-conversion surgery, including a precise treatment of
starting-site states and of the vertex boundary used for blocking.
\end{remark}

\section{The strict monotonicity: general profiles and conjectures}\label{sec:conjectures}
\setcounter{equation}{0}
\noindent We use two subsections to discuss respectively general profiles and conjectures in short-range and long-range settings.

\subsection{The short-range situation}\label{subsec:short-range}
\noindent Note that percolations $\Pp_\lambda^{\G}$ and $\Pp_\lambda^{\Hh}$ are nearest-neighbour percolations; and short-range percolation can be viewed as a nearest-neighbour one
by changing suitably the underlying graph structure. The exponential parametrisation is not the only one for which the preceding
arguments apply. We prove a general version in Theorem \ref{thm:general-profile} with the threshold conventions
and constant dependence made explicit.

Let $J\subseteq\mathbb R$ be a nonempty interval.  For
$\X\in\{\G,\Hh\}$ and $e\in E(\X)$, let $p_e^{\X}:\ J\longrightarrow(0,1)$
be an edge-opening profile, and $\Pp_t^{\X}$ the corresponding product
measure. Write $\cC_{\X}^{t}(A)$ for the open cluster of $A$ in this
configuration. The same finite-path Harris--FKG argument used after
Definition~\ref{def:lambda-c} shows that the critical threshold below is
independent of the chosen root.  For a root $o\in V(\X),$ define
\begin{equation}\label{eq:general-tc}
 t_c\left(\X,p^{\X}\right)
 :=\inf\Bigl\{t\in J:\ \Pp_t^{\X}\left(\left|\cC_{\X}^{t}(o)\right|=\infty\right)>0\Bigr\}\ \left(\inf\varnothing:=+\infty\right),
\end{equation}
where the infimum is taken in the extended real line. Likewise set
\begin{equation}\label{eq:general-tconn}
 t_{\mathrm{conn}}\left(\X,p^{\X}\right)
 :=\inf\Bigl\{t\in J:\
   \lim_{\ell\to\infty}\inf_{x,y\in V(\X)}
   \Pp_t^{\X}\left[B_\ell^{\X}(x)\leftrightarrow B_\ell^{\X}(y)\right]=1\Bigr\}\ \left(\inf\varnothing:=+\infty\right).
\end{equation}
The limit in \eqref{eq:general-tconn} exists because the quantity inside the limit is
increasing in $\ell$. Say that one of these thresholds lies in the
interior $\operatorname{int}J$ of $J$ if it is a finite real number in $\operatorname{int}J$.

Assume the following profile hypotheses.
\begin{enumerate}[label=\textnormal{(P\arabic*)},leftmargin=2.8em]
  \item Every $p_e^{\X}$ is increasing and continuously differentiable on
  $J$;
  \item for every compact interval $I\Subset J,$ there exists $\rho_I>0$ such
  that for $\X\in\{\G,\Hh\}$,
  \[
    \rho_I\le p_e^{\X}(t)\le1-\rho_I \ \ (e\in E(\X),\ t\in I);
  \]
  \item for every compact interval $I\Subset J$,
  \begin{equation}\label{eq:MI}
    M_I:=\sup_{e\in E(\Hh),\,t\in I}\left(p_e^{\Hh}\right)'(t)<\infty;
  \end{equation}
  \item if $f\in E(\G)$ lifts $e\in E(\Hh)$, then
  \[
    p_f^{\G}(t)\ge p_e^{\Hh}(t)\ \ (t\in J).
  \]
\end{enumerate}

For $t\in J$ and $s\in[0,1]$, let $\Pp^{\Hh}_{t,s}$ be the product of the
physical $\Hh$-edge law at parameter $t$ and independent $\Ber(s)$ enhancement
marks, and $C^{t,s}_{\Hh}(A)$ the exploratory enhanced cluster
started from $A$.

\begin{theorem}[General-profile extension]\label{thm:general-profile}
Assume \textnormal{(P1)--(P4)}.  Then the following statements hold with the
same geometric hypotheses as the corresponding exponential-model results.
\begin{enumerate}[label=\textnormal{(\roman*)}]
\item Under the graph and map hypotheses of Theorem~\ref{thm:pc-monotone},
 $t_c\left(\G,p^{\G}\right)\le t_c\left(\Hh,p^{\Hh}\right).$

\item Under the graph and map hypotheses of Theorem~\ref{thm:conn-monotone},
  $t_{\mathrm{conn}}\left(\G,p^{\G}\right)
  \le t_{\mathrm{conn}}\left(\Hh,p^{\Hh}\right).$
\item Under the geometric hypotheses of Theorem~\ref{thm:main}, if
$t_c\left(\G,p^{\G}\right)<+\infty$ and either
$t_c\left(\Hh,p^{\Hh}\right)=+\infty$ or
$t_c\left(\Hh,p^{\Hh}\right)\in\operatorname{int}J$, then
$t_c\left(\G,p^{\G}\right)<t_c\left(\Hh,p^{\Hh}\right).$
\item Under the geometric hypotheses of Theorem~\ref{thm:conn-strict},
assume $t_{\mathrm{conn}}\left(\G,p^{\G}\right)<+\infty,$
and in addition either $t_{\mathrm{conn}}\left(\Hh,p^{\Hh}\right)=+\infty$ or
$t_{\mathrm{conn}}\left(\Hh,p^{\Hh}\right)\in\operatorname{int}J$.  Then
  $t_{\mathrm{conn}}\left(\G,p^{\G}\right)
  <t_{\mathrm{conn}}\left(\Hh,p^{\Hh}\right).$
\end{enumerate}
\end{theorem}

\begin{proof}
The non-strict comparisons use only monotonicity and (P4).  The shared-uniform
exploration from Lemma~\ref{lem:connection-lifting} remains valid with $p_{\cdot}^{\Hh}(t)$ and $p_{\cdot}^{\G}(t)$, and therefore proves (i) and (ii)
exactly as before. Now to isolate the two compact-interval inputs needed for
strictness.

\smallskip
\noindent\emph{Claim 1 (compact-interval lifted enhancement).}
Fix a compact interval $I\Subset J$ and assume the geometric hypotheses of
Theorem~\ref{thm:main}.  There exist $r\ge1$ and $s_I\in(0,1)$ such that for
every $t\in I$, every $s\in(0,s_I]$, and every root lift $o'\in\pi^{-1}(o)$,
\begin{equation}\label{eq:general-profile-coupling}
 C^{t,s}_{\Hh}(o)\preceq \pi\Bigl(\cC^t_{\G}(o')\Bigr).
\end{equation}
The same constants work for every finite starting set $A_0\subset V(\G)$:
\begin{equation}\label{eq:general-profile-finite-start}
 C^{t,s}_{\Hh}\left(\pi(A_0)\right)
 \preceq \pi\Bigl(\cC^t_{\G}(A_0)\Bigr).
\end{equation}

\noindent\emph{Proof of Claim 1.}
Fix the integer $N$ and the deterministic reservation maps of Section~\ref{sec:weighted-lifted-coupling}.
Replace a physical edge having opening probability $p_e^{\X}(t)$ by $N$
independent auxiliary copies with probability
\begin{equation}\label{eq:general-copy-profile}
 \widehat p_{e,N}^{\X}(t)
 :=1-\Bigl(1-p_e^{\X}(t)\Bigr)^{1/N}.
\end{equation}
Then
\[
 1-\Bigl(1-\widehat p_{e,N}^{\X}(t)\Bigr)^N=p_e^{\X}(t),
\]
and (P4) is preserved because $p\mapsto1-(1-p)^{1/N}$ is increasing.  By
(P2), if $\rho_I$ is the finite-energy constant on $I$, then
\[
 1-(1-\rho_I)^{1/N}
 \le \widehat p_{e,N}^{\X}(t)
 \le 1-\rho_I^{1/N},
\]
so every auxiliary coordinate, in either state, has probability at least
\begin{equation}\label{eq:general-copy-finite-energy}
 \widehat\rho_{I,N}
 :=\min\Bigl\{1-(1-\rho_I)^{1/N},\,\rho_I^{1/N}\Bigr\}>0.
\end{equation}
The deterministic bookkeeping of Lemma~\ref{lem:fresh-coordinate} is unchanged.
Apply Lemma~\ref{lem:adaptive-deferred} to the adaptive choices.  The proof
of Lemma~\ref{lem:product-law} then applies with the edge-probabilities in
\eqref{eq:general-copy-profile}; thus the completed lifted configuration has
the required product law.  Every enhancement trial requests at most the same
constant $K$ fresh coordinates as in Section~\ref{sec:weighted-lifted-coupling}, and therefore its conditional
success probability is at least $\widehat\rho_{I,N}^{K}$.  Set
\[
 s_I:=\tfrac12\widehat\rho_{I,N}^{K}.
\]
The sequential thinning proof of Lemma~\ref{lem:alpha-law} then produces an
i.i.d.\,\,$\Ber(s)$ enhancement field for every $0<s\le s_I$, uniformly in
$t\in I$.  This proves \eqref{eq:general-profile-coupling}.  The reservation
maps and all freshness statements are independent of the starting set, so the
componentwise argument of Proposition~\ref{prop:finite-start} gives
\eqref{eq:general-profile-finite-start} with the same constants.
\hfill$\triangleleft$

\smallskip
\noindent\emph{Claim 2 (compact-interval differential comparison).}
Fix $r$, a compact $I\Subset J$, and a compact enhancement interval
$J_{\rm en}\Subset(0,1)$.  Write $\theta_L(t,s)$ for the root event in Section~\ref{sec:strict};
and for the finite-set version, write
$$\theta_L(t,s):=\theta_L^{A,B}(t,s)=\Pp_{t,s}^{\Hh}\left(\cE_L^{A,B}\right),\ \cE_L:=\cE_L^{A,B}.$$
There is a constant
$c_I>0$, independent of the volume and of the finite starting/target sets,
such that the following inequality holds for either choice of $\theta_L$
\begin{equation}\label{eq:general-profile-differential}
  \partial_s\theta_L(t,s)
  \ge c_I\,\partial_t\theta_L(t,s),\ \ (t,s)\in I\times J_{\rm en}.
\end{equation}

\noindent\emph{Proof of Claim 2.}
The pivotal surgeries of Lemmas~\ref{lem:surgery} and
\ref{lem:finite-set-pivotal-conversion} are deterministic and depend only on
the graph geometry and on the edge/enhancement states, so they are unchanged.
Their finite-energy counting is performed on physical $\Hh$-edges and on
enhancement marks.  Hypothesis (P2) and compactness of $J_{\rm en}$ therefore give a
uniform constant $c_0>0$ such that
\[
 \partial_s\theta_L(t,s)
 \ge \frac{c_0}{C_R}
     \sum_e\Pp_{t,s}(e\text{ is edge-pivotal for $\cE_L$}).
\]
On the other hand, Lemma~\ref{lem:inhomogeneous-russo} and (P3) give
\[
 \partial_t\theta_L(t,s)
 =\sum_e \left(p_e^{\Hh}\right)'(t)\,
   \Pp_{t,s}(e\text{ is edge-pivotal for $\cE_L$})
 \le M_I\sum_e\Pp_{t,s}(e\text{ is edge-pivotal for $\cE_L$}).
\]
With $M_I^*:=\max\{1,M_I\}$, one may take
$c_I=c_0/(C_RM_I^*)>0$.  This proves
\eqref{eq:general-profile-differential}; the finite-set version is identical
because Lemma~\ref{lem:finite-set-pivotal-conversion} has a uniform surgery
radius and multiplicity bound.
\hfill$\triangleleft$
\smallskip

We now prove (iii). Put $t_c:=t_c\left(\Hh,p^{\Hh}\right)$. When $t_c=+\infty,$ (iii) holds trivially. So suppose that
$t_c<\infty$ and belongs to $\operatorname{int}J.$ Choose a compact interval $I\Subset J$ containing
$t_c$ in its interior and let $s_I$ be given by Claim~1.  Apply Claim~2 on a rectangle
$I\times[\eta,s_I]$, where $0<\eta<s_I/2$, and let $c_I$ be the resulting
constant.  Since the supercritical set in \eqref{eq:general-tc} is an upper
set, every parameter strictly larger than $t_c$ is supercritical.  Choose
$\varepsilon>0$ so small that $t_c\pm\varepsilon\in I$ and
$\eta+2\varepsilon/c_I\le s_I$.  The ordinary model at
$t_c+\varepsilon$ then percolates, and hence so does the enhanced model there
for every positive enhancement parameter.  Integrating
\eqref{eq:general-profile-differential} along the segment of slope
${\rm d}t/{\rm d}s=-c_I$ from
$(t_c+\varepsilon,\eta)$ to
$(t_c-\varepsilon,\eta+2\varepsilon/c_I)$ gives, after passage to the
infinite-volume limit,
\[
 \Pp^{\Hh}_{t_c-\varepsilon,\,\eta+2\varepsilon/c_I}
 \Bigl(\left|C^{t_c-\varepsilon,\,\eta+2\varepsilon/c_I}_{\Hh}(o)\right|=\infty\Bigr)>0.
\]
Claim~1, with $s=\eta+2\varepsilon/c_I\le s_I$, therefore implies
\[
 t_c\left(\G,p^{\G}\right)\le t_c-\varepsilon<t_c=t_c\left(\Hh,p^{\Hh}\right).
\]

For (iv), we only need to consider the case that
$t_{\mathrm{conn}}\left(\Hh,p^{\Hh}\right)<\infty$ lies in $\operatorname{int}J$.
Under the geometric hypotheses of Theorem~\ref{thm:conn-strict},
Lemma~\ref{lem:automatic-tameness} supplies tame fibres, so Claim~1 and its
finite-start version apply.  Fix a compact interval $I\Subset J$ around
$t_{\mathrm{conn}}\left(\Hh,p^{\Hh}\right)$, let $r$ and $s_I$ be given by Claim~1, and
put $s:=s_I/2$.  Apply Claim~2 on $I\times[s/2,s]$, obtaining a constant
$c_I>0$.  Choose $\delta>0$ small enough that $2\delta/s\le c_I$ and that an
interval of length $\delta$ centred sufficiently near
$t_{\mathrm{conn}}\left(\Hh,p^{\Hh}\right)$ remains in $I$.  For any $t$ with
$[t,t+\delta]\subset I$ and $u\in[s/2,s]$, set
\[
 \mu(u):=t+\delta-\frac{2\delta}{s}(u-s/2).
\]
Then $\mu(s/2)=t+\delta$, $\mu(s)=t$, and
$\mu'(u)=-2\delta/s\ge-c_I$. By \eqref{eq:general-profile-differential},
\[
 \frac{{\rm d}}{{\rm d}u}\theta_L^{A,B}(\mu(u),u)\ge0
\]
for every finite-volume connection event of Section~\ref{sec:uniqueness}.  Passing to the limit
exactly as in Proposition~\ref{prop:finite-set-enhancement} yields, uniformly
in the finite sets $A,B$ satisfying the $r$-nice and separation conditions,
\begin{equation}\label{eq:general-profile-finite-set-shift}
 \Pp^{\Hh}_{t,s}[A\Rightarrow B]
 \ge \Pp^{\Hh}_{t+\delta}[A\leftrightarrow B].
\end{equation}
The constant $\delta$ is independent of $A$ and $B$.

Combining \eqref{eq:general-profile-finite-start} and
\eqref{eq:general-profile-finite-set-shift}, the rest of the proof of
Theorem~\ref{thm:conn-strict} is purely geometric: use the explicitly verified
$r$-nice target set, the bounded-fibre metric comparison, and the
cofinal-radius argument.  Choose $t<t_{\mathrm{conn}}\left(\Hh,p^{\Hh}\right)<t+\delta$ with
$[t,t+\delta]\subset I$.  Because the defining set in
\eqref{eq:general-tconn} is an upper set, $t+\delta$ satisfies the uniform
connection criterion on $\Hh$.  The preceding comparison therefore yields
\[
 t_{\mathrm{conn}}\left(\G,p^{\G}\right)
 <t_{\mathrm{conn}}\left(\Hh,p^{\Hh}\right).
\]
\end{proof}

We now propose a few conjectures as follows.

\begin{conjecture}[Finite-type stability]\label{conj:stability}
Every uniformly elliptic weighted quasi-transitive exponential model has
property $\mathsf{Stab}$, and particularly
\[
 \lambda_{\mathrm{conn}}(\X,w)=\lambda_u(\X,w).
\]
\end{conjecture}

If Conjecture~\ref{conj:stability} holds, Corollary
\ref{cor:pu-conditional} is an unconditional weighted extension of
\cite[Theorems~1.3-1.4]{Ma-Se2019}.  The exponential-clock coupling
does provide the uniform sprinkling bound
\begin{eqnarray*}
 &&\Pp(e\text{ is added between }\lambda\text{ and }\lambda'\mid
       e\text{ closed at }\lambda)\\
 &&=1-\exp\left\{-(\lambda'-\lambda)w(e)\right\}
   \ge1-\exp\left\{-(\lambda'-\lambda)a\right\}>0,\ \forall\,\lambda'>\lambda,
\end{eqnarray*}
but a complete anisotropic analogue of the stability argument is unproven. In particular, the bound alone does not establish the inhomogeneous
cluster-containment/merging conclusion that yields (S1), and does not prove
the large-ball connectivity characterization in (S2).

\begin{conjecture}[Weak coverings without disjoint tree lifting]
\label{conj:weak-cover}
Suppose $\pi:V(\G)\to V(\Hh)$ is a weak covering map with tame fibres and the
weights are uniformly elliptic and compatible.  If
$\lambda_c(\G,w_{\G})<\infty$, then
\[
 \lambda_c(\G,w_{\G})<\lambda_c(\Hh,w_{\Hh}).
\]
\end{conjecture}

This is already open in the homogeneous setting; it is the weighted counterpart
of \cite[Question~5.2]{Ma-Se2019}.

\begin{conjecture}[Connection threshold under weak coverings]
\label{conj:weak-connection}
Let $\G$ and $\Hh$ be quasi-transitive graphs and $\pi:\ V(\G)\to V(\Hh)$ a non-injective weak covering map with tame and
bounded fibres.  Under uniformly elliptic compatible weights,
\[
 \lambda_{\mathrm{conn}}(\G,w_{\G})<
 \lambda_{\mathrm{conn}}(\Hh,w_{\Hh})\ \mbox{if}\  \lambda_{\mathrm{conn}}(\G,w_{\G})<\infty.
\]
\end{conjecture}

Together with Conjecture~\ref{conj:stability}, this would yield the weighted
analogue of \cite[Question~5.3]{Ma-Se2019}. The obstruction is again
geometric: the present enhancement coupling has no uniformly duplicated exit
without the disjoint lifting hypothesis.


\begin{proposition}[Restricted lower ellipticity]\label{prop:ellipticity}
In Theorems~\ref{thm:main} and~\ref{thm:conn-strict}, retain their geometric
hypotheses, compatibility, finiteness of the left-hand threshold, positive
weights on $\G$, a uniform upper bound on all weights, and uniform
ellipticity on $\Hh$.
Fix the deterministic geometric choices of
Section~\ref{sec:weighted-lifted-coupling}. Let $\mathcal F_{\mathrm{en}}$ be the union, over
every possible reached lift $x$, of the edges in the chosen spanning tree
of $Z(x,r)$ and of both candidate exits for every
$v\in S_{r+1}^{\Hh}(\pi(x))$. The uniform lower bound on all $\G$-weights may be
replaced by
\[
 a_{\mathrm{en}}:=\inf_{f\in\mathcal F_{\mathrm{en}}}w_{\G}(f)>0.
\]
Thus the strict comparisons remain valid even if
$\inf\limits_{f\in E(\G)}w_{\G}(f)=0$.
\end{proposition}

\begin{proof}
The geometric input, reservation maps and fresh-coordinate bookkeeping
do not depend on the weights. The set $\mathcal F_{\mathrm{en}}$ is fixed before sampling
and contains every physical edge that can request a fresh coordinate in an
enhancement, regardless of which candidate exit is selected after the
ordinary exploration.
For $\lambda\ge\lambda_0>0$, every such request has success probability
at least
\[
 q_{\mathrm{en}}=1-\exp\left(-\lambda_0a_{\mathrm{en}}/N\right)>0.
\]
There are at most $K$ requests in one trial, with the same uniform
geometric bound $K$ as before. Hence its conditional success probability
is at least $q_{\mathrm{en}}^K$, and the thinning construction works
with $s=\tfrac12q_{\mathrm{en}}^K$.

The ordinary coupling uses only compatibility. Completion of the lifted
product measure needs positivity of each edge weight but no uniform lower
bound on edges that are not requested in an enhancement. Consequently
Propositions~\ref{prop:coupling} and~\ref{prop:finite-start} hold under
the stated restricted condition. All pivotal surgeries and differential
inequalities are on $\Hh$, whose uniform ellipticity is unchanged.
The proofs of the two strict comparisons therefore go through with this
new common enhancement parameter.
\end{proof}

The role of the restricted infimum is to supply a position-independent
success probability for coordinates actually consumed by enhancement.
If it vanishes, the local success probabilities may tend to zero, and this
argument does not give a single $s>0$. The base lower bound is still used
when the pivotal surgeries in Sections~\ref{sec:strict}
and~\ref{sec:uniqueness} force finitely many physical base edges open.

\begin{example}[An example with vanishing weights on unused fibre edges]\label{ex:vanishing-weight}
The restriction in Proposition~\ref{prop:ellipticity} is weaker than a
lower bound on all covering edges. To see this, take
$\Hh=\mathbb Z^2$ with unit weights and
$\G=\mathbb Z^2\times K_3$ with layer set $\{0,1,2\}$ and the layer
projection. Give all horizontal edges and the fibre edges joining layers
$0$ to $1$ and $1$ to $2$ weight one. Enumerate the lattice vertices as
$z_1,z_2,\ldots$ and give the remaining fibre edge
$\{(z_n,0),(z_n,2)\}$ weight $1/(n+1)$. All weights are positive and
at most one, compatibility holds, but their infimum on $\G$ is zero.

The projection is a strong covering, its fibres have diameter one, and
the unweighted graphs are quasi-transitive. We may choose $R=1$ in
tameness and hence $r=1$ in Lemma~\ref{lem:geometry}. For every
$x=(z,i)$, choose its companion lift in a layer adjacent to $i$ in the
path $0$--$1$--$2$. The one-edge path to this companion has weight one,
and the two lifts of the base tree can be chosen in these two layers.
In this product graph,
 $$Z(x,1)=B_1^{\mathbb Z^2}(z)\times\{0,1,2\}.$$
{\it Indeed}, any vertex of the right-hand side lies within distance two of
$x$, and the induced subgraph is connected.

Inside $Z(x,1)$, the horizontal edges and the fibre path $0$--$1$--$2$ form a connected
spanning subgraph consisting entirely of unit-weight edges. Fix a spanning
tree of this subgraph for each $x$. Both candidate exits furnished by
the two layer lifts are horizontal and have weight $1$. With these
deterministic geometric choices,
\[
 w_{\G}(f)=1\ \text{for every}\ f\in\mathcal F_{\mathrm{en}},
 \ a_{\mathrm{en}}=1.
\]
Thus none of the fibre edges with weights tending to zero is requested
by enhancement.

The required finiteness of the left-hand thresholds also holds. The
spanning subgraph $\mathbb Z^2\times P_3$ has only unit weights; where
$P_3$ is the path on the 3 layers. This subgraph has finite
critical and uniform connection intensities.
For the latter assertion use homogeneous uniqueness on this amenable
quasi-transitive graph and Lemma~\ref{lem:connection-criterion}(i).
Adding edges preserves the large-ball criterion, since
$B_\ell^{\mathbb Z^2\times P_3}(x)\subseteq B_\ell^{\G}(x)$.
By Proposition~\ref{prop:ellipticity},
\begin{eqnarray}\label{eq:strict-example}
 \lambda_c(\G,w_{\G})<\lambda_c\left(\mathbb Z^2,\mathbf1\right),\
 \lambda_{\mathrm{conn}}(\G,w_{\G})
 <\lambda_{\mathrm{conn}}\left(\mathbb Z^2,\mathbf1\right).
\end{eqnarray}
Here the non-collapsed edges already inherit the base lower bound by
compatibility. The genuine relaxation concerns redundant edges inside
the fibres, which can be excluded from the fixed enhancement geometry.

Note \eqref{eq:strict-example} also follow by domination from the
unit-weight spanning subgraph $\mathbb Z^2\times P_3$ and the homogeneous
covering theorem of \cite{Ma-Se2019} $\left(\text{$\lambda_c=\lambda_{\mathrm{conn}}$ on amenable
quasi-transitive graphs $\mathbb Z^2\times P_3$ and $\mathbb Z^2$}\right).$
\end{example}



\begin{conjecture}[Finite free quotients for uniqueness]
\label{conj:finite-quotient}
Let a finite nontrivial group $\Gamma$ act freely on a quasi-transitive graph $\G$, and
$\Hh=\G/\Gamma$. Suppose that $(\G,w_{\G})$ and
$(\Hh,w_{\Hh})$ are weighted quasi-transitive, the weights are uniformly elliptic and compatible with the quotient map,
and $\lambda_u(\G,w_\G)<\infty.$ Then
$$\lambda_u(\G,w_\G)<\lambda_u(\Hh,w_{\Hh}).$$
\end{conjecture}

The homogeneous compactness argument replaces the quotient map by a new
bounded-fibre weak covering.  In the weighted setting the replacement map need
not preserve compatibility, which is the additional unresolved obstruction.

\subsection{The long-range setting}\label{subsec:long-range}
\noindent
{\bf Long-range percolation on transitive graphs.} Let $(\X,E(\X))$ be a locally finite connected infinite transitive graph with vertex set $\X$ and edge set $E(\X).$ Fix a vertex $o\in\X$ and a group ${\rm S}\subseteq {\rm Aut}(\X)$ acting transitively on $\X.$
Let $\X^{[2]}=\{\{x,y\}:\ x,y\in \X,\ x\not=y\}.$ Recall \eqref{eq:ball} and \eqref{eq:weighted-probability}. For any weight function $w:\ \X^{[2]}\rightarrow [0,\infty]$,
call the percolation $\Pp_{\lambda,w}^{\X}$ ($\lambda\geq 0$) a long-range percolation on $(\X,E(\X)).$  When
$$w(\{x,y\})=w(\{sx,sy\}),\ \{x,y\}\in\X^{[2]},\ s\in {\rm S},$$
say $w$ is ${\rm S}$-invariant. Write $\mathscr{W}:=\mathscr{W}(\X,{\rm S})$ for the set of all ${\rm S}$-invariant weight functions $w:\ \X^{[2]}\rightarrow [0,\infty]$ such that
$\sum\limits_{x\in\X\setminus B_r^{\X}(o)}w(\{o,x\})<\infty\ \mbox{for large enough}\ r.$
Let $\mathscr{W}_{\rm f}:=\mathscr{W}_{\rm f}(\X,{\rm S})$ be the subset of $\mathscr{W}$ consisting of all finite weight functions $w:\ \X^{[2]}\rightarrow [0,\infty).$
For any two weight functions $w$ and $w'$ on $\X^{[2]},$ write $w'<w$ if
$$w'(e)\leq w(e)\ \mbox{for any}\ e\in \X^{[2]},\ \mbox{and}\ w'(f)<w(f)\ \mbox{for some}\ f\in\X^{[2]}.$$

By Bethuelsen and M\"{o}nch (2025) \cite[Theorems 2.2 and 2.3]{Bet-Mo2025}, the following strict monotonicity holds for long-range percolations on $(\X,E(\X))$: Given any $w\in\mathscr{W}$ and $w'\in\mathscr{W}_{\rm f}$ with $w'<w$ and $\lambda_c(\X,w)<\infty.$ Then
\begin{eqnarray}\label{eq:long-range-strict}
\lambda_c(\X,w)<\lambda_c(\X,w').
\end{eqnarray}

Naturally the following monotonic continuity (locality) conjecture arises:
\begin{conjecture}\label{conj:continuity}
In the transitive case, given any sequence $\{w_n\}_{n=1}^{\infty}\subseteq\mathscr{W}_{\rm f}$ and $w\in\mathscr{W}$ such that
\begin{eqnarray*}
w_n(e)\leq w_{n+1}(e),\ \lambda_c(\X,w_n)<\infty,\ n\in\mathbb{N},\ \lim\limits_{k\rightarrow\infty}w_k(e)=w(e),\ e\in\X^{[2]},\
 \mbox{and}\ \lambda_c(\X,w)>0.
\end{eqnarray*}
Then
 $$\lim\limits_{n\rightarrow\infty}\lambda_c(\X,w_n)=\lambda_c(\X,w).$$
\end{conjecture}
\noindent For the locality of the critical percolation probabilities of nearest-neighbour Bernoulli bond percolations, see P. Easo and T. Hutchcroft \cite{Ea-Hu2023} and references therein.
And for the locality (continuity) of the critical parameters of long-range percolations on $\mathbb{Z}^d$, refer to J. B\"{a}umler \cite[Theorem 1.6]{Ba2026}.
\smallskip

\begin{remark} Conjecture~\ref{conj:continuity} may fail without $\lambda_c(\X,w)>0.$ See the following example: For any $n\ge 1,$ define a weight on $\mathbb{Z}^{[2]}$ by
\[
 w_n(\{x,y\}):=
 \begin{cases}
 n,&|x-y|=1,\\
 |x-y|^{-2},&|x-y|\ge 2.
 \end{cases}
\]
Clearly $\{w_n\}_{n=1}^{\infty}\subset\mathscr{W}_{\rm f}$ converges increasingly pointwise to a $w\in\mathscr{W}$. Then
\[
 \lambda_c(\mathbb Z,w_n)<\infty\ (n\ge 1),\ \lim_{n\to\infty}\lambda_c(\mathbb Z,w_n)=1>0=\lambda_c(\mathbb Z,w).
\]
{\it Indeed}, consider the inverse-square model with nearest-neighbour opening probability $1-e^{-h}$ and other edge probabilities
$1-e^{-\beta/|x-y|^2}$. If $\theta(\beta,h)$ is its percolation probability, then by \cite{Ai-Ne1986, Du-Ga2024},
\[
 \theta(\beta,h)>0\ \Longrightarrow\ \beta\theta(\beta,h)^2\ge 1,
 \ \
 \mbox{and}\ \beta>1\ \Longrightarrow\
 \theta(\beta,h)>0\text{ for some finite }h.
\]
For finite $\beta,h$, the probability that the origin is isolated is
\[
 \exp\left(-2h-2\beta\sum_{k=2}^{\infty}k^{-2}\right)>0.
\]
Hence $\theta(\beta,h)<1$, and there is no percolation for $\beta\le1$.
The model with weight $w_n$ at intensity $\lambda$ has
$(\beta,h)=(\lambda,n\lambda)$, so
$\lambda_c(\mathbb Z,w_n)\ge 1$. Each $\lambda_c(\mathbb Z,w_n)$ is finite: choose $h_2$ with $\theta(2,h_2)>0$ and use
any $\lambda\ge\max\{2,h_2/n\}$.

For each fixed $\lambda>1$, choose $h<\infty$ giving
$\theta(\lambda,h)>0$. Then $n\lambda\ge h$ for all sufficiently large
$n$, and monotonicity gives percolation for $w_n$ at $\lambda$.
So the decreasing sequence of critical intensities converges to
$1$.

The limiting weight $w$ equals $\infty$ on nearest-neighbour pairs and
$|x-y|^{-2}$ on the other pairs. At every positive intensity its
nearest-neighbour edges are all open. Thus $\lambda_c(\mathbb{Z},w)=0.$
\end{remark}
\vskip 2mm

Restrict to infinite finitely generated groups $\X$. For a symmetric probability $\mu$ on $\X$, define a weight function $w_{\mu}$ as follows:
$$w_{\mu}(\{x,y\}):=\mu\left(x^{-1}y\right),\ \{x,y\}\in \X^{[2]}.$$
Write $\lambda_c(\X,\mu)$ for $\lambda_c(\X,w_\mu).$ Then the following monotonicity holds:

\begin{theorem}\label{thm:long-range quotient} Let $\pi:\ \X\rightarrow \Y$ be a surjective homomorphism of infinite finitely generated groups. If $\mu_{\Y}$ is a symmetric probability on $\Y,$ then every symmetric probability $\mu_{\X}$ on $\X$ with $\pi_*\mu_{\X}=\mu_{\Y}$ satisfies
$$\lambda_c(\X,\mu_{\X})\leq\lambda_c(\Y,\mu_{\Y}).$$
\end{theorem}

\begin{proof} Given any $\lambda\in (0,\infty).$ Consider percolations $\Pp_{\lambda,w_{\mu_\X}}^{\X}$ and $\Pp_{\lambda,w_{\mu_{\Y}}}^{\Y}.$  Notice $\Pp_{\lambda,w_{\mu_\X}}^{\X}$
can be defined equivalently as follows: Define a random multi-graph $\X_{\mu_\X}(\lambda)$ on $\X$ by letting the number of parallel edges between any two elements $x,y\in \X$ be an independent Poisson random variable with mean $\lambda\mu_{\X}\left(x^{-1}y\right).$ Keeping only one parallel edge (if exists) between any two vertices, one gets the
percolation  $\Pp_{\lambda,w_{\mu_\X}}^{\X}$. And $\Pp_{\lambda,w_{\mu_{\Y}}}^{\Y}$ can also be defined in this way.

Consider the following coupling. Suppose a quotient vertex $v\in\Y$ has been assigned a lift vertex $x\in\X.$ For a quotient neighbor $w=vq$, all lift-edges from $x$ projecting to $\{v,w\}$ have aggregate Poisson mean
$$\lambda\sum\limits_{\pi (g)=q}\mu_{\X}(g)=\lambda\mu_{\Y}(q).$$
Explore the quotient root cluster breadth first. Whenever an untested edge has one discovered and one undiscovered endpoint, expose this family of lift-edges. Declare the quotient open exactly when at least one is open, and then assign the undiscovered endpoint one such open lift. Families used for distinct quotient edges are disjoint, so these declarations are independent and have exactly the quotient law. Quotient collisions merely discard possible new lift vertices. Hence an infinite quotient cluster produces an infinite lift cluster.
\end{proof}

\begin{conjecture}\label{conj:long-range quotient} In the situation of Theorem~\ref{thm:long-range quotient}, if $\X$ is not isomorphic to $\Y$ and $\lambda_c(\X,\mu_{\X})<\infty,$ then
$$\lambda_c(\X,\mu_{\X})<\lambda_c(\Y,\mu_{\Y}).$$
\end{conjecture}
\noindent The above conjecture should hold in more general covering map setting.\\

Now assume $(\X,E(\X))$ is a locally finite connected infinite quasi-transitive graph, and group ${\rm S}\subseteq {\rm Aut}(X)$ acts quasi-transitively on $\X.$ Write $\mathscr{W}_{\rm qu}$ for the set of weight functions $w:\ \X^{[2]}\rightarrow [0,\infty]$ such that $\left(\X^{[2]},w\right)$ is quasi-transitive under the ${\rm S}$-action, and for any $x\in\X,$
 $$\sum\limits_{y\in\X\setminus B_r^{\X}(x)}w(\{x,y\})<\infty\ \mbox{for large enough}\ r.$$
Let $\mathscr{W}_{{\rm qu,f}}$ be the subset of $\mathscr{W}_{\rm qu}$ consisting of all finite weight functions $w:\ \X^{[2]}\rightarrow [0,\infty).$
\smallskip

\begin{remark}
The following natural question arises: In the quasi-transitive setting, for any $w\in \mathscr{W}_{\rm qu}$ with $\lambda_c(\X,w)<\infty,$ and any $w'\in \mathscr{W}_{{\rm qu,f}}$ with $w'<w,$
 $$\mbox{does}\ \lambda_c(\X,w)<\lambda_c(\X,w')\ \mbox{hold}?$$
Recall \eqref{eq:long-range-strict} in the transitive case. The following example shows that {\it a quasi-transitive perturbation may not change the thresholds.}

Attach one pendant vertex $v^*$ to every vertex $v$ of the
nearest-neighbour lattice $\mathbb Z^2$. Translations act on the resulting
graph $\X$ with two vertex orbits. Put $w=1$ on each lattice edge and each
pendant edge, and put $w'=1$ on lattice edges and $w'=1/2$ on pendant
edges. Both weights are zero on all other pairs. They are invariant,
finite and summable, and $w'<w$.

In every configuration, a finite lattice cluster of size $m$ acquires at
most $m$ pendant vertices. Pendant edges cannot join distinct lattice
clusters. Thus an infinite cluster exists if and only if the lattice
restriction has an infinite cluster. The latter has the same law for the
two weights, and its critical intensity is
$-\log(1-p_c(\mathbb Z^2))\in(0,\infty)$.
This proves
\[
 0<\lambda_c(\X,w')=\lambda_c(\X,w)<\infty.
\]

The changed edge orbit in this example consists of finite decorations.
{\it A strictness extension to quasi-transitive weights must therefore control
how the perturbed orbit contributes to connections to infinity.}
\end{remark}

\begin{conjecture}\label{conj:strict-quasi-transitive}
In the quasi-transitive setting, for any $w\in \mathscr{W}_{\rm qu}$ with $\lambda_c(\X,w)<\infty$ and any $w'\in \mathscr{W}_{{\rm qu,f}}$ with $w'<w$ such that
$\left\{f\in\X^{[2]}:\ w(f)\not=w'(f)\right\}$ induces a graph on $\X$ with an infinite connected component,
 $$\lambda_c(\X,w)<\lambda_c(\X,w').$$
\end{conjecture}

\renewcommand\thesection{\appendixname~\Alph{section}}
\setcounter{section}{0}
\section{Proofs for two covering monotonicities}\label{app-1}
\setcounter{equation}{0}
\renewcommand\thetheorem{A.\arabic{theorem}}
\renewcommand\theequation{A.\arabic{equation}}
\noindent
This appendix proves Theorems~\ref{thm:pc-monotone} and \ref{thm:conn-monotone}. Recall $A\leftrightarrow B$ denotes the event that an open path meets both vertex sets $A,B$.

\begin{lemma}[Connection and uniqueness]\label{lem:connection-criterion}
Let $(\X,w)$ be an infinite weighted quasi-transitive graph.
\begin{enumerate}[label=\textnormal{(\roman*)}]
  \item If the model has almost surely a unique infinite cluster at $\lambda$,
  then
  \begin{equation}\label{eq:uniform-connection-supercritical}
    \lim_{\ell\to\infty}\inf_{x,y\in V(\X)}
    \Pp_\lambda^{\X}\left[B_\ell^{\X}(x)\leftrightarrow B_\ell^{\X}(y)\right]=1.
  \end{equation}
  \item If the model has property $\mathsf{Stab}$, then
  $\lambda_{\mathrm{conn}}(\X,w)=\lambda_u(\X,w).$
\end{enumerate}
\end{lemma}

\begin{proof}
Choose a quasi-transitive weight-preserving subgroup
$\Gamma\le\Aut(\X,w)$.  Suppose that
$I_\lambda$ is the almost surely unique infinite cluster. For any $x\in V(\X)$,
\[
  \Pp_\lambda^{\X}\left[B_\ell^{\X}(x)\cap I_\lambda\ne\varnothing\right]\longrightarrow 1\ (\ell\rightarrow\infty).
\]
There are only finitely many $\Gamma$-orbits of $x$, so the convergence is
uniform.  If both $B_\ell^{\X}(x)$ and $B_\ell^{\X}(y)$ meet $I_\lambda$, they are openly
connected. By Harris-FKG inequality,
\begin{align*}
 \Pp_\lambda^{\X}\left[B_\ell^{\X}(x)\leftrightarrow B_\ell^{\X}(y)\right]
 &\ge \Pp_\lambda^{\X}\left[B_\ell^{\X}(x)\leftrightarrow\infty,\
                   B_\ell^{\X}(y)\leftrightarrow\infty\right]\\
 &\ge \Pp_\lambda^{\X}\left[B_\ell^{\X}(x)\leftrightarrow\infty\right]\,
       \Pp_\lambda^{\X}\left[B_\ell^{\X}(y)\leftrightarrow\infty\right],
\end{align*}
which proves (i).

Assume $\mathsf{Stab}$. If $\lambda>\lambda_u$, the definition of the infimum
provides a $\mu<\lambda$ at which uniqueness holds; (S1) gives uniqueness at
$\lambda$, and part (i) gives the connection limit.  Hence
$\lambda_{\mathrm{conn}}\le\lambda_u$.  Conversely, if the connection limit
holds at $\lambda$, (S2) gives uniqueness at every $\lambda'>\lambda$, so
$\lambda_u\le\lambda$.  Taking the infimum over such $\lambda$ proves the
reverse inequality.
\end{proof}

\begin{remark}\label{rem:homogeneous-stability}
For homogeneous Bernoulli percolation $\X[p]$  on a quasi-transitive graph $\X$, uniqueness
is monotone in $p\in [0,1]$, $p_u(\X)$ admits a large-ball connectivity characterization (\cite{Sc1999, Ha-Pe1999}), and
$$p_{\mathrm{conn}}(\X)=p_u(\X).$$
Here $p_{\mathrm{conn}}(\X)$ is the homogeneous large-ball connection threshold. While for inhomogeneous Bernoulli percolation, we propose Conjecture \ref{conj:stability}.
\end{remark}

\begin{lemma}[Weighted connection lifting]\label{lem:connection-lifting}
Let $\pi:V(\G)\to V(\Hh)$ be a weak covering map and suppose the weights are
compatible with $\pi$.  For all finite $A,B\subseteq V(\Hh)$ and all
$\lambda\ge0$,
\begin{equation}\label{eq:connection-lifting}
 \Pp_\lambda^{\G}\left[\pi^{-1}(A)\leftrightarrow\pi^{-1}(B)\right]
 \ge \Pp_\lambda^{\Hh}[A\leftrightarrow B].
\end{equation}
\end{lemma}

\begin{proof}
For each $a\in A$, choose one lift $\widetilde a\in\pi^{-1}(a)$ and put
$\widetilde A=\{\widetilde a:a\in A\}$.  Explore the $\Hh$-open cluster of $A$
one boundary edge at a time and simultaneously lift an exploration forest to
$\G$, starting from $\widetilde A$.  Whenever a previously unrevealed base edge
$e=\{u,v\}\in E(\Hh)$ is tested from a reached vertex $u\in V(\Hh)$, select, by the weak lifting
property, an edge $f\in E(\G)$ above $e$ from a reached lift of $u$.  Use a fresh uniform
variable $U$ and declare
\[
  e\text{ open if }U\le 1-\exp\{-\lambda w_{\Hh}(e)\},\  f\text{ open if }U\le1-\exp\{-\lambda w_{\G}(f)\}.
\]
Thus, by the compatibility, every open edge accepted by the base exploration has an accepted open lift.  Coordinates
not selected by the exploration are filled with independent uniforms.  The
selected lifted edges are distinct, because each projects to the base edge
currently being tested and every base edge is tested only once.
 Consequently
the principle of deferred decisions (Lemma \ref{lem:adaptive-deferred}) gives the correct product marginals on
both graphs. Thus one can get a coupling for
$\Pp_{\lambda}^{\G}$ and $\Pp_{\lambda}^{\Hh}$ such that if the base exploration reaches $B$, the lifted forest connects
$\widetilde A$ to some vertex of $\pi^{-1}(B)$. Since
$\widetilde A\subseteq\pi^{-1}(A)$, this proves
\eqref{eq:connection-lifting}.  
\end{proof}

\begin{proof}[Proof of Theorem~\ref{thm:pc-monotone}]
Fix roots $o\in V(\Hh)$ and $o'\in\pi^{-1}(o)$.  Run the exploration in the
proof of Lemma~\ref{lem:connection-lifting} from $\{o\}$ and its chosen lift
$\{o'\}$.  Under that coupling, every vertex reached by the base open cluster
has a reached lift connected to $o'$.  Hence
\[
  \left|\cC_{\Hh}^{\lambda}(o)\right|=\infty \ \Longrightarrow\ \left|\cC_{\G}^{\lambda}(o')\right|=\infty.
\]
This implies that any intensity supercritical for $\Hh$ is supercritical for $\G$, and the theorem is true.
\end{proof}

\begin{lemma}[Metric comparison]\label{lem:metric-comparison}
If $\pi$ has fibre-diameter bound $K$, then for all $x,y\in V(\G)$,
\begin{equation}\label{eq:metric-comparison}
 d_{\G}(x,y)-K\le d_{\Hh}(\pi(x),\pi(y))\le d_{\G}(x,y),
\end{equation}
and for every $\ell\ge0$,
\begin{equation}\label{eq:preimage-ball}
 \pi^{-1}\left(B_\ell^{\Hh}(\pi(x))\right)\subseteq B_{\ell+K}^{\G}(x).
\end{equation}
\end{lemma}

\begin{proof}
The upper bound in \eqref{eq:metric-comparison} is the $1$-Lipschitz property.
Lift a geodesic from $\pi(x)$ to $\pi(y)$ starting at $x$.  Its endpoint $y'$ is
in the fibre of $y$, so $d_{\G}(y',y)\le K$.  This gives the lower bound.
The same argument, with an arbitrary vertex of
$\pi^{-1}\left(B_\ell^{\Hh}(\pi(x))\right)$ in place of $y$, proves
\eqref{eq:preimage-ball}.
\end{proof}

\begin{proof}[Proof of Theorem~\ref{thm:conn-monotone}]
Fix any intensity $\lambda$ at which the uniform connection limit holds on
$\Hh$. Then
\[
 \lim_{\ell\to\infty}\inf_{u,v\in V(\Hh)}
 \Pp_\lambda^{\Hh}\left[B_\ell^{\Hh}(u)\leftrightarrow B_\ell^{\Hh}(v)\right]=1.
\]
For $x,y\in V(\G)$, by Lemma~\ref{lem:connection-lifting} $\left(\mbox{with
$A=B_\ell^{\Hh}(\pi(x))$ and $B=B_\ell^{\Hh}(\pi(y))$}\right),$ and \eqref{eq:preimage-ball}, we have
\[
 \Pp_\lambda^{\G}\left[B_{\ell+K}^{\G}(x)\leftrightarrow B_{\ell+K}^{\G}(y)\right]
 \ge
 \Pp_\lambda^{\Hh}\left[B_\ell^{\Hh}(\pi(x))\leftrightarrow
                     B_\ell^{\Hh}(\pi(y))\right].
\]
Recall Definition~\ref{def:lambda-conn}. Taking the infimum over $x,y$, we get
\[
 a_{\ell+K}^{\G}(\lambda)\ge a_\ell^{\Hh}(\lambda)\longrightarrow 1\ (\ell\rightarrow\infty),
\]
which implies $\lim\limits_{R\to\infty}a_R^{\G}(\lambda)=1.$
Thus $\lambda_{\mathrm{conn}}(\G,w_{\G})\le \lambda_{\mathrm{conn}}(\Hh,w_{\Hh})$ by taking the infimum.
\end{proof}

\renewcommand\thesection{\appendixname~B}
\setcounter{section}{0}
\section{Proof of Lemma \ref{lem:surgery}}\label{app-2}
\setcounter{equation}{0}
\renewcommand\thetheorem{B.\arabic{theorem}}
\renewcommand\theequation{B.\arabic{equation}}
\noindent
\begin{proof} This lemma is \cite[Lemma 2.4]{Ma-Se2019}. For completeness, a detailed proof is given here according to that of \cite[Lemma 2.4]{Ma-Se2019}.
We identify an enhancement configuration with the set of vertices whose mark
is one. Write $\omega^{e=i}$ for the configuration obtained by forcing the
state of $e$ to be $i$.  Fix deterministic enumerations of the relevant edges
and vertices and use them for all tie-breaking choices below.

For the event $\cE_L$, edge $e$ is pivotal exactly when
\begin{equation}\label{eq:root-e-pivotal}
  \cE_L\left(\omega^{e=1},\alpha\right)=1,
  \
  \cE_L\left(\omega^{e=0},\alpha\right)=0.
\end{equation}
Here $\cE_L\left(\omega^{e=1},\alpha\right)=1$ means $\cE_L\left(\omega^{e=1},\alpha\right)$ occurs, while $\cE_L\left(\omega^{e=0},\alpha\right)=0$ means $\cE_L\left(\omega^{e=0},\alpha\right)$ does not occur.

Accordingly, at the first step of the deterministic surgery, we \emph{normalize}
the input by replacing $\omega$ with $\omega^{e=1}$.
All subsequent choices are made from $\omega^{e=1}$. Thus the resulting surgery is defined on every edge-pivotal
input, irrespective of its actual $e$-state, and the normalization merely
increases the uniform bound on the number of changed local coordinates by at
most one. Erase the enhancement marks in $B_R^{\Hh}(e)$ one at a time. If the event first
fails when the mark at a vertex $z$ is erased, then $z$ is already
enhancement-pivotal and we are done. Otherwise put
\[
  \alpha_0:=\alpha\setminus B_R^{\Hh}(e);
\]
and by monotonicity,
\begin{equation}\label{eq:root-e-remains-pivotal}
  \cE_L\left(\omega^{e=1},\alpha_0\right)=1,
  \
  \cE_L\left(\omega^{e=0},\alpha_0\right)=0.
\end{equation}

Resolve each ordinary open-cluster closure in the successful exploration
$\left(\omega^{e=1},\alpha_0\right)$ by a fixed edge-by-edge exploration. Since $e$ is
pivotal, the state of $e$ must become relevant before the first visit to
$S_L^{\Hh}(o)$: either an ordinary exploration uses $e$, or $e$
belongs to $E\left(B_r^{\Hh}(c)\right)$ for an enhancement centre $c$ that is actually used.
In the latter case the ordinary closure preceding the enhancement reaches every
vertex of $B_r^{\Hh}(c)$ and hence reaches both endpoints of $e$. Thus at least one
endpoint $x$ of $e$ is reached before $S_L^{\Hh}(o)$ is reached. By \eqref{eq:event-measurable}, no vertex outside $B_{L-1}^{\Hh}(o)$
can be reached before $S_L^{\Hh}(o)$ is reached. Hence we may choose $x$ so that
\begin{equation}\label{eq:root-reached-endpoint}
  x\in B_{L-1}^{\Hh}(o).
\end{equation}

We now choose the centre $z$ as follows. Suppose first $d_{\Hh}(o,e)>r$.  Then $d_{\Hh}(o,x)>r$.  Put
$d:=d_{\Hh}(o,x)$ and $t:=\max\{0,d+r-(L-1)\}.$ Since $r+1\le d\le L-1$ and $L\ge2r+2$, one has $0\le t\le r$.
Choose $z$ on a geodesic from $x$ to $o$ at distance $t$ from $x$.
Then
\begin{equation}\label{eq:root-z-ball}
  x\in B_r^{\Hh}(z)\subseteq B_{L-1}^{\Hh}(o),
  \
  o\notin B_r^{\Hh}(z).
\end{equation}
{\it Indeed}, $d_{\Hh}(o,z)=d-t>r$, while for every $q\in B_r^{\Hh}(z)$, $d_{\Hh}(o,q)\le d_{\Hh}(o,z)+r=d-t+r\le L-1.$

Suppose next $d_{\Hh}(o,e)\le r$.  Choose an endpoint $z$ of $e$ satisfying
$d_{\Hh}(o,z)\le r$.  Then
\begin{equation}\label{eq:root-near-z}
  o\in B_r^{\Hh}(z),
  \
  B_r^{\Hh}(z)\cap S_L^{\Hh}(o)=\varnothing,
\end{equation}
because $B_r^{\Hh}(z)\subseteq B_{2r}^{\Hh}(o)$ and $L\ge2r+2$.

In either case $d_{\Hh}(z,e)\le r$. Since $B_{2r+2}^{\Hh}(z)\subseteq B_R^{\Hh}(e)$ for $R=3r+2$,
\begin{equation}\label{eq:root-local-marks-erased}
  \alpha_0\cap B_{2r+2}^{\Hh}(z)=\varnothing.
\end{equation}

Set
\[
  \widetilde\omega:=\omega^{e=1}\setminus E\left(B_{r+1}^{\Hh}(z)\right);\ \text{that is, all edges with both endpoints in $B_{r+1}^{\Hh}(z)$ are forced closed.}
\]
Consider first the case $d_{\Hh}(o,e)>r$. The successful exploration in
$\left(\omega^{e=1},\alpha_0\right)$ reaches $x\in B_r^{\Hh}(z)$, whereas the root $o$ is outside
$B_r^{\Hh}(z)$.  We claim that
$$\text{some vertex $u\in S_{r+1}^{\Hh}(z)$ is still reached by the
enhanced exploration in $\left(\widetilde\omega,\alpha_0\right)$.}$$
{\it Actually}, if the root already
lies on $S_{r+1}^{\Hh}(z)$, take $u=o$.  Otherwise let $n$ be the first alternating
step of the successful exploration at which $S_{r+1}^{\Hh}(z)$ is reached.  Before
that time an enhancement also cannot jump directly into $B_r^{\Hh}(z)$: if an
enhancement centred at $c$ inserted a vertex $v\in B_r^{\Hh}(z)$, then
\[
  d_{\Hh}(c,z)\le d_{\Hh}(c,v)+d_{\Hh}(v,z)\le (r+1)+r=2r+1,
\]
so $c\in B_{2r+2}^{\Hh}(z)$, contradicting
\eqref{eq:root-local-marks-erased}. The first entrance into
$S_{r+1}^{\Hh}(z)$ cannot be caused by an enhancement either: if an enhancement
centred at $c$ first inserted a vertex of $S_{r+1}^{\Hh}(z)$, then
\[
  d_{\Hh}(c,z)\le (r+1)+(r+1)=2r+2,
\]
again contradicting \eqref{eq:root-local-marks-erased}. Hence the first
entrance is an ordinary one. Before that entrance no vertex of $B_{r+1}^{\Hh}(z)$ has been
reached, so a first open path to $S^{\Hh}_{r+1}(z)$ may be chosen whose final edge
enters $S_{r+1}^{\Hh}(z)$ from outside $B_{r+1}^{\Hh}(z)$. That path uses no edge of
$E\left(B_{r+1}^{\Hh}(z)\right)$ and therefore remains open in $\widetilde\omega$.

Choose such a $u$ according to the fixed order and choose a neighbour
$v\in S_r^{\Hh}(z)$ of $u$. Define
\begin{equation}\label{eq:root-surgery-far}
  \omega'
  :=\widetilde\omega\cup E\left(B_r^{\Hh}(z)\right)\cup\{\{u,v\}\}.
\end{equation}
In the case $d_{\Hh}(o,e)\le r$, no attaching edge is needed, and we define
\begin{equation}\label{eq:root-surgery-near}
  \omega':=\widetilde\omega\cup E\left(B_r^{\Hh}(z)\right).
\end{equation}
In both cases set $\alpha':=\alpha_0$.
\smallskip

We first show that $\cE_L$ fails when $z$ is inactive.  

\smallskip
\noindent\emph{Claim 1 (inactive containment).}
\begin{equation}\label{eq:root-inactive-containment}
 C^*:=C_o\left(\omega',\alpha'\right)
  \subseteq C_o\left(\widetilde\omega,\alpha_0\right)\cup B_r^{\Hh}(z).
\end{equation}

\noindent\emph{Proof of Claim 1.}
Let $\left(D'_n\right)_{n\ge0}$ be the alternating explored sets for
$\left(\omega',\alpha'\right)$. To prove
\begin{equation}\label{eq:root-inactive-stage}
  D'_n\subseteq C^*\cup B_r^{\Hh}(z),\ n\ge 0.
\end{equation}
The assertion is immediate at $n=0$ because $o\in C^*$. Suppose it holds at
an even stage $2m$, just before an ordinary closure, and let $y\notin B_r^{\Hh}(z)$
be joined to $D'_{2m}$ by an $\omega'$-open path. Outside
$E\left(B_{r+1}^{\Hh}(z)\right)$ the configurations $\omega'$ and $\widetilde\omega$ agree.
Every edge from $B_r^{\Hh}(z)$ to $S_{r+1}^{\Hh}(z)$ is closed in $\omega'$ except, in the
far case ($d_{\Hh}(o,e)>r$), the attaching edge $\{u,v\}$; moreover $u\in C^*$. Therefore each
excursion of the path through $B_r^{\Hh}(z)$ can only return through the same
attaching edge, or can leave the ball through $u$. After deleting such local
excursions, the portion leading to $y$ is a $\widetilde\omega$-open path
starting from a vertex of $C^*$. Since $C^*$ is closed under ordinary
$\widetilde\omega$-connectivity, $y\in C^*$. Hence
$D'_{2m+1}\subseteq C^*\cup B_r^{\Hh}(z)$.

At the following enhancement step, a centre in $B_{2r+2}^{\Hh}(z)$ has mark zero by
\eqref{eq:root-local-marks-erased}. If a primed reached centre $c$ lies
outside $B_{2r+2}^{\Hh}(z),$ then $c\in C^*$ by the odd-stage inclusion, and every edge
of $E\left(B_r^{\Hh}(c)\right)$ has the same state in $\omega'$ and
$\widetilde\omega$: otherwise $B_r^{\Hh}(c)\cap B_{r+1}^{\Hh}(z)\not=\varnothing$ and $d_{\Hh}(c,z)\le 2r+1$, a contradiction! Thus
any enhancement performed there by the primed process is also available from
$c$ in the exploration defining $C^*$, so all vertices inserted outside
$B_r^{\Hh}(z)$ already belong to $C^*$.  This proves the next even-stage inclusion.
\hfill$\triangleleft$\\

Since $e$ is incident to a vertex of $B_r^{\Hh}(z)$, one has
$e\in E\left(B_{r+1}^{\Hh}(z)\right)$ and hence $\widetilde\omega\le\omega^{e=0}.$
By \eqref{eq:root-e-remains-pivotal}, the exploration in
$(\widetilde\omega,\alpha_0)$ does not meet $S_L^{\Hh}(o)$. The ball $B_r^{\Hh}(z)$ is
also disjoint from $S_L^{\Hh}(o)$ by \eqref{eq:root-z-ball} in the far case ($d_{\Hh}(o,e)>r$) and by
\eqref{eq:root-near-z} in the near case ($d_{\Hh}(0,e)\leq r$). Thus
\begin{equation}\label{eq:root-z-off}
  \cE_L\left(\omega',\alpha'^{z=0}\right)=0.
\end{equation}

\smallskip

We next justify the resumption of the successful exploration stage by stage.

\smallskip
\noindent\emph{Claim 2 (active resumption).}
\begin{equation}\label{eq:root-resumption}
  C_o\left(\omega^{e=1},\alpha_0\right)
  \subseteq \widehat C:= C_o\left(\omega',\alpha'^{z=1}\right).
\end{equation}

\noindent\emph{Proof of Claim 2.}
Let $(F_n)_{n\ge0}$ be the alternating explored sets for
$(\omega^{e=1},\alpha_0)$. Firstly, we claim $B_{r+1}^{\Hh}(z)\subseteq\widehat C$. {\it In fact}, in the far case ($d_{\Hh}(o,e)>r$) the reached vertex $u$ is joined through
the open attaching edge $\{u,v\}$ to the completely open ball $B_r^{\Hh}(z)$; in the
near case ($d_{\Hh}(o,e)>r$) the root $o$ already lies in that ball. Hence the exploration reaches
$z$. Since all edges of $B_r^{\Hh}(z)$ are open, activating $z$ inserts
$S_{r+1}^{\Hh}(z)$, and therefore every vertex of $B_{r+1}^{\Hh}(z)$ is reached.

Secondly, we prove by induction on $n$ that $F_n\subseteq\widehat C$. The case $n=0$ is clear.  If
$F_{2m}\subseteq\widehat C$ and $y\in F_{2m+1}$, choose an
$\omega^{e=1}$-open path from $F_{2m}$ to $y$. Whenever this path runs inside
$B_{r+1}^{\Hh}(z)$, all of those vertices are already in $\widehat C$; every edge of
the path outside $E\left(B_{r+1}^{\Hh}(z)\right)$ has the same state in $\omega'$ as in
$\omega^{e=1}$.  Decomposing the path at its visits to $B_{r+1}^{\Hh}(z)$ therefore
shows $y\in\widehat C$.

If $F_{2m+1}\subseteq\widehat C$ and an enhancement of the original process is
performed at a centre $c\in F_{2m+1}$, then $c\notin B_{2r+2}^{\Hh}(z)$ because all
marks there vanish in $\alpha_0$.  For such $c$, the state of every edge in
$E\left(B_r^{\Hh}(c)\right)$ is unchanged by the surgery, and its mark is the same in
$\alpha_0$ and $\alpha'^{z=1}$. Thus the same enhancement is available in the
primed process and all vertices it inserts belong to $\widehat C$. Hence
$F_{2m+2}\subseteq\widehat C$. \hfill$\triangleleft$\\

By \eqref{eq:root-e-remains-pivotal}, $C_o\left(\omega^{e=1},\alpha_0\right)\cap S_L^{\Hh}(o)\not=\emptyset;$ hence $\cE_L\left(\omega',\alpha'^{z=1}\right)=1.$
Together with \eqref{eq:root-z-off}, this proves that $z$ is enhancement-pivotal.

Apart from the possible normalization of the coordinate $e$, all edge changes
occur in $E\left(B_{r+1}^{\Hh}(z)\right)$ and all erased enhancement coordinates lie in
$B_R^{\Hh}(e)$.  Since $d_{\Hh}(z,e)\le r$, every modified coordinate is contained in
$B_R^{\Hh}(e)$. Bounded degree gives a uniform bound, depending only on $r$ and the
degree bound, on the number of coordinates in this neighbourhood; enlarging
that bound by one covers the normalization step.  All choices were made
deterministically as required.
\end{proof}

\begin{figure}[H]
\centering
\begin{tikzpicture}[scale=0.86, every node/.style={font=\small}]
  \coordinate (z) at (0,0);
  \draw (z) circle (1.15);
  \draw[dashed] (z) circle (1.72);
  \draw[densely dotted] (z) circle (2.55);
  \draw[dash dot] (z) circle (3.20);
  \fill (z) circle (2pt) node[below left] {$z$};

  \node at (-0.10,0.88) {$B_r^{\Hh}(z)$};
  \node[anchor=west] at (1.78,-0.75) {$B_{r+1}^{\Hh}(z)$};
  \node[anchor=east] at (-2.62,1.40) {$B_{2r+2}^{\Hh}(z)$};
  \node[anchor=east] at (-3.28,2.10) {$B_R^{\Hh}(e)$};
  \node[anchor=west] at (2.62,-1.62) {local marks erased};

  \coordinate (v) at (1.10,0.22);
  \coordinate (u) at (1.70,0.34);
  \fill (v) circle (1.8pt) node[below] {$v$};
  \fill (u) circle (1.8pt) node[above] {$u$};
  \draw[very thick] (u)--(v) node[midway,below=2pt] {$\{u,v\}$};

  \coordinate (e1) at (0.35,-1.28);
  \coordinate (e2) at (0.98,-1.18);
  \draw[very thick] (e1)--(e2);
  \fill (e1) circle (1.4pt);
  \fill (e2) circle (1.4pt);
  \node[below=2pt] at (0.67,-1.28) {$e$};

  \coordinate (o) at (4.10,1.95);
  \fill (o) circle (2pt) node[above right] {$o$};
  \draw[thick] (o) .. controls (3.45,1.65) and (2.55,0.80) .. (u);
  \node[above,rotate=-23] at (2.80,1.35) {ordinary entrance};
\end{tikzpicture}
\caption{Schematic of the ``far'' surgery in Lemma~\ref{lem:surgery}.  The
edge $\{u,v\}$ attaches the previously reached exterior part to the fully open
ball $B_r^{\Hh}(z)$. When $z$ is inactive, the surgery blocks passage through the
ball; when $z$ is active, the enhancement reaches $S_{r+1}^{\Hh}(z)$ and the
successful exploration can be resumed. 
}
\label{fig:root-surgery}
\end{figure}
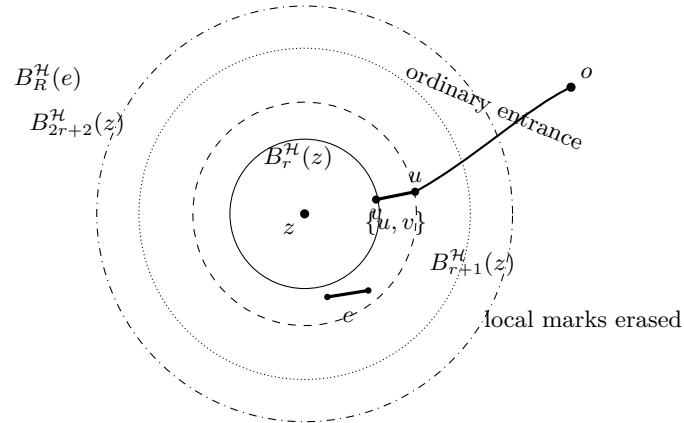

\renewcommand\thesection{\appendixname~C}
\setcounter{section}{0}
\section{Proof of Lemma \ref{lem:finite-set-pivotal-conversion}}\label{app-3}
\setcounter{equation}{0}
\renewcommand\thetheorem{C.\arabic{theorem}}
\renewcommand\theequation{C.\arabic{equation}}
\noindent
\begin{proof}
This lemma is \cite[Lemma 3.3]{Ma-Se2019}. For completeness, a detailed proof is given here according to that of \cite[Lemma 3.3]{Ma-Se2019}. The proof is similar to that of Lemma~\ref{lem:surgery} in \ref{app-2}.
Since all coordinates outside $B_L^{\Hh}(o)$ are forced to zero in
\eqref{eq:finite-volume-event}, {\it we suppress the truncation superscript during
the proof.}

Note edge-pivotality does not prescribe the actual state of $e$. It means precisely that
\begin{equation}\label{eq:e-pivotal-finite-set}
 \cE_L^{A,B}\left(\omega^{e=1},\alpha\right)=1,
 \
 \cE_L^{A,B}\left(\omega^{e=0},\alpha\right)=0.
\end{equation}
We therefore normalize every pivotal input by first replacing $\omega$ with
$\omega^{e=1}$. The normalization lies inside
$B_R^{\Hh}(e)$ and is included in the uniform local-coordinate count below, so the
surgery is deterministic for both possible actual states of the pivotal edge.
Erase the enhancement marks in $B_R^{\Hh}(e)$ one at a time, following the fixed
vertex order. If the event first fails when a mark at $z$ is erased, then $z$
is enhancement-pivotal and we are done.  Otherwise put $\alpha_0:=\alpha\setminus B_R^{\Hh}(e);$
then
\begin{equation}\label{eq:e-remains-pivotal-finite-set}
 \cE_L^{A,B}\left(\omega^{e=1},\alpha_0\right)=1,
 \ \mbox{and by monotonicity}\ \cE_L^{A,B}\left(\omega^{e=0},\alpha_0\right)=0.
\end{equation}

We first note that, in the successful exploration associated with
$\left(\omega^{e=1},\alpha_0\right)$, at least one endpoint of $e$ is reached before the
target set $B$ is reached.  {\it Indeed}, this follows from the following facts:
\begin{enumerate}
\item[{\bf (a)}] At least one endpoint of $e$ is relevant to
the exploration; otherwise changing the state of $e$ will not affect an ordinary
closure.
\item[{\bf (b)}]  Edge $e$ cannot affect an enhancement that is actually used. If conversely a
reached centre $c$ satisfies $e\in E\left(B_r^{\Hh}(c)\right)$ and $E\left(B_r^{\Hh}(c)\right)\subseteq\omega$,
then the ordinary closure preceding that enhancement reaches every vertex of
$B_r^{\Hh}(c)$, particularly both endpoints of $e$; and thus the state of $e$ is not pivotal, a contradiction.
\item[{\bf (c)}]  If a reached endpoint belongs to $B$, then the other endpoint
is also reached through the open edge $e$ before the event can depend on the
state of $e$.
\end{enumerate}
Therefore we choose, using the fixed order, a reached endpoint
$x\in e$ with $x\notin B;$ and write $e=\{x,y\}.$

\smallskip
We now distinguish the following two cases.

\smallskip
\noindent\emph{Case 1: $d_{\Hh}(e,A)>r$.}
Then $d_{\Hh}(x,A)>r$. We can choose a vertex $z$ such that
\begin{equation}\label{eq:z-far-properties}
 x\in B_r^{\Hh}(z),\  B_r^{\Hh}(z)\subseteq B_L^{\Hh}(o)\setminus(A\cup B).
\end{equation}
We verify this carefully, including the boundary case. If $d_{\Hh}(x,B)>r$ and $d_{\Hh}\left(x,S_L^{\Hh}(o)\right)\ge r$, take $z=x$. Then
$B_r^{\Hh}(x)\cap(A\cup B)=\varnothing$ because both $d_{\Hh}(x,A)$ and $d_{\Hh}(x,B)$ exceed
$r$. Also $B_r^{\Hh}(x)\subseteq B_L^{\Hh}(o)$: if a path of length at most $r$ from $x$
leaves $B_L^{\Hh}(o)$, it will meet $S_L^{\Hh}(o)$ before leaving, producing a point of
$S_L^{\Hh}(o)$ at distance strictly less than $r$ from $x$.

Suppose next $d_{\Hh}(x,B)\le r$. By $r$-niceness of $B,$ there is a vertex
$z$ with $x\in B_r^{\Hh}(z)$ and $B_r^{\Hh}(z)\cap B=\varnothing.$ Choose $b\in B$ with $d_{\Hh}(x,b)\le r.$
If there is $a\in A\cap B_r^{\Hh}(z),$ then
$d_{\Hh}(a,b)\le d_{\Hh}(a,z)+d_{\Hh}(z,x)+d_{\Hh}(x,b)\le 3r,$
contrary to $d_{\Hh}(A,B)>3r$. Thus $B_r^{\Hh}(z)\cap A=\varnothing$. If there is
$q\in B_r^{\Hh}(z)\cap S_L^{\Hh}(o)$, then
\[
 d_{\Hh}\left(b,S_L^{\Hh}(o)\right)\le d_{\Hh}(b,x)+d_{\Hh}(x,z)+d_{\Hh}(z,q)\le 3r,\ \text{contrary to $d_{\Hh}\left(B,S_L^{\Hh}(o)\right)>3r$}.
\]
Hence $B_r^{\Hh}(z)\cap S_L^{\Hh}(o)=\varnothing$. Since $B_r^{\Hh}(z)$
contains the already reached vertex $x\in B_L^{\Hh}(o)$ and is connected, it cannot
contain a vertex outside $B_L^{\Hh}(o)$ without crossing $S_L^{\Hh}(o)$. So $B_r^{\Hh}(z)\subseteq B_L^{\Hh}(o)$ and \eqref{eq:z-far-properties} holds.

It remains to consider $d_{\Hh}\left(x,S_L^{\Hh}(o)\right)<r$. This subcase automatically has
$d_{\Hh}(x,B)>r$: otherwise
$$d_{\Hh}\left(B,S_L^{\Hh}(o)\right)\le d_{\Hh}(B,x)+d_{\Hh}\left(x,S_L^{\Hh}(o)\right)<2r,\ \text{contradicting
\eqref{eq:finite-volume-separation}.} $$
Let $d_0=d_{\Hh}(o,x)\le L$ and put $t:=\max\{0,d_0+r-L\}\in[0,r].$ Choose $z$ on a geodesic from $x$ to $o$ at distance $t$ from $x$. Then
$d_{\Hh}(x,z)=t\le r$, so $x\in B_r^{\Hh}(z)$, and $d_{\Hh}(o,z)=d_0-t\le L-r.$
Consequently every $q\in B_r^{\Hh}(z)$ satisfies $d_{\Hh}(o,q)\le L$, and $B_r^{\Hh}(z)\subseteq B_L^{\Hh}(o)$. Moreover, for every $D\in\{A,B\},$ the separation
condition and $d_{\Hh}\left(x,S_L^{\Hh}(o)\right)<r$ imply
$ d_{\Hh}(x,D)\ge d_{\Hh}\left(D,S_L^{\Hh}(o)\right)-d_{\Hh}\left(x,S_L^{\Hh}(o)\right)>3r-r=2r.$
If there is $q\in D\cap B_r^{\Hh}(z)$, then
$d_{\Hh}(x,q)\le d_{\Hh}(x,z)+d_{\Hh}(z,q)\le 2r$, a contradiction. Thus $B_r^{\Hh}(z)\cap(A\cup B)=\varnothing.$
\vskip 2mm

\noindent\emph{Case 2: $d_{\Hh}(e,A)\le r$.}
Choose an endpoint $x_A$ of $e$ satisfying $d_{\Hh}(x_A,A)\le r$ and put $z=x_A$.
There is $a\in A$ with $d_{\Hh}(a,z)\le r$, so $A\cap B_r^{\Hh}(z)\ne\varnothing$. If there is
$b\in B\cap B_r^{\Hh}(z)$, then
\[
 d_{\Hh}(A,B)\le d_{\Hh}(a,b)\le d_{\Hh}(a,z)+d_{\Hh}(z,b)\le 2r,\ \text{contrary to $d_{\Hh}(A,B)>3r$};
\]
hence $B_r^{H}(z)\cap B=\varnothing$. If there is $q\in B_r^{\Hh}(z)\cap S_L^{\Hh}(o)$, then
$d_{\Hh}\left(A,S_L^{\Hh}(o)\right)\le d_{\Hh}(a,q)\le 2r,$
contrary to \eqref{eq:finite-volume-separation}. The ball $B_r^{\Hh}(z)$ contains
$a\in A\subseteq B_L^{\Hh}(o)$ and is connected, so the absence of an intersection
with $S_L^{\Hh}(o)$ implies $B_r^{\Hh}(z)\subseteq B_L^{\Hh}(o)$. Thus all local modifications
remain inside the truncated volume.
\smallskip

In either case $d_{\Hh}(z,e)\le r$, and since $B_{2r+2}^{\Hh}(z)\subseteq B_R^{\Hh}(e)$ for $R=3r+2$,
\begin{equation}\label{eq:local-marks-erased}
 \alpha_0\cap B_{2r+2}^{\Hh}(z)=\varnothing.
\end{equation}

Set $\widetilde\omega:=\omega^{e=1}\setminus E\left(B_{r+1}^{\Hh}(z)\right);$
that is, all edges with both endpoints in $B_{r+1}^{\Hh}(z)$ are forced closed.  In
Case~1 we claim that
$$\text{there is a vertex $u\in S_{r+1}^{\Hh}(z)$ which is reached
by the enhanced exploration from $A$ in $\left(\widetilde\omega,\alpha_0\right)$.}$$
{\it In fact}, the successful exploration in $\left(\omega^{e=1},\alpha_0\right)$ reaches $x\in B_r^{\Hh}(z)$, so
it must first reach $S_{r+1}^{\Hh}(z)$.  If the starting set already meets this
sphere, choose the first such vertex.  Otherwise let $n$ be the first
alternating step for which the explored set meets $S_{r+1}^{\Hh}(z)$.  Before that
time an enhancement cannot jump directly into $B_r^{\Hh}(z)$: if an enhancement
centred at $c$ inserted $v\in B_r^{\Hh}(z)$, then
$d_{\Hh}(c,z)\le d_{\Hh}(c,v)+d_{\Hh}(v,z)\le 2r+1,$
so \eqref{eq:local-marks-erased} forces its mark to be zero.  By the same erased-mark condition, the first
entrance into $S_{r+1}^{\Hh}(z)$ cannot be an enhancement either: a centre $c$
capable of inserting a vertex of $S_{r+1}^{\Hh}(z)$ satisfies $d_{\Hh}(c,z)\le 2r+2$.
Thus the first entrance occurs during an ordinary closure.  Before that entrance no vertex of
$B_{r+1}^{\Hh}(z)$ has been reached, and a first ordinary path to the sphere can be
chosen so that its final edge enters $S_{r+1}^{\Hh}(z)$ from outside
$B_{r+1}^{\Hh}(z)$. That path uses no edge of $E\left(B_{r+1}^{\Hh}(z)\right)$ and therefore remains
open in $\widetilde\omega$.

Choose such a $u$ according to the fixed order, and choose a neighbour
$v\in S_r^{\Hh}(z)$ of $u$.  In Case~1 define
\begin{equation}\label{eq:finite-set-surgery-far}
 \omega'
 :=\widetilde\omega\cup E\left(B_r^{\Hh}(z)\right)\cup\{\{u,v\}\}.
\end{equation}
In Case~2 no attaching edge is needed; define
\begin{equation}\label{eq:finite-set-surgery-near}
 \omega':=\widetilde\omega\cup E\left(B_r^{\Hh}(z)\right).
\end{equation}
In both cases set $\alpha':=\alpha_0$.
\smallskip

We first show the event fails when $z$ is inactive.
\smallskip

\noindent\emph{Claim 1 (inactive containment).}
\begin{equation}\label{eq:inactive-cluster-identity}
 C_A(\omega',\alpha')
 \subseteq C_A(\widetilde\omega,\alpha_0)\cup B_r^{\Hh}(z).
\end{equation}

\noindent\emph{Proof of Claim 1.}
Put $C^*:=C_A(\widetilde\omega,\alpha_0)$ and let $(D'_n)_{n\ge0}$ be the alternating explored sets started from $A$ for
$(\omega',\alpha')$.  We prove by induction that $D'_n\subseteq C^*\cup B_r^{\Hh}(z),\ n\ge 0.$

At stage zero this is clear: $A\subseteq C^*$ (note that $A\cap B_r^{\Hh}(z)=\varnothing$ in Case~1,
while possibly $A\cap B_r^{\Hh}(z)\not=\varnothing$ in Case~2). Suppose
the inclusion holds before an ordinary closure and let a primed open path end
at $y\notin B_r^{\Hh}(z)$.  Outside $E\left(B_{r+1}^{\Hh}(z)\right)$, $\omega'$ and
$\widetilde\omega$ agree.  In Case~1 the only reopened edge from $B_r^{\Hh}(z)$ to
$S_{r+1}^{\Hh}(z)$ is $\{u,v\}$ and $u\in C^*$; in Case~2 no such edge is reopened.
Deleting excursions inside $B_r^{\Hh}(z)$ therefore leaves a
$\widetilde\omega$-open path from $C^*$ to $y$.  Hence $y\in C^*$, proving the
odd-stage inclusion.

At an enhancement step, any centre in $B_{2r+2}^{\Hh}(z)$ has mark zero by
\eqref{eq:local-marks-erased}.  A reached centre outside that ball belongs to
$C^*$ by the preceding inclusion, and its entire $r$-ball has unchanged edge
states; an affected $r$-ball would imply $d_{\Hh}(c,z)\le 2r+1$. Thus any primed enhancement outside $B_r^{\Hh}(z)$ is also available in the exploration
defining $C^*$, and its inserted vertices already lie in $C^*$. The next even-stage inclusion holds. \hfill$\triangleleft$
\smallskip

Now $e\in E\left(B_{r+1}^{\Hh}(z)\right)$, since $x\in B_r^{\Hh}(z)$ and $e$ is incident to $x$.
Hence $\widetilde\omega\le\omega^{e=0}$.  By
\eqref{eq:e-remains-pivotal-finite-set}, the exploration in
$\left(\widetilde\omega,\alpha_0\right)$ does not reach $B$.  Since
$B_r^{\Hh}(z)\cap B=\varnothing$, \eqref{eq:inactive-cluster-identity} gives
\begin{equation}\label{eq:z-off-fails}
 \cE_L^{A,B}\left(\omega',\alpha'^{z=0}\right)=0.
\end{equation}

We next activate $z$. In Case~1 the reached vertex $u$ is joined through the
open edge $\{u,v\}$ to the open ball $B_r^{\Hh}(z)$; in Case~2 the
starting set already meets that ball. Thus the exploration reaches $z$, and
since any edge of $B_r^{\Hh}(z)$ is open, the enhancement at $z$ inserts
$S_{r+1}^{\Hh}(z)$. So any vertex of $B_{r+1}^{\Hh}(z)$ is reached.
\smallskip

\noindent\emph{Claim 2 (active resumption).}
\begin{equation}\label{eq:finite-set-resumption}
 C_A\left(\omega^{e=1},\alpha_0\right)
 \subseteq \widehat C:=C_A\left(\omega',\alpha'^{z=1}\right).
\end{equation}

\noindent\emph{Proof of Claim 2.}
Let $(F_n)_{n\ge0}$ be the alternating explored sets started from $A$ for
$\left(\omega^{e=1},\alpha_0\right).$ We have already shown $B_{r+1}^{\Hh}(z)\subseteq\widehat C$.
To prove by induction that $F_n\subseteq\widehat C$ for any $n$.  At an ordinary step,
take an $\omega^{e=1}$-open path from $F_{2m}$ to a vertex of $F_{2m+1}$.
All vertices of the path lying in $B_{r+1}^{\Hh}(z)$ are already in $\widehat C$,
and any edge outside $E\left(B_{r+1}^{\Hh}(z)\right)$ has the same state in $\omega'$ and
$\omega^{e=1}$.  Splitting the path at its visits to $B_{r+1}^{\Hh}(z)$ therefore
shows that its endpoint lies in $\widehat C$.

At an enhancement step, any centre used by the original process has nonzero
mark in $\alpha_0$ and therefore lies outside $B_{2r+2}^{\Hh}(z)$.  Its $r$-ball
edge states are unchanged by the surgery and its mark is unchanged in the
primed configuration.  Hence any enhancement available in the original
process is also available from the already reached centre in $\widehat C$.
This gives $F_{2m+2}\subseteq\widehat C$ and completes the induction.
\hfill$\triangleleft$
\smallskip

By \eqref{eq:e-remains-pivotal-finite-set}, $C_A\left(\omega^{e=1},\alpha_0\right)\cap B\not=\varnothing;$ hence
$\cE_L^{A,B}\left(\omega',\alpha'^{z=1}\right)=1.$ Together with \eqref{eq:z-off-fails}, this proves that $z$ is
enhancement-pivotal. The rest of the proof is the same as that of Lemma~\ref{lem:surgery}.
\end{proof}

\end{document}